\documentclass[12pt]{elsart}
\usepackage{amsfonts,amsmath}

\journal{
}
\date{}
\begin{document}
\begin{frontmatter}
\title
{$q-$Pascal's triangle and irreducible representations of the
braid group $B_3$ in arbitrary dimension}

\thanks[dec]{The second author would like to thank the Max-Planck-Institute
of Mathematics, Bonn and the Institute of Applied Mathematics,
University of Bonn for the hospitality. The partial financial
support by the DFG project 436 UKR 113/87 is gratefully
acknowledged.}
\author{Sergio Albeverio}
\address{
Institut f\"ur Angewandte Mathematik, Universit\"at Bonn,
Wegelerstr. 6, D-53115 Bonn, Germany; SFB 611, IZKS; Bonn, BiBoS,
Bielefeld--Bonn, Germany; CERFIM, Locarno; Accademia di
Architettura, USI, Mendrisio, Switzerland \\
E-mail: albeverio@uni.bonn.de}

\author{Alexandre Kosyak\corauthref{cor}}
\ead{kosyak01@yahoo.com, kosyak@imath.kiev.ua}
\address{Max-Planck-Institut f\"ur Mathematik, Vivatsgasse 7, D-53111 Bonn, Germany}
\address{Institute of Mathematics, Ukrainian National Academy of Sciences,
3 Tereshchenkivs'ka, Kyiv, 01601, Ukraine,\\
E-mail: kosyak01@yahoo.com, kosyak@imath.kiev.ua \\
tel.: 38044 2346153 (office), 38044 5656758 (home), fax: 38044
2352010 }

\corauth[cor]{Corresponding author}

\newpage

\begin{abstract}
We construct a $\left[\frac{n+1}{2}\right]+1$ parameters family of
irreducible representations of the Braid group $B_3$ in arbitrary
dimension $n\in {\mathbb N}$, using a $q-$deformation of the
Pascal triangle. This construction extends in particular results
by S.P.~Humphries [8], who constructed  representations of the
braid group $B_3$ in arbitrary dimension using the classical
Pascal triangle. E.~Ferrand [7] obtained  an equivalent
representation of $B_3$ by considering two special operators in
the space ${\mathbb C}^n[X].$ Slightly more general
representations were given by I.~Tuba and H.~Wenzl [11]. They
involve  $[\frac{n+1}{2}]$ parameters (and also use the classical
Pascal triangle). The latter authors also gave  the complete
classification of all simple representations of $B_3$ for
dimension $n\leq 5$.  Our construction generalize all mentioned
results and throws a new light on  some of them. We also study the
irreducibility and the equivalence of the representations.

 In \cite{Kos07q} we establish the connection between the constructed
representation of  the braid group $B_3$ and  the highest weight
modules of $U({\mathfrak sl}_{2})$ and quantum group
$U_q({\mathfrak sl}_{2})$.
\end{abstract}

\begin{keyword}
Braid group, ${\rm SL }(2,{\mathbb Z})$, representations,
classification, Pascal's triangle, $q-$Pascal's triangle,
$q-$binomial coefficient, Gaussian polynomials, quantum groups

\MSC 20F36 \sep (05A30, 11B56, 17B37)
\end{keyword}
\end{frontmatter}
\sloppy
\newcommand{\tr}{\mathrm{tr}\,}
\newcommand{\rank}{\mathrm{rank}\,}
\newcommand{\diag}[1]{\mathrm{diag}\,(#1)}
\renewcommand{\Im}{\mathrm{Im}\,}
\maketitle
\newpage
\tableofcontents
\section{Introduction} Let $B_3$ be Artin's braid group, given
by the generators $\sigma_1$ and $\sigma_2$ and the relation
$\sigma_1\sigma_2\sigma_1=\sigma_2\sigma_1\sigma_2$ \cite{Art26}.
Let ${\mathbb C}$ be the field of complex numbers, ${\mathbb
N}:=\{0,1,2,...\}$ and let ${\rm Mat}(n,{\mathbb C})$ be the set
of complex $n\times n$ matrices. We construct irreducible
representations of the braid group $B_3$ in the space ${\mathbb
C}^{n+1}$ for an arbitrary $n\in {\mathbb N}$ using the $q-$Pascal
triangle, i.e. the $q-$deformation of the usual Pascal triangle.

Let $\left(\begin{smallmatrix} n\\
k
\end{smallmatrix}\right)_q,\,\,k,n\in{\mathbb
N}$ be the $q-$binomial coefficients (see the definition in
Section 2). For the matrix $A=(a_{km})_{0\leq k,m,\leq n}\in {\rm
Mat}(n+1,{\mathbb C})$ we set $A^\sharp=(a^\sharp_{ij})$ where
$a^\sharp_{ij}=a_{n-i,n-j}$. In all sections except Section 9 we
index row and columns of the matrix $A\in {\rm Mat}(n+1,{\mathbb
C})$ starting from $0$.
For any  nonzero $q\in {\mathbb C}$ let the  matrix $\Lambda_n(q)$
and the numbers $q_n$ be defined as follows:
\begin{equation}
\label{La(q)}
 \Lambda_n(q)=
{\rm diag}\left(q_{rn} \right)_{r=0}^n,\text{\,\,where\,\,}
q_{rn}:=\frac{q_rq_{n-r}}{q_n}=q^{-(n-r)r}\text{\,\,and\,\,}
q_n=q^{\frac{(n-1)n}{2}},\,\,n\in{\mathbb N}.
\end{equation}
For an arbitrary  $n\in {\mathbb N}$ and a complex diagonal matrix
$\Lambda={\rm diag}(\lambda_0,\lambda_1,...,\lambda_n)$ we define
our {\bf representation of the group $B_3$ in the space ${\mathbb
C}^{n+1}$} by  following formulas
\begin{equation}
\label{Rep(q)}
\sigma_1\mapsto\sigma_1^{\Lambda}:=\sigma_1(q,n)\Lambda\quad\text{\,and\,}\,\,
\sigma_2\mapsto\sigma_2^{\Lambda}:=\Lambda^\sharp\sigma_2(q,n)=
\Lambda^\sharp(\sigma_1^{-1}(q^{-1},n))^\sharp,
\end{equation}
where $\sigma_1(q,n)=(\sigma_1(q,n)_{km})_{0\leq k,m\leq n}$ and
$\sigma_2(q,n)$ are defined by
\begin{equation}
\label{si_1(q)} \sigma_1(q,n)_{km}=
\sigma_1(q)_{km}=\left(\begin{smallmatrix} n-k\\
n-m
\end{smallmatrix}\right)_q,\,\,0\leq k,m\leq n,\,\,\,
\sigma_2(q,n)=(\sigma_1^{-1}(q^{-1},n))^\sharp,
\end{equation}
(as usually, we set $\left(\begin{smallmatrix} n\\
k
\end{smallmatrix}\right)_q=0$ for $k>n$), and the matrix $\Lambda=\Lambda_n$ satisfies the
following condition:
\begin{equation}
\label{cond_q} \lambda_0\lambda_n
\Lambda_n(q)=\Lambda_n\Lambda_n^\sharp \text{\quad or\quad}
\lambda_0\lambda_n\frac{q_rq_{n-r}}{q_n}=\lambda_r\lambda_{n-r},\,\,0\leq
r\leq n.
\end{equation}

The aim of this article is to show that formulas (\ref{Rep(q)})
give an $\left[\frac{n+1}{2}\right]+1$ parameters family of $B_3$
representations in dimension $n+1$, for any $n\in{\mathbb N}$
(Theorem 1) and study the irreducibility (Theorem 3,4) and the
equivalence (Theorem  5).

In Section 2 we introduce the main objects and give the main
statements. In Section 3 we present the result of S.P.~Humphries
\cite{Hum00} who has constructed a representation of $B_3$
equivalent with the particular case of our representation when
$q=1$ and $\Lambda=I$
\begin{equation}
\label{Hum1} \sigma_1\mapsto\sigma_1(1,n),\quad
\sigma_2\mapsto\sigma_2(1,n).
\end{equation}
In his representations  the classical Pascal  triangle plays a
basic role. In Section 4 we mention the result of E.~Ferrand
\cite{Fer05} who has constructed a representation of $B_3$
equivalent with the representation (\ref{Hum1}). For this
E.~Ferrand has considered two operators $\Phi:p(x)\mapsto p(x+1)$
and $\Psi:p(x)\mapsto (1-x)^np(\frac{x}{1-x})$ in the space
${\mathbb C}^n[X]$ of the polynomials of degree $\leq n$
satisfying the relation $\Phi\Psi\Phi=\Psi\Phi\Psi$. The
$q-$analogue of the mentioned results is given in Section 13.

In Section 5 we show that this representation is closely connected
with the morphism  $\rho:B_3\mapsto {\rm SL}(2,{\mathbb Z})$ (see
(\ref{TW1}) below) and the $n$th symmetric power of the natural
representation $\pi: {\rm SL}(2,{\mathbb Z})\mapsto {\rm
SL}(2,{\mathbb Z})$. Section 6 recalls shortly the results of
I.~Tuba and H.~Wenzl \cite{TubWen01}. Firstly, they showed that
\begin{equation}
\label{TW0}  \sigma_1\mapsto\sigma_1(1,n)\Lambda,\quad
\sigma_2\mapsto\Lambda^\sharp\sigma_2(1,n)
\end{equation}
is a representation of $B_3$ in arbitrary dimension $n+1$ if
$\lambda_r\lambda_{n-r}=c$ for some constant $c$ (see Remark 6.2,
Section 2). Secondly, they gave the complete classification of all
irreducible representations of $B_3$ for dimension $\leq 5$.

Our motivation for the present study was to generalize the results
and formulas of the mentioned authors to the case of an arbitrary
dimension $n\in{\mathbb N}$. We have realized that not only the
classical Pascal triangle may be used for constructing of the
representations of $B_3$ but also  the $q-$deformed Pascal
triangle. The conditions $\lambda_r\lambda_{n-r}=c$ on $\Lambda$
in the classical case should be replaced in the deformed case by
some rather nontrivial conditions (see (\ref{cond_q})) connecting
the matrix $\Lambda$ with some canonical diagonal matrix
$\Lambda_n(q)$, depending on $q$ and $n$. We prove that the
representations of $B_3$ given by (\ref{Rep(q)}) coincide with the
representations of I.~Tuba and H.~Wenzl \cite{TubWen01} for $n=4$,
are equivalent with them for $n=2,3,5$, and generalize them for an
arbitrary dimension $n$ (Remark 6.3, Section 2). This is explained
in Section 10.

In Section 7 we present the results of S.P.~Humphries and the
results of I.~Tuba and H.~Wenzl in a form which is convenient for
our extensions. In Section 8 we show how the Pascal (resp.
$q-$Pascal ) triangle appears as the operators $\exp T_1$ (resp.
$\exp_{(q)} T_{(q)}$) associated with some operators $T_1$ (resp.
$T_{(q)}$). The irreducibility and the equivalence of our
representations (Theorem 3, 4 and 5) are studied in Section 9.

In Section 11 we give the proof of the Theorem 1 i.e. that
(\ref{Rep(q)}) is a representation. In Section 12 we prove some
combinatorial identities for $q-$binomial coefficients. These
identities are an essential part in the proof of the Theorem 1.
They generalize  the well-known combinatorial identities for
classical binomial coefficients (see \cite{GasRah}) used by
S.P.~Humphries to prove that (\ref{Hum1}) is a representation of
$B_3$.

Let us also mention that in  the article of S.~Albeverio and
S.~Rabanovich \cite{AlbRab07} a class of unitary irreducible
representations of $B_3$ by $n\times n$ matrices for every $n\geq
3$ was constructed. Using tensor products of these representations
and the reduced Barrau  representations \cite{Jon87} these authors
also find a class of irreducible unitary representations of $B_4$.

In \cite{ForWSV03} E.Formanek et al. gave the {\it complete
classification} of all {\it simple representations of $B_n$} for
{\it dimension} $\leq n$. To know more on the braid groups and its
applications see \cite{Bir93,BirBre93}.

\section{Main objects}
The  Pascal triangle consists of the
binomial coefficients $\left(\begin{smallmatrix} n\\
k
\end{smallmatrix}\right):=C^k_n,\,\,
k,n\in{\mathbb N},$ defined by
\begin{equation}
\label{n!} C_n^k:=\frac{n!}{k!(n-k)!},\,\,0\leq k\leq n,
C_n^k=0,\,\,k>n ,\,\,\text {where}\,\,n!=1\cdot2\cdot...\cdot n.
\end{equation}
We recall that the binomial coefficients may also be defined  by
induction, using the relations
$$
C_n^0=C_n^n=1,\,\,n\in{\mathbb N},\quad
C_{n+1}^k=C_n^{k-1}+C_n^k,\quad 1\leq k\leq n.
$$

 We also consider the $q-$Pascal triangle consisting of
$q-$binomial coefficients $
\left(\begin{smallmatrix} n\\
k
\end{smallmatrix}\right)_q
=C^k_n(q),\,\,0\leq k\leq
n,\,\,n\in{\mathbb N}$ defined as follows (see
\cite{And76,KacChen01,Kas95,KliSch97} )
\begin{equation}
\label{C_n^k(q)}
\left(\begin{smallmatrix} n\\
k
\end{smallmatrix}\right)_q:=C_n^k(q)=\frac{(q;q)_n}{(q;q)_k(q;q)_{n-k}},\,\,0\leq
k\leq n,\quad C_n^k(q)=0,\,\, k>n,
\end{equation}
where  $(a;q)_n$ denotes the standard $q-$shifted factorial
\cite{And76,AndAsk85,KliSch97}
\begin{equation}
\label{(a;q)_n} (a;q)_n=(1-a)(1-aq)(1-aq^2)\ldots(1-aq^{n-1}).
\end{equation}
The $q-$binomial coefficients were first studied by Gauss and have
come to be known as {\it Gaussian polynomials} (see \cite[Ch. 3.3
]{And76}). Another definition \cite[Ch.IV.2]{Kas95} of the Gauss
polynomials is following. For any integer $n>0$ we define the
associated $q-${\it integer} $(n)_q$ by
\begin{equation}
\label{(n)}
(n)_q=1+q+\cdots+q^{n-1}=\frac{1-q^n}{1-q}.
\end{equation}
Define the $q-${\it factorial} of $n$ by $(0)!_q=1$ and
\begin{equation}
\label{(n)!}
(n)!_q=(1)_q(2)_q\dots(n)_q=\frac{(1-q)(1-q^2)\dots(1-q^n)}{(1-q)^n},
\end{equation}
when $n>0$. We define the Gaussian polynomials for $0\leq k\leq n$
by
\begin{equation}
\label{GP} C_n^k(q)=\frac{(n)!_q}{(k)!_q(n-k)!_q}.
\end{equation}
 They
can also be defined recursively, using the following  relations
(\cite[Ch.IV.2]{Kas95}, see also \cite{KacChen01})
$C_n^0(q)=C_n^n(q)=1,\,\,n\in{\mathbb N}:$
\begin{equation}
\label{GP1} C_{n+1}^k(q)=C_n^{k-1}(q)+q^kC_n^k(q),\,\,
C_{n+1}^k(q)=q^{n-k}C_n^{k-1}(q)+C_n^k(q),\,\, 1\leq k\leq n.
\end{equation}
We can also obtain the Gaussian polynomials in the following way.
Let us set
$$
(1+x)^k_q:=(1+x)(1+xq)(1+xq^2)\ldots(1+xq^{n-1}).
$$  We have (see
\cite{And76})
\begin{equation}\label{GP2}
(1+x)^k_q=\sum_{r=0}^kq^{r(r-1)/2}C_k^r(q)x^r
=\sum_{r=0}^kq^{r(r-1)/2}
\left(\begin{smallmatrix} k\\
r
\end{smallmatrix}\right)_qx^r.
\end{equation}
As an example we make explicit the corresponding $q-$Pascal
triangle for  $n=5$:
$$
\begin{smallmatrix}
 &&&&&1&&&&&\\
&&&&1&&1&&&&\\
&&&1&&1+q&&1&&&\\
&&1&&1+q+q^2&&
C^2_3(q)&&1&&\\
&1&&(1+q)(1+q^2)&&(1+q^2)(1+q+q^2)&&C^3_4(q)&&1&\\
1&&1+q+q^2+q^3+q^4&&
(1+q^2)(1+q+q^2+q^3+q^4)&&
C^3_5(q)&&C^4_5(q)&&1\\
\end{smallmatrix}
$$

{\bf Notations}. For an $n\times n$ matrix $A=(a_{ij})$ we set
$A^t$ (resp. $A^s$ and $A^\sharp$) where
\begin{equation}
\label{sharp} A^t=(a^t_{ij}),\, a^t_{ij}=a_{ji},\text{\,(resp\,}\,
A^s=(a^s_{ij}),\,\,a^s_{ij}=a_{n-j,n-i};\,\,
A^\sharp=(a^\sharp_{ij}),\,\, a^\sharp_{ij}=a_{n-i,n-j}).
\end{equation}
The operation $A\to A^\sharp$ means composing the transposition
with respect to the main diagonal ($A\to A^t$)  with the
transposition with respect to the auxiliary (subsidiary) diagonal
($A\to A^s$) i.e. $A^\sharp=(A^t)^s=(A^s)^t$.

Let us consider the $(n+1)\times(n+1)$ matrix $S(q)$
defined as follows:
\begin{equation}
\label{S(q)} S(q)=(S(q)_{km}),\quad\text{\,\,where\,\,}
S(q)_{km}=q^{-1}_k(-1)^k\delta_{k+m,n},\,S:=S(1),
\end{equation}
$$
S(q)=\left(\begin{smallmatrix}
0&...&0&0&0&1\\
0&...&0&0&-1&0\\
0&...&0&q^{-1}&0&0\\
0&...&-q^{-3}&0&0&0\\
&...&&&\\
(-1)^{n}q^{-\frac{(n-1)n}{2}}&...&0&0&0&0
\end{smallmatrix}\right)\text{\,\,\,and\,\,\,}
S=\left(\begin{smallmatrix}
0&...&0&0&0&1\\
0&...&0&0&-1&0\\
0&...&0&1&0&0\\
0&...&-1&0&0&0\\
&...&&&&\\
(-1)^{n}&...&0&0&0&0
\end{smallmatrix}\right).
$$
\begin{thm}
\label{braid.t.1} Formulas (\ref{Rep(q)}) define a representation
of $B_3$ in the space ${\mathbb C}^{n+1}$ for an arbitrary $n\in
{\mathbb N}$ i.e.
\begin{equation}
\label{Br_n(q)}
\sigma_1^{\Lambda}\sigma_2^{\Lambda}\sigma_1^{\Lambda}=
\sigma_2^{\Lambda}\sigma_1^{\Lambda}\sigma_2^{\Lambda}=\lambda_0\lambda_nS(q)\Lambda,
\end{equation}
moreover
\begin{equation}
\label{si_2(q)} (\sigma_1^{-1}(q^{-1}))^\sharp_{km}=
\left\{\begin{array}{cc}
0,&\text{\,if\,}\,\,0< k< m\leq n,\\
 (-1)^{k+m}q_{k-m}^{-1}C_{k}^{m}(q^{-1}),&\text{\,if\,}\,\,0\leq m\leq k\leq n.
\end{array}\right.
\end{equation}
\end{thm}
\begin{rem}
1. Let us set
\begin{equation}
\label{D_n(q)} D_n(q)={\rm diag}(q_r)_{r=0}^n.
\end{equation}
 We have by (\ref{La(q)}) and (\ref{S(q)})
\begin{equation}
\label{L=D_n(q)} \Lambda(q)=q_n^{-1}D_n(q)D_n^\sharp(q),\quad
S(q)=D_n^{-1}(q)S,
\end{equation}
so if we take $\Lambda=D_n(q)$ or $\Lambda=D_n^\sharp(q)$ the
relation (\ref{cond_q}) is satisfied, hence
\begin{equation}
\label{RepD(q)}
\sigma_1\mapsto\sigma_1^D(q,n):=\sigma_1(q,n)D_n^\sharp(q)
\quad\text{\,and\,}\quad \sigma_2\mapsto
\sigma_2^D(q,n):=D_n(q)\sigma_2(q,n)
\end{equation}
also gives a representation of the braid group $B_3.$\\
2. The general form of the matrix $\Lambda_n$ satisfying
(\ref{cond_q}) is following:  $\Lambda_n=D_n^\sharp(q)\Lambda_n'$
or $\Lambda_n=D_n(q)\Lambda_n'$ where $\Lambda_n'={\rm
diag}(\lambda_0',\lambda_1',...,\lambda_n')$ with
$\Lambda_n'(\Lambda_n')^\sharp=cI$ for some constant $c$.
\end{rem}
Using Remark 2 we shall write our representation in the following
form
\begin{equation}
\label{Rep(q)'}
\sigma_1\mapsto\sigma_1^{\Lambda}(q,n)=\sigma_1(q,n)D_n^\sharp(q)\Lambda_n,\,\,
\sigma_2\mapsto\sigma_2^{\Lambda}(q,n)=\Lambda_n^\sharp
D_n(q)\sigma_2(q,n),\,\, \Lambda_n\Lambda_n^\sharp=cI.
\end{equation}

{\bf Definition}. {\it We say that the representation is {\bf
subspace irreducible} or {\bf ireducible} (resp. {\bf operator
irreducible}) when there no nontrivial invariant close {\bf
subspaces} for all operators of the representation (resp.  there
no nontrivial bounded {\bf operators} commuting with all operators
of the representation).}

Let us define for $n,r,q,\lambda$ such that $n\in{\mathbb
N},\,\,0\leq r\leq n,\,\,\lambda\in{\mathbb
C}^{n+1},\,\,q\in{\mathbb C}$ the following operators
\begin{equation}
\label{F(nu)}
 F_{r,n}(q,\lambda)=\exp_{(q)}\left(\sum_{k=0}^{n-1}
(k+1)_qE_{kk+1}\right)-
q_{n-r}\lambda_r(D_n(q)\Lambda_n^\sharp)^{-1},
\end{equation}
where $\exp_{(q)} X=\sum_{m=0}^\infty X^m/(m)!_q$.
For the matrix $C\in{\rm Mat}(n+1,{\mathbb C})$ we
denote by
$$
M^{i_1i_2...i_r}_{j_1j_2...j_r}(C),\,\,({\rm resp.\,\,}
A^{i_1i_2...i_r}_{j_1j_2...j_r}(C)),\,\, 0\leq i_1<...<i_r\leq
n,\,\, 0\leq j_1<...<j_r\leq n
$$
its  minors (resp. the cofactors) with $i_1,i_2,...,i_r$ rows and
$j_1,j_2,...,j_r$ columns.
\begin{thm}
\label{t.IrrB_3} The representation of the group $B_3$ defined by
(\ref{Rep(q)'})
have the following properties:\\
1) for $q=1,\,\,\Lambda_n=1$, it is  {\rm subspace irreducible} in arbitrary dimension $n\in{\mathbb N }$;\\
2) for $q\not=1,\,\,\Lambda_n={\rm diag}(\lambda_k)_{k=0}^n\not=1$
it is {\rm operator irreducible} if and only if for any $0\leq
r\leq \left[\frac{n}{2}\right]$ there exists $0\leq
i_0<i_i<...<i_r\leq n $ such that
\begin{equation}
\label{M^i_j...=0} M^{i_0i_i...i_{n-r-1}}_{r+1r+2...n} (
F_{r,n}^s(q,\lambda))\not=0;
\end{equation}
3) for $q\not=1,\,\,\Lambda_n=1$ it is  {\it subspace irreducible}
if and only if $(n)_q\not=0$.\\
 The representation has $[\frac{n+1}{2}]+1$ free parameters.
\end{thm}
Let us denote by $\sigma^\Lambda(q,n)$ the representation of $B_3$
defined by (\ref{Rep(q)'}).
\begin{thm}
\label{SIrr} The representation $\sigma^\Lambda(q,n)$ is {\rm
subspace irreducible} for $n=1$ if and only if
$\Lambda_1\not=\lambda_0(1,\alpha)$ where $\alpha^2-\alpha+1=0$.

\end{thm}
{\bf Problem.} {\it To find a criteria of the subspace
irreducibility for all representations $\sigma^\Lambda(q,n)$. Some
particular cases are studied in Section 8}.
\begin{thm}
\label{EQ} If two represenations $\sigma^\Lambda(q,n)$ and
$\sigma^{\Lambda'}(q',n)$  are equivalent i.e.
$$
\sigma_i^\Lambda(q,n)C=C\sigma_i^{\Lambda'}(q',n),\,\,i=1,2
$$
for some $C\in {\rm GL}(n+1,{\mathbb C})$ then $q/q'=1$ for $n=2m$
and $(q/q')^2=1$ for $n=2m-1.$
\end{thm}
\begin{rem} 1. In the particular case where $\Lambda=I$ and $q=1$ Theorem 1
gives the result of S.P.~Humphries \cite{Hum00} (see Section 3).\\
2. When $q=1$ and $\Lambda={\rm
diag}(\lambda_0,\lambda_1,...,\lambda_n)$ we obtain the example of
I.~Tuba and H.~Wenzl \cite{TubWen01} (see Section 5, Example 1).\\
3. The representations of $B_3$ given by (\ref{Rep(q)'}) coincide
with the representations of I.~Tuba and H.~Wenzl \cite{TubWen01}
for $n=4$, are equivalent with them in the dimension $n=2,3,5$,
and generalize them for an arbitrary dimension $n$.\\
4. Using Theorem 3 and 4 we give in  Section 9.5 examples of
representations of $B_3$ that are operator irreducible but are
subspace reducible.
\end{rem}
\begin{thm}In particular using result of \cite{TubWen01} (Sections
2.4-2.7) we conclude that all irreducible representations of $B_3$
for dimension $\leq 5$ are given by (\ref{Rep(q)'}).
\end{thm}
\section{Pascal's triangle and representations of $B_3$.
Results of Humphries } Following  S.P.~Humphries \cite{Hum00}, for
fixed $n\geq 1$ we let $\Sigma_1=\Sigma_1(n)$ and
$\Sigma_2=\Sigma_2(n)$ (respectively) be the following
$(n+1)\times(n+1)$ lower and upper (respectively) triangular
matrices:
$$
\left(\begin{array}{ccccccc}
1&&&&&&\\
1&1&&&&&\\
1&2&1&&&&\\
1&3&3&1&&&\\
1&4&6&4&1&&\\
&&&&&...& \\
 \left(\begin{smallmatrix}
 n\\
 0
\end{smallmatrix}\right)&
 \left(\begin{smallmatrix}
 n\\
 1
\end{smallmatrix}\right)&
 \left(\begin{smallmatrix}
 n\\
 2
\end{smallmatrix}\right)&
 \left(\begin{smallmatrix}
 n\\
 3
\end{smallmatrix}\right)&
 \left(\begin{smallmatrix}
 n\\
 4
\end{smallmatrix}\right)&...&
 \left(\begin{smallmatrix}
 n\\
 n
\end{smallmatrix}\right)
\end{array}\right),
\quad \left(\begin{array}{ccccccc}
 \left(\begin{smallmatrix}
 n\\
 n
\end{smallmatrix}\right)&...&
 \left(\begin{smallmatrix}
 n\\
 4
\end{smallmatrix}\right)&
 \left(\begin{smallmatrix}
 n\\
 3
\end{smallmatrix}\right)&
 \left(\begin{smallmatrix}
 n\\
 2
\end{smallmatrix}\right)&
 \left(\begin{smallmatrix}
 n\\
 1
\end{smallmatrix}\right)
&
 \left(\begin{smallmatrix}
 n\\
 0
\end{smallmatrix}\right)\\
&&&&&...& \\
&&1&4&6&4&1\\
&&&1&3&3&1\\
&&&&1&2&1\\
&&&&&1&1\\
&&&&&&1\\
\end{array}\right),
$$
(Thus we make the convention that a blank indicates the zero
entry). Let $E=E_n$ be the $(n+1)\times (n+1)$ permutation matrix
corresponding to the permutation $(0 n)(1 n-1)(2 n-2)...$. S.P.~
Humphries shows that
\begin{equation}\label{Hum0}
\sigma_1\mapsto\Sigma_1,\quad
\sigma_2\mapsto\Sigma_2^{-1}
\end{equation}
 gives a representation of $B_3$ using  the following lemmas.

{\bf Lemma 4.1} {\it We have $E\Sigma_1E^{-1}=\Sigma_2$. Further
$$
\Sigma_2^{-1}=\left(\begin{array}{rrrrrrr}
 \left(\begin{smallmatrix}
 n\\
 n
\end{smallmatrix}\right)&...&(-1)^{n-4}
 \left(\begin{smallmatrix}
 n\\
 4
\end{smallmatrix}\right)&(-1)^{n-3}
 \left(\begin{smallmatrix}
 n\\
 3
\end{smallmatrix}\right)&(-1)^{n-2}
 \left(\begin{smallmatrix}
 n\\
 2
\end{smallmatrix}\right)&(-1)^{n-1}
 \left(\begin{smallmatrix}
 n\\
 1
\end{smallmatrix}\right)&
 (-1)^{n}\left(\begin{smallmatrix}
 n\\
 0
\end{smallmatrix}\right)\\
&...&&&&& \\
&&1&-4&6&-4&1\\
&&&1&-3&3&-1\\
&&&&1&-2&1\\
&&&&&1&-1\\
&&&&&&1\\
\end{array}\right)
$$
There is a similar expression for $\Sigma_1^{-1}$, namely
$\Sigma_1^{-1}=E^{-1}\Sigma_2^{-1}E$.}

{\bf Lemma 4.2} {\it We have }
$$
\Sigma_1\Sigma_2^{-1}=\left(\begin{array}{rrrrrrr}
                         1&
 -\left(\begin{smallmatrix}
 n\\
 1
\end{smallmatrix}\right)&
 \left(\begin{smallmatrix}
 n\\
 2
\end{smallmatrix}\right)&
- \left(\begin{smallmatrix}
 n\\
 3
\end{smallmatrix}\right)&
 \left(\begin{smallmatrix}
 n\\
 4
\end{smallmatrix}\right)&...&1\\
&&&&&...& \\
1&-4&6&-4&1&&\\
1&-3&3&-1&&&\\
1&-2&1&&&&\\
1&-1&&&&&\\
1&&&&&&
\end{array}\right)
$$

{\bf Lemma 4.3} {\it We have
$\Sigma_1\Sigma_2^{-1}\Sigma_1=\Sigma_2^{-1}\Sigma_1\Sigma_2^{-1}=
(-1)^nG_n$, where $G_n$ is the $(n+1)\times (n+1)$ matrix ${\rm
diag}(1,-1,1,...)E_n$.}

\section{Pascal's triangle in the space ${\mathbb C}^n[X]$ and  results of E.~Ferrand}
In the work of E.~Ferrand \cite{Fer05} the Pascal triangle appears
in the following way. Denote by $\Phi$ the endomorphism of the
space ${\mathbb C}^n[X]$ of polynomials of degree $n$ with complex
coefficients, which maps a polynomial $p(X)$ to the polynomial
$p(X+1)$. Denote by $\Psi$ the endomorphism of ${\mathbb C}^n[X]$
which maps a a polynomial $p(X)$ to
$(1-X)^np\left(\frac{X}{1-X}\right)$

{\bf Theorem} \cite{Fer05}. {\it $\Phi$ and $\Psi$ verify a
braid-like relation $\Phi\Psi\Phi=\Psi\Phi\Psi$.}

One can verify that  $\Phi$ (resp. $\Psi$) in the canonical basis
$1,X,X^2,...,X^n$ of ${\mathbb C}^n[X]$ have
 the form
\begin{equation}
\label{Fer1} \Phi=\Sigma_2(n)^s=\sigma_1^s(1,n),\quad
\Psi=(\Sigma_1(n)^{-1})^s=(\sigma_2^{-1}(1,n))^s,
\end{equation}
where the notations $A^s$ is defined in (\ref{sharp}).
For the operator $\Phi$ we have
\begin{equation}
\label{Fer2} X^k\stackrel{\Phi}{\mapsto}
 (1+X)^k=\sum_{r=0}^kC_k^rX^r
\end{equation}
hence $\Phi_{rk}=C_k^r$ and we get the first part of (\ref{Fer1})
if we compare (\ref{Fer2}) with (\ref{si_1(q)}). For the operator
$\Psi$ we get
\begin{equation}
\label{Fer3}
X^k\stackrel{\Psi}{\mapsto}(1-X)^{n-k}X^k=\sum_{r=0}^{n-k}(-1)^rC_{n-k}^rX^{r+k}=
\sum_{t=k}^{n}(-1)^{k+t}C_{n-k}^{t-k}X^t
\end{equation}
 if we set $r+k=t$. Hence $\Psi_{tk}=(-1)^{k+t}C_{n-k}^{t-k}$ and
we get the second part of (\ref{Fer1}) if we compare (\ref{Fer3})
with (\ref{si_2(q)}). Since $\Phi\Psi\Phi =\Psi\Phi\Psi$ we have
another proof of the braid relations given in \cite{Hum00}: $
\Sigma_1\Sigma_2^{-1}\Sigma_1=\Sigma_2^{-1}\Sigma_1\Sigma_2^{-1}.
$

\section{Pascal's triangle as the symmetric power}
The representation of $B_3$ given by E.~Ferrand  can be obtained
in the following way. There is a morphism $\rho:B_3\mapsto {\rm SL
}(2,{\mathbb Z})$ of the group $B_3$ in ${\rm SL }(2,{\mathbb Z})$
defined by (\ref{TW1}) below. Let us consider the natural
representation $\pi$ of the group ${\rm SL }(2,{\mathbb Z})$ in
the space ${\mathbb C}^1[X]$ defined as follows
$$
(\pi_gf)(x)=(cx+d)f\left(\frac{ax+b}{cx+d}\right),\text{\,\,
where\,\,} g=\left(\begin{smallmatrix}
a&b\\
c&d
\end{smallmatrix}\right)\in {\rm SL }(2,{\mathbb Z}).
$$
We show that (see (\ref{si_1(q)}) for the notation
$\sigma_1(1,n)$)
\begin{equation}
\label{Fer4} {\rm
Sym}^n(\pi)\circ\rho(\sigma_k)=\sigma_k(1,n),\,\,k=1,2,\,\,n\in{\mathbb
N},
\end{equation}
where ${\rm Sym}^n(\pi)$ is the symmetric power of the
representation $\pi.$   We have
$$
\rho(\sigma_1)=\sigma_1(1,1)=\left(\begin{smallmatrix}
1&1\\
0&1
\end{smallmatrix}\right),\quad
\rho(\sigma_2)=\sigma_2(1,1)=\left(\begin{smallmatrix}
1&0\\
-1&1
\end{smallmatrix}\right).
$$
Then (\ref{Fer4}) is transformed into
\begin{equation}\label{Fer5}
 {\rm Sym}^n(\sigma_1(1,1))=\sigma_1(1,n).
\end{equation}
Let us take the basis $e_0,\,e_1$ of the space $V:={\mathbb
C}^1[X]\simeq {\mathbb C}^2$. In the space $V\otimes V$ with the
basis $e_{km}:=e_k\otimes e_m$ ordered as follows $
e_{00},\,\,e_{01},\,\,e_{10},\,\,e_{11},$ we have (see, e.g.,
\cite[Ch. 2]{Kas95} for the definition of the tensor product of
two operators)
$$
\left(\begin{smallmatrix}
1&1\\
0&1
\end{smallmatrix}\right)\otimes\left(\begin{smallmatrix}
1&1\\
0&1
\end{smallmatrix}\right)=
\left(\begin{smallmatrix}
1&1&1&1\\
0&1&0&1\\
0&0&1&1\\
0&0&0&1
\end{smallmatrix}\right).
$$
The symmetric basis in the space ${\rm Sym}^2(V)\subset V\otimes
V$ is as follows
\begin{equation}
\label{sym.b.2} e_{0}^s=e_{00}=e_0\otimes
e_0,\,\,e_{1}^s=e_{01}+e_{10}=e_0\otimes e_1+e_1\otimes
e_0,\,\,e_{2}^s=e_{11}=e_1\otimes e_1.
\end{equation}
The symmetric basis in the space ${\rm Sym}^n(V)\subset
\underbrace{V\otimes...\otimes V}$ for $n\in{\mathbb N}$ is
\begin{equation}
\label{sym.b.n}e_k^s=\frac{1}{k!(n-k)!}\sum_{\sigma\in
S_{n+1}}\sigma(e_k),\,\,0\leq k\leq n, \text{\,\,
where\,\,}e_k=e_0\otimes...\otimes e_0\otimes \underbrace{
e_1\otimes...\otimes e_1},
\end{equation}
and $\sigma(e_{i_0}\otimes e_{i_1}\otimes...\otimes e_{i_n})
=(e_{\sigma(i_0)}\otimes e_{\sigma(i_1)}\otimes...\otimes
e_{\sigma(i_n)})$ for $\sigma\in S_{n+1}.$ Since
$\left(\begin{smallmatrix}
1&1\\
0&1
\end{smallmatrix}\right)e_0=e_0$ and $\left(\begin{smallmatrix}
1&1\\
0&1
\end{smallmatrix}\right)e_0=e_0+e_1$ we have for the operator
$\left(\begin{smallmatrix}
1&1\\
0&1
\end{smallmatrix}\right)\otimes \left(\begin{smallmatrix}
1&1\\
0&1
\end{smallmatrix}\right)$ in the symmetric basis:
$$
\left(\begin{smallmatrix}
1&1\\
0&1
\end{smallmatrix}\right)\otimes \left(\begin{smallmatrix}
1&1\\
0&1
\end{smallmatrix}\right)(
e_0\otimes e_1)=e_0\otimes e_1=e_{0}^s,\quad
\left(\begin{smallmatrix}
1&1\\
0&1
\end{smallmatrix}\right)\otimes \left(\begin{smallmatrix}
1&1\\
0&1
\end{smallmatrix}\right)(
e_0\otimes e_1+e_1\otimes e_0)
$$
$$
=e_0\otimes (e_0+e_1)+(e_0+e_1)\otimes e_0=2(e_0\otimes
e_1)+e_0\otimes e_0+e_1\otimes e_0=2e_0^s+e_1^s,
$$
$$
\left(\begin{smallmatrix}
1&1\\
0&1
\end{smallmatrix}\right)\otimes \left(\begin{smallmatrix}
1&1\\
0&1
\end{smallmatrix}\right)(e_1\otimes e_1)=
(e_0+e_1)\otimes(e_0+e_1)=e_{00}+e_{01}+e_{10}+e_{11}=e_{0}^s+e_{1}^s+e_{2}^s,
$$
hence the operator $\left(\begin{smallmatrix}
1&1\\
0&1
\end{smallmatrix}\right)\otimes \left(\begin{smallmatrix}
1&1\\
0&1
\end{smallmatrix}\right)$ in the symmetric basis has the form
$$
\left(\begin{smallmatrix}
1&2&1\\
0&1&1\\
0&0&1\\
\end{smallmatrix}\right)=\sigma_1(1,2).
$$
The proof of the relations (\ref{Fer5}) for general $n\in{\mathbb
N}$ is similar, indeed we have
$$
{\rm Sym}^n(\left(\begin{smallmatrix}
1&1\\
0&1
\end{smallmatrix}\right))e_k^s=\sum_{r=0}^k \left(\begin{smallmatrix}
 n-k\\
 n-r
\end{smallmatrix}\right) e_r^s.
$$
This proves (\ref{Fer5}).

\section{Results of I.~Tuba and H.~Wenzl}
Consider the braid group $B_3$ given by the generators $\sigma_1$
and $\sigma_2$ and the relation
$\sigma_1\sigma_2\sigma_1=\sigma_2\sigma_1\sigma_2$.  $B_3$ maps
surjectively onto ${\rm SL}(2,{\mathbb Z})$ via the map $\rho$
given by
\begin{equation}
\label{TW1}
 \sigma_1\mapsto\left(\begin{smallmatrix}
1&1\\
0&1
\end{smallmatrix}\right),\quad
\sigma_2\mapsto\left(\begin{smallmatrix}
1&0\\
-1&1
\end{smallmatrix}\right).
\end{equation}
It is easy to check that this is a homomorphism. Moreover we have
\begin{equation}
\label{Br2}
 \sigma_1\sigma_2\sigma_1=\sigma_2\sigma_1\sigma_2\mapsto
S= \left(\begin{smallmatrix}
0&1\\
-1&0
\end{smallmatrix}\right),\quad
\sigma_1\sigma_2\mapsto\left(\begin{smallmatrix}
0&1\\
-1&1
\end{smallmatrix}\right),\quad
\sigma_2\sigma_1\mapsto\left(\begin{smallmatrix}
1&1\\
-1&0
\end{smallmatrix}\right).
\end{equation}
{\bf Example 1} \cite{TubWen01}.  Let $V$ be a $(d+1)-$dimensional
vector space with a basis labeled by $0,1,...,d$, and let
$\lambda_0,\lambda_1,...,\lambda_d$ be parameters satisfying
$\lambda_i\lambda_{d-i}=c$ for a fixed constant $c$. Set
$\overline{i}=d-i.$ Then in \cite{TubWen01} it is shown that a
$\left[\frac{d+1}{2}\right]$ parameter family of representations
of $B_3$ is given by the matrices
\begin{equation}
\label{PasTW}
A=\left( \left(\begin{smallmatrix} \overline{i}\\
\overline{j}
\end{smallmatrix}\right)\lambda_i
\right)_{ij},\quad
B=\left((-1)^{i+j} \left(\begin{smallmatrix} i\\
j
\end{smallmatrix}\right)\lambda_{\overline{i}}
\right)_{ij}.
\end{equation}
The proof consists in checking that  $ABA=BAB=S$ with $S$ being
the skew-diagonal matrix defined by $s_{ij}=(
\delta_{i,\overline{j}}(-1)^i\lambda_{\overline{i}})$. This in
turn can be derived by the identity
$$
\sum_{k=0}^d(-1)^{i+k}\left(\begin{smallmatrix} i\\
k
\end{smallmatrix}\right)\left(\begin{smallmatrix} \overline{k}\\
\overline{j}
\end{smallmatrix}\right)=(-1)^i\left(\begin{smallmatrix} d-i\\
d-j
\end{smallmatrix}\right)=(-1)^i\left(\begin{smallmatrix} \overline{i}\\
\overline{j}
\end{smallmatrix}\right)
$$
(cf.  \cite[p.8 eq. (5)]{Rio68} see also  (\ref{Bin2}) below).
Another result in \cite{TubWen01} is as follows.

{\bf Proposition 2.5} . Let $V$ be a simple $B_3$ module of
dimension $n=2,3$. Then there exist a basis for $V$ for which
$\sigma_1$ and $\sigma_2$ act as follows ($\lambda=(\lambda_k)_k$)
\begin{equation}
\label{TW2} \sigma_1\mapsto\sigma_1^{
\lambda}:=\left(\begin{smallmatrix}
\lambda_1&\lambda_1\\
0&\lambda_2&
\end{smallmatrix}\right),\quad
\sigma_2\mapsto\sigma_2^{ \lambda}:=\left(\begin{smallmatrix}
\lambda_2&0\\
-\lambda_2&\lambda_1
\end{smallmatrix}\right)\text{\,\,for\,\,} n=2.
\end{equation}
\begin{equation}
\label{TW3}
\sigma_1\mapsto\sigma_1^{\lambda}=\left(\begin{smallmatrix}
\lambda_1&\lambda_1\lambda_3\lambda_2^{-1}+\lambda_2&\lambda_2\\
0&\lambda_2&\lambda_2\\
0&0&\lambda_3
\end{smallmatrix}\right),\quad
\sigma_2\mapsto\sigma_2^{ \lambda}:=\left(\begin{smallmatrix}
\lambda_3&0&0\\
-\lambda_2&\lambda_2&0\\
\lambda_2&-\lambda_1\lambda_3\lambda_2^{-1}-\lambda_2&\lambda_1
\end{smallmatrix}\right)\text{\,\,for\,\,} n=3.
\end{equation}
Let us set $D=\sqrt{\lambda_2\lambda_3/\lambda_1\lambda_4}$. All
simple modules for $n=4$ are following:
\begin{equation}
\label{TW41}
\sigma_1\mapsto\sigma_1^{\lambda}=\left(\begin{smallmatrix}
\lambda_1&(1+D^{-1}+D^{-2})\lambda_2&(1+D^{-1}+D^{-2})\lambda_3&\lambda_4\\
0&\lambda_2&(1+D^{-1})\lambda_3&\lambda_4\\
0&0&\lambda_3&\lambda_4\\
0&0&0&\lambda_4
\end{smallmatrix}\right),
\end{equation}
\begin{equation}
\label{TW42}
\sigma_2\mapsto\sigma_2^{\lambda}=\left(\begin{smallmatrix}
\lambda_4&0&0&0\\
-\lambda_3&\lambda_3&0&0\\
D\lambda_2&-(D+1)\lambda_2&\lambda_2&0\\
-D^3\lambda_1&(D^3+D^2+D)\lambda_1&-(D^2+D+1)\lambda_1&\lambda_1
\end{smallmatrix}\right).
\end{equation}
Let us set
$\gamma=(\lambda_1\lambda_2\lambda_3\lambda_4\lambda_5)^{1/5}$.
All simple modules for $n=5$ are following:
\begin{equation}
\label{TW5} \sigma_1\mapsto\sigma_1^{
\lambda}=\left(\begin{smallmatrix} \lambda_1 &
(1+\frac{\gamma^2}{\lambda_2\lambda_4})(\lambda_2+\frac{\gamma^3}{\lambda_3\lambda_4})
&
(\frac{\gamma^2}{\lambda_3}+\lambda_3+\gamma)(1+\frac{\lambda_1\lambda_5}{\gamma^2})&
(1+\frac{\lambda_2\lambda_4}{\gamma^2})
(\lambda_3+\frac{\gamma^3}{\lambda_2\lambda_4})&\frac{\gamma^3}{\lambda_1\lambda_5}
\\
0&\lambda_2&\frac{\gamma^2}{\lambda_3}+\lambda_3+\gamma&
\frac{\gamma^3}{\lambda_1\lambda_5}+\lambda_3+\gamma&\frac{\gamma^3}{\lambda_1\lambda_5}
\\
0&0&\lambda_3&\frac{\gamma^3}{\lambda_1\lambda_5}+\lambda_3
&\frac{\gamma^3}{\lambda_1\lambda_5}\\
0&0&0&\lambda_4&\lambda_4\\
0&0&0&0&\lambda_5
\end{smallmatrix}\right).
\end{equation}
The formula for $\sigma_2^{ \lambda}$ is not given in
\cite{TubWen01}. In section 9 we show that $\sigma_2^{
\lambda}=C^{-1}\sigma_2^{\Lambda}C$ where $C={\rm
diag}(1,1,1,q^{-1}\frac{\lambda_3}{\lambda_4},q^{-1}\frac{\lambda_3}{\lambda_5})$
and $\sigma_2^{\Lambda}=\Lambda^\sharp\sigma_2(q,4)$
(see(\ref{Rep(q)})--(\ref{si_1(q)})).
\section{Pascal's triangle}
We shall rewrite the results and the proof of S.P.~Humphries
\cite{Hum00} in a slightly different form, using also the results
of I.~Tuba and H.~Wenzl \cite{TubWen01} in order to generalize
Humphries'  result to the case of a $q-$Pascal triangle.
S.P.~Humphries uses the representation of $B_3$
$$\sigma_1\mapsto\Sigma_1,\quad \sigma_2\mapsto\Sigma_2^{-1}$$
I.~Tuba and H.~Wenzl use another representation of $B_3$ (in the
notations of S.P.~Humphries)
$$
\sigma_1\mapsto\Sigma_2,\quad \sigma_2\mapsto\Sigma_1^{-1}.
$$
Obviously these two representations are isomorphic and the
isomorphism is given by
$$
\sigma_1\mapsto\sigma_2^{-1}\text{\,\,and\,\,}\sigma_2\mapsto\sigma_1^{-1}.
$$
We shall use the form of representation given by I.~Tuba and
H.~Wenzl.

In the general case (for arbitrary $n\in{\mathbb N}$) we put,
using Pascal's triangle
\begin{equation}
\label{Rep(1)} \sigma_1\mapsto
\sigma_1(1):=\sigma_1(1,n):=\Sigma_2(n),\quad \sigma_2\mapsto
\sigma_2(1):=\sigma_2(1,n):=\Sigma_1^{-1}(n),
\end{equation}
where (see (\ref{Rep(q)}) and (\ref{si_1(q)}))
$\sigma_1(1)=(\sigma_1(1)_{km})_{0\leq k,m\leq n},\,\
\sigma_2(1)=(\sigma_2(1)_{km})_{0\leq k,m\leq n}$ and
\begin{equation}
\label{si_1(1)} \sigma_1(1)_{km}=\left\{\begin{array}{cc}
C_{n-k}^{n-m},&\text{\,if\,}\,\,
0\leq k\leq m\leq n,\\
0,&\text{\,if\,}\,\,0<m<k\leq n,
\end{array}\right.
\end{equation}
\begin{equation}
\label{si_2(1)} \sigma_2(1)_{km}= \left\{\begin{array}{cc}
0,&\text{\,if\,}\,\,0< k< m\leq n,\\
 (-1)^{k+m}C_{k}^{m},&\text{\,if\,}\,\,0\leq m\leq k\leq n.
\end{array}\right.
\end{equation}
\begin{thm}{\rm \cite{Hum00,TubWen01}} For $\sigma_1(1)$ and $\sigma_2(1)$ defined by
(\ref{si_1(1)}) and (\ref{si_2(1)}),  $\Lambda=I$ and arbitrary
$n\in{\mathbb N}$ we have
\begin{equation}
\label{Br_n(1)}
\sigma_1(1)\sigma_2(1)\sigma_1(1)=\sigma_2(1)\sigma_1(1)\sigma_2(1)=
S= \left(\begin{smallmatrix}
0&...&0&0&1\\
0&...&0&-1&0\\
0&...&1&0&0\\
&...&&&\\
(-1)^n&...&0&0&0
\end{smallmatrix}\right)
\end{equation}
moreover,  we have
\begin{equation}
\sigma_2(1)=(\sigma_1^{-1}(1))^\sharp.
\end{equation}
\end{thm}
\begin{pf} The identity (\ref{Br_n(1)}) is equivalent with
\begin{equation}
\label{Braid(1)}
\sigma_1(1)\sigma_2(1)=S\sigma_1^{-1}(1)=\sigma_2^{-1}(1)S.
\end{equation}
We have in particular
\begin{equation}
\label{trg(1)} \sigma_1(1)\sigma_2(1)=(\sigma_{km}^{12})_{0\leq
k,m\leq n},\text{\,where\,}
\sigma_{km}^{12}= \left\{\begin{array}{cl} 0,&\text{\,if\,\,}0\leq k+m<n,\\
(-1)^{n-m}C_k^{n-m},&\text{\,if\,\,}k+m\geq n,
\end{array}\right.
\end{equation}
and
\begin{equation}
\sigma_2(1)\sigma_1(1)=(\sigma_{km}^{21})_{0\leq k,m\leq
n},\text{\,where\,} \sigma_{km}^{21}= \left\{\begin{array}{cl}
(-1)^kC^m_{n-k},&\text{\,if\,\,}
0\leq k+m\leq n, \\
0,&\text{\,if\,\,}k+m> n.
\end{array}\right.
\end{equation}

We have
$\sigma_1(1)_{km}=C_{n-k}^{n-m}.$ To prove that $
\sigma_1^{-1}(1)_{km}=(-1)^{k+m}C_{n-k}^{n-m}$ we observe that
$$
(\sigma_1(1)\sigma_1^{-1}(1))_{km}=\sum_{r=k}^{n}\sigma_1(1)_{kr}\sigma_1^{-1}(1)_{rm}
=\sum_{r=k}^{n}C_{n-k}^{n-r} (-1)^{r+m}C_{n-r}^{n-m}
$$
$$
= \sum_{r=0}^n(-1)^{r+m}\left(\begin{smallmatrix}
 n-k\\
 n-r
\end{smallmatrix}\right)
\left(\begin{smallmatrix}
 n-r\\
 n-m
\end{smallmatrix}\right)=
 \sum_{r=0}^n(-1)^{(n-r)+(n-m)}\left(\begin{smallmatrix}
 n-k\\
 n-r
\end{smallmatrix}\right)
\left(\begin{smallmatrix}
 n-r\\
 n-m
\end{smallmatrix}\right)
= \delta_{km},
$$
(where in the latter step we have used the well-known identity
(\ref{Bin1}), Section 11 below). Analogously $
\sigma_2(1)_{km}=(-1)^{k+m}C_{k}^{m}$. To prove
$\sigma_2^{-1}(1)_{km}=C_{k}^{m}$ we observe that
$$
(\sigma_2(1)\sigma_2^{-1}(1))_{km}=\sum_{r=0}^{n}\sigma_2(1)_{kr}\sigma_2^{-1}(1)_{rm}
=\sum_{r=0}^{n}(-1)^{k+r}C_{k}^{r} C_{r}^{m}=
$$
$$
\sum_{i=0}^n(-1)^{k+r}\left(\begin{smallmatrix}
 k\\
 r
\end{smallmatrix}\right)\left(\begin{smallmatrix}
 r\\
 m
\end{smallmatrix}\right)=
\delta_{km},
$$
(using again (\ref{Bin1}) in the last step). Further the identity
(\ref{Braid(1)})
$$
 \sigma_1(1)\sigma_2(1)=S\sigma_1^{-1}(1),\quad
\sigma_1(1)\sigma_2(1)=\sigma_2^{-1}(1)S,
$$
means
$$
(\sigma_1(1)\sigma_2(1))_{km}=(S\sigma_1^{-1}(1))_{km}\text{\quad
and\quad} (\sigma_1(1)\sigma_2(1))_{km}=(\sigma_2^{-1}(1)S)_{km}.
$$
But
$$
(\sigma_1(1)\sigma_2(1))_{km}=\sum_{r=k}^{n}\sigma_1(1)_{kr}\sigma_2(1)_{rm}=
\sum_{r=k}^{n} C_{n-k}^{n-r}(-1)^{r+m}C_{r}^{m}.
$$
Since $S=(S_{km})$, where $S_{km}=(-1)^k\delta_{k+m,n}$ (see
(\ref{S(q)})), we get
$$
(S\sigma_1^{-1}(1))_{km}=S_{k,n-k}\sigma_1^{-1}(1)_{n-k,m}=
(-1)^k(-1)^{n-k+m}C^{n-m}_{n-(n-k)}
$$
$$
=(-1)^{n+m}C^{n-m}_{k},\,\,
(\sigma_2^{-1}(1)S)_{km}=\sigma_2^{-1}(1)_{k,n-m}S_{n-m,m}=C^{n-m}_k(-1)^{n-m},
$$
so
\begin{equation}
(S\sigma_1^{-1}(1))_{km}=(\sigma_2^{-1}(1)S)_{km}.
\end{equation}
 Finally the identity (\ref{Braid(1)}) is
equivalent with the following
$$
\sum_{r=k}^{n} C_{n-k}^{n-r}(-1)^{r+m}C_{r}^{m}
=(-1)^{n-m}C_k^{n-m} \text{\,\,or \,\,}
\sum_{r=0}^n(-1)^{n-r}\left(\begin{smallmatrix}
 n-k\\
 n-r
\end{smallmatrix}\right)
\left(\begin{smallmatrix}
 r\\
 m
\end{smallmatrix}\right)=
\left(\begin{smallmatrix}
 k\\
 n-m
\end{smallmatrix}\right),
$$
which is easily proven (in the latter step we have used
(\ref{Bin2}) below).
\qed\end{pf}

\section{Pascal's triangle as  $\exp T$ }
We give here some useful presentation for Pascal's (resp.
$q-$Pascal's) triangle as operators of the form $\exp T$ (resp.
$\exp_{(q)} T_q$) of some operators $T$ (resp. $T_q$). Let us
consider $\sigma_1(1,n)$ and $\sigma_1(1,n)^s.$ Since by
(\ref{si_1(q)}) we have
$\sigma_1(q)_{km}=\left(\begin{smallmatrix} n-k\\
n-m
\end{smallmatrix}\right)_q=C_{n-k}^{n-m}(q)$ then by (\ref{sharp})
we get
\begin{equation}\label{si1}
\sigma_1(1,n)^s_{km}=C_{m}^{k}(q).
\end{equation}
In the space of infinite matrices let us consider two operators:
\begin{equation}
\label{T} T_1:=\sum_{k\in{\mathbb Z}}(k+1)E_{kk+1},\,\,
T_{(q)}:=\sum_{k\in{\mathbb Z}}(k+1)_qE_{kk+1},
\end{equation}
where $(n)_q$ is defined by (\ref{(n)}) and $E_{km}$, are infinite
matrix  with $1$ at the place $k,m\in {\mathbb Z}$ and zeros
elsewhere. Consider the $\exp T$ and $\exp_{(q)} T_q$ of these
operators, namely
\begin{equation}
\label{expT} \exp T=\sum_{m=0}^\infty\frac{1}{m!}T^m,\,\,
\exp_{(q)} T_{(q)}=\sum_{m=0}^\infty\frac{1}{(m)!_q}T^m_q.
\end{equation}
Let us denote by $P_n$ the projector from the space of all
infinite matrices onto the subspace ${\rm Mat}(n+1,{\mathbb
C})=\{A=\sum_{0\leq k,m\leq n}a_{km}E_{km}\}$.
\begin{lem}
\label{l.Pas.exp} We have
\begin{equation}
\label{Pas.exp} P_n\exp T_1P_n=\sigma_1(1,n)^s,\quad P_n\exp_{(q)}
T_{(q)}P_n=\sigma_1(q,n)^s.
\end{equation}
\end{lem}
\begin{pf} If we set
$$
T(\nu)=\sum_{k\in{\mathbb Z}}\nu_{k+1}E_{kk+1},
$$
where $\nu_k\in {\mathbb C},\,\,k\in {\mathbb Z},$ we then have
$$
T(\nu)^m=\sum_{k\in{\mathbb
Z}}\nu_{k+1}\nu_{k+2}...\nu_{k+m}E_{kk+m}.
$$
Hence we get
$$
\exp
T_1=\sum_{m=0}^\infty\frac{1}{m!}T_1^m=\sum_{m=0}^\infty\sum_{k\in{\mathbb
Z}}\frac{(k+1)(k+2)...(k+m)}{m!}E_{kk+m}.
$$
Finally $(\exp
T_1)_{kk+m}=\frac{(k+1)(k+2)...(k+m)}{m!}=C_{k+m}^k=\sigma_1(1,n)^s_{kk+m}.$
Similarly we have
$$
\exp_{(q)}
T_{(q)}=\sum_{m=0}^\infty\frac{1}{(m)!_q}T^m_{(q)}=\sum_{m=0}^\infty\sum_{k\in{\mathbb
Z}}\frac{(k+1)_q(k+2)_q...(k+m)_q}{(m)!_q}E_{kk+m},
$$
hence $(\exp_{(q)}
T_{(q)})_{kk+m}=\frac{(k+1)_q(k+2)_q...(k+m)_q}{(m)!_q}=
C_{k+m}^k(q)=\sigma_1(q,n)^s_{kk+m}.$  \qed\end{pf}

\section{Irreducibility and equivalence of the representations}

\subsection{Operator irreducibility}
{\bf Theorem 3} {\it The representation of the group $B_3$ defined
by (\ref{Rep(q)'}) have the following properties:\\
1) for $q=1,\,\,\Lambda_n=1$, it is {\rm subspace irreducible} in arbitrary dimension $n\in{\mathbb N }$;\\
2) for $q=1,\,\,\Lambda_n={\rm diag}(\lambda_k)_{k=0}^n\not=1$ it
is {\rm operator irreducible} if and only if for any $0\leq r\leq
\left[\frac{n}{2}\right]$ there exists $0\leq i_0<i_i<...<i_r\leq
n $ such that (see (\ref{M^i_j...=0}))
$$
M^{i_0i_i...i_{n-r-1}}_{r+1r+2...n}(F_{r,n}^s(q,\lambda))\not=0;\\
$$
3) for $q\not=1,\,\,\Lambda_n=1$ it is  {\rm subspace irreducible}
if and only if $(n)_q\not=0$.\\
The representation has $[\frac{n+1}{2}]+1$ free parameters. }
%
%

We study the irreducibility of the representation (\ref{Rep(q)})
if the following cases:\\
1) $q=1$ and $\Lambda=I$ (the Humphries case);\\
2) $q=1$ and $\Lambda\not=I$
(the Tuba and Wenzl Example 1, Section 6);\\
3) $q\not=1$ and $\Lambda=I$;\\
4) $q\not=1$ and $\Lambda\not=I$.\\
{\bf Case 1)}. Let us set $T_n=\sum_{k=0}^{n-1}(n-k)E_{kk+1}$. By
Lemma \ref{l.Pas.exp} we conclude that $\sigma_1(1,n)=\exp T_n $
and $\sigma_2(1,n)=\exp(- T_n^\sharp)$.
\begin{rem}
\label{(T,T*)=Mat}
The subalgebra $\{T_n,T_n^\sharp \}$ generated by the operators
$T_n$ and $T_n^\sharp$ coincide with the algebra ${\rm
Mat}(n+1,{\mathbb C})$.
\end{rem}
\begin{rem}
\label{mat()-irr} The algebra ${\rm Mat}(n,{\mathbb C})$ of all
matrices in the space ${\mathbb C}^n$ is {\bf irreducible}.
\end{rem}
 We use the following two lemmas describing the
commutant of the operator $S(q)\Lambda_n$ (see (\ref{S(q)})) and
the commutant of a strictly upper triangular matrix $\beta$
defined as follows:
\begin{equation}
\label{beta} \beta=\sum_{k=0}^{n-1}\beta_{kk+1}E_{kk+1}=
\left(\begin{smallmatrix}
0&\beta_{01}&0&...&0\\
0&0 &\beta_{12}&0...&0\\
&&&...&\\
0&0&0&...&\beta_{n-1n}\\
0&0&0&...&0\\
\end{smallmatrix}\right).
\end{equation}
Let us fix for an operator $A=\sum_{0\leq k,m\leq n}a_{km}E_{km}$
the following decomposition
\begin{equation}
\label{decom} A=\sum_{r=-n}^nA_k,\text{\,\,
where\,\,}A_k:=\sum_{r=0}^{n-k}a_{rr+k}E_{rr+k},\,\,\,
A_{-k}:=\sum_{r=0}^{n-k}a_{r+kr}E_{r+kr},\,k\geq 0.
\end{equation}
\begin{lem}
\label{l.(beta)'} Let an operator $A\in {\rm Mat(n+1,{\mathbb
C})}$ commute  with $\beta$ defined by (\ref{beta}) and
$\beta_{kk+1}\not=0$ for all $0\leq k\leq n-1$. Then $A$ is also
upper triangular, moreover
\begin{equation}
\label{(beta)'} A=a_{0}I+\sum_{k=1}^na_k\beta^k.
\end{equation}
\end{lem}
\begin{pf}
We have
$$
(\beta A)_{km}=\beta_{kk+1}a_{k+1m},\quad 0\leq k\leq
n-1,\quad(\beta A)_{nm}=0,\,\,\,0\leq m\leq n,
$$
and
$$
(A\beta)_{km}=a_{km-1}\beta_{m-1m},\quad 1\leq m\leq n,\quad (
A\beta)_{k0}=0,\,\,\,0\leq k\leq n.
$$ Hence we have
\begin{gather}
\label{beta1}
\beta_{kk+1}a_{k+1m}=a_{km-1}\beta_{m-1m},\text{\quad or\quad}
\beta_{k-1k}a_{km}=a_{k-1,m-1}\beta_{m-1m},\\
\text{\quad for\quad} 0\leq k,m-1\leq n-1,\nonumber\\
\label{beta2} \beta_{k-1k}a_{k0}=0,\quad 1\leq k\leq n,\quad
a_{nm}\beta_{mm+1}=0,\quad 0\leq m\leq n-1.
\end{gather}
Using (\ref{beta1}) and (\ref{beta2}) we conclude that $a_{km}=0$
for $0\leq m<k\leq n$. Indeed let us take  $m=k$ in (\ref{beta1}),
then we get $\beta_{k-1k}a_{kk}=a_{k-1,k-1}\beta_{k-1k}$ or
$a_{kk}=a_{k-1,k-1}$ hence $a_{kk}=a_{00}$ for all $0\leq k\leq
n$. Finally we conclude that $A_0=a_{00}I$.

Similarly if we take $m=k+1$  in (\ref{beta1}) we get
$\beta_{k-1k}a_{kk+1}=a_{k-1k}\beta_{kk+1}$ or
$\frac{a_{k-1k}}{\beta_{k-1k}}=\frac{a_{kk+1}}{\beta_{kk+1}}=:a_1$
hence $A_1=a_1\beta$.  If we take $m=k+2$  in (\ref{beta1}) we get
using the relation $(\beta^2)_{kk+2}=\beta_{kk+1}\beta_{k+1k+2}$
$$
\frac{a_{kk+2}}{a_{k-1k+1}}=\frac{\beta_{k+1k+2}}{\beta_{k-1k}}=
\frac{(\beta^2)_{kk+2}}{(\beta^2)_{k-1k+1}},\text{\,\,\,so\,\,\,}
\frac{a_{k-1k+1}}{(\beta^2)_{k-1k+1}}=\frac{a_{kk+2}}{(\beta^2)_{kk+2}}=:a_2,
$$
hence $A_2=a_2\beta^2$.  If we put $m=k+r$ we get using the
relation $(\beta^r)_{kk+r}=\beta_{kk+1}...\beta_{k+r-1,k+r}$
$$
\frac{a_{kk+r}}{a_{k-1,k+r-1}}=\frac{\beta_{k+r-1,k+r}}{\beta_{k-1k}}=
\frac{(\beta^r)_{kk+r}}{(\beta^r)_{k-1,k+r-1}},\text{\,\,\,so\,\,\,}
\frac{a_{k-1,k+r-1}}{(\beta^r)_{k-1,k+r-1}}=\frac{a_{kk+r}}{(\beta^r)_{kk+r}}=:a_r,
$$
hence $A_r=a_r\beta^r$. This proves Lemma \ref{l.(beta)'}.
\qed\end{pf}
\begin{lem}
\label{l.(S)'} Let an operator $A\!\in\!{\rm Mat(n+1,{\mathbb
C})}$ commute  with $S(q)\Lambda$
 (see (\ref{S(q)})), then
\begin{equation}
\label{(S)'}
q_{n-k}^{-1}\lambda_{n-k}a_{km}=(-1)^{k+m}q_{k}^{-1}a_{n-k,n-m}\lambda_{m}\,\text{\,\,
where\,\,} A=(a_{km})_{0\leq k,m\leq n}.
\end{equation}
\end{lem}
\begin{pf}  We have (see (\ref{S(q)}))
$$(S(q)\Lambda
A)_{km}=S(q)_{k,n-k}\lambda_{k}a_{n-k,m}=(-1)^kq_{k}^{-1}\lambda_{k}a_{n-k,m}$$
and
$$
(AS(q)\Lambda)_{km}=a_{k,n-m}S(q)_{n-m,m}\lambda_{m}=(-1)^{n-m}q_{n-m}^{-1}a_{k,n-m}\lambda_{m}.
$$
Since $S(q)\Lambda A=AS(q)\Lambda$ we get (\ref{(S)'}).
\qed\end{pf}
 To prove the irreducibility of the representation
$\sigma_1\mapsto\sigma_1(1,n)\quad \sigma_2\mapsto\sigma_2(1,n)$
let us suppose that an operator $A$ commute with $\sigma_1(1,n)$
and $\sigma_2(1,n)$. If we set
$\beta:=(\sigma_1(1,n)-I)_1=\sum_{k=0}^{n-1}(n-k)E_{kk+1}$, the
first term in the decomposition (\ref{decom}) of the operator
$\sigma_1(1,n)-I$, by Lemma \ref{l.Pas.exp} we conclude that
$\sigma_1(1,n)=\exp \beta.$ Since $A$ commutes with
$\sigma_1(1,n)$ then $A$ commutes with $\beta=\ln\sigma_1(1,n),$
where
$$\beta=\ln\sigma_1(1,n)=\sum_{r=1}^n\frac{(-1)^r}{r}(\sigma_1(1,n)-I)^r,$$
and since
$\beta_{kk+1}=(\sigma_1(1)-I)_{kk+1}=C_{n-k}^{n-k-1}=(n-k)\not=0$
for $0\leq k\leq n-1$ we conclude by Lemma \ref{l.(beta)'} that
$A$ is upper triangular, moreover
$$
A= a_{0}I+\sum_{k=1}^na_k\beta^k= a_{0}I+\sum_{0\leq k<m\leq
n}a_{km}E_{km},
$$
i.e. $a_{km}=0$ for $k>m.$ Since $A$ commute with
$S=\sigma_1(1,n)\sigma_2(1,n)\sigma_1(1,n)$ (see (\ref{Br_n(1)}))
by Lemma \ref{l.(S)'} we get $a_{km}=(-1)^{k+m}a_{n-m,n-k}$, so
$a_{km}=0$ for $k<m$. Finally by (\ref{(beta)'}) we conclude that
$A=a_{0}I.$ Thus the irreducibility of the representation
$\sigma_1\mapsto\sigma_1(1,n),\,\,\sigma_2\mapsto\sigma_2(1,n)$ is
proved in the case 1).
\begin{rem} In fact, the representation is irreducible not only in
the operator sense (i.e. that only the trivial operators commute
with the representations) but also in the usual sense (i.e. that
there are no nontrivial invariant subspaces for operators of the
representations). It follows from the fact that formulas
$$
X\mapsto \rho_n(X)=(\sigma_1(1,n)-I)_1,\,Y\mapsto
\rho_n(Y)=(\sigma_2(1,n)-I)_{-1},\,H\mapsto [\rho_n(X),\rho_n(Y)]
$$
define the irreducible representation of the universal enveloping
algebra $U(\mathfrak{sl}_2)$ of the Lie algebra $\mathfrak{sl}_2$
(see \cite[Theorem V.4.4.]{Kas95}). Recall \cite{Kas95} that
$U(\mathfrak{sl}_2)$ is the associative algebra generated by three
generators $X,\,Y,\,H$ with the relations (\ref{sl_2}).
\begin{equation}\label{sl_2}
[H,X]=2X,\,\,[H,Y]=-2Y,\,\,[X,Y]=H,
\end{equation}
\end{rem}
{\bf Case 2). Idea of the proof.}   Let $A$ commute with
$\sigma_1^\Lambda(1,n)$ and $\sigma_2^\Lambda(1,n)$ hence by
Theorem \ref{braid.t.1}, relation (\ref{Br_n(q)}), $A$ commute
with $S(q)\Lambda$. By Lemma \ref{l.case2} (below) $A$ is upper
triangular, hence by Lemma \ref{l.(S)'} $A$ is diagonal so
$[A,\Lambda_n]=0$, hence $[A,\sigma_1(1,n)]=[A,\sigma_2(1,n)]=0$
and we are in the case 1), $n\in{\mathbb N}$. i.e. $A$ is trivial.
\begin{lem}
\label{l.case2} Let an operator $A\in {\rm Mat(n+1,{\mathbb C})}$
commute  with $\sigma_1(1,n)\Lambda_n$ where $\Lambda_n={\rm
diag}(\lambda_k)_{k=0}^n$
then $A$ is also upper triangular, i.e.
\begin{equation}
\label{case2} A=\sum_{0\leq k\leq m\leq n}a_{km}E_{km}
\end{equation}
if for any $0\leq r\leq \left[\frac{n}{2}\right]$ there exists
$0\leq i_0<i_i<...<i_r\leq n $ such that 
(\ref{M^i_j...=0})
$$
 M^{i_0i_i...i_{n-r-1}}_{r+1r+2...n} (F_{r,n}^s(1,\lambda))\not=0.
$$
\end{lem}
\begin{pf}
{\bf Let} $n=1$ and $[A,\sigma_1^{\Lambda}]
=[A,\sigma_2^{\Lambda}] =0$ where
$$
A=\left(\begin{smallmatrix}
a_{00}&a_{01}\\
a_{10}&a_{11}
\end{smallmatrix}\right),\quad
\sigma_1^{\Lambda}=\left(\begin{smallmatrix}
1&1\\
0&1
\end{smallmatrix}\right)
\left(\begin{smallmatrix}
\lambda_0&0\\
0&\lambda_1
\end{smallmatrix}\right)=\left(\begin{smallmatrix}
\lambda_0&\lambda_1\\
0&\lambda_1
\end{smallmatrix}\right)
,\quad \sigma_2^{\Lambda}=\left(\begin{smallmatrix}
\lambda_1&0\\
0&\lambda_0
\end{smallmatrix}\right)\left(\begin{smallmatrix}
1&0\\
-1&1
\end{smallmatrix}\right)=\left(\begin{smallmatrix}
\lambda_1&0\\
-\lambda_0&\lambda_0
\end{smallmatrix}\right).
$$
The relation $A\sigma_1^{\Lambda}=\sigma_1^{\Lambda}A$ gives us
$$
\left(\begin{smallmatrix}
a_{00}\lambda_0&a_{00}\lambda_1+a_{01}\lambda_1\\
a_{10}\lambda_0&a_{10}\lambda_1+a_{11}\lambda_1
\end{smallmatrix}\right)=
\left(\begin{smallmatrix}
\lambda_0a_{00}+\lambda_1a_{10}&\lambda_0a_{01}+\lambda_1a_{11}\\
\lambda_1a_{10}&\lambda_1a_{11}
\end{smallmatrix}\right)\Rightarrow
\left\{\begin{smallmatrix}
\lambda_1a_{10}=0\\
(\lambda_1-\lambda_0)a_{10}=0.\\
\end{smallmatrix}\right.
$$
Since $\lambda_1\not=0$ hence $a_{10}=0$.

{\bf Let} $n=2$ and
$[A,\sigma_1^{\Lambda}]=[A,\sigma_2^{\Lambda}]=0$ where
$$
A=\left(\begin{smallmatrix}
a_{00}&a_{01}&a_{02}\\
a_{10}&a_{11}&a_{12}\\
a_{20}&a_{21}&a_{22}
\end{smallmatrix}\right),\quad
\sigma_1^{\Lambda}=\sigma_1(1,2)\Lambda_2=\left(\begin{smallmatrix}
1&2&1\\
0&1&1\\
0&0&1\\
\end{smallmatrix}\right)
\left(\begin{smallmatrix}
\lambda_0&0&0\\
0&\lambda_1&0\\
0&0&\lambda_2
\end{smallmatrix}\right)=\left(\begin{smallmatrix}
\lambda_0&2\lambda_1&\lambda_2\\
0&\lambda_1&\lambda_2\\
0&0&\lambda_2
\end{smallmatrix}\right),
$$
 The relation $\sigma_1^{\Lambda}A=A\sigma_1^{\Lambda}$ gives
us
$$
\left(\begin{smallmatrix}
\lambda_0a_{00}+2\lambda_1a_{10}+\lambda_2a_{20}&
\lambda_0a_{01}+2\lambda_1a_{11}+\lambda_2a_{21}&
\lambda_0a_{02}+2\lambda_1a_{12}+\lambda_2a_{22}\\
\lambda_1a_{10}+\lambda_2a_{20}&\lambda_1a_{11}+\lambda_2a_{21}&\lambda_1a_{12}+\lambda_2a_{22}\\
\lambda_2a_{20}&\lambda_2a_{21}&\lambda_2a_{22}
\end{smallmatrix}\right)=
$$
$$
\left(\begin{smallmatrix}
a_{00}\lambda_0&(2a_{00}+a_{01})\lambda_1&(a_{00}+a_{01}+a_{02})\lambda_2\\
a_{10}\lambda_0&(2a_{10}+a_{11})\lambda_1&(a_{10}+a_{11}+a_{12})\lambda_2\\
a_{20}\lambda_0&(2a_{20}+a_{21})\lambda_1&(a_{20}+a_{21}+a_{22})\lambda_2\\
\end{smallmatrix}\right).
$$
If we compare the first columns we get
$$
\left\{\begin{smallmatrix}
2\lambda_1a_{10}+\lambda_2a_{20}=0\\
(\lambda_1-\lambda_0)a_{10}+\lambda_2a_{20}=0\\
(\lambda_2-\lambda_0)a_{20}=0
\end{smallmatrix}\right.
 \text{\quad or \quad}
[\sigma_1(1,2)\Lambda_2-\lambda_0I] a^{(0)}=0,\text{\,where\,}
a^{(0)}=\left(\begin{smallmatrix} 0\\
a_{10}\\
a_{20}\\
\end{smallmatrix}\right).
$$
Let $ a^{(0)}=0$. If we compare the second columns we get
$$
%
%
 \left\{\begin{smallmatrix}
\lambda_2a_{21}=0\\
(\lambda_2-\lambda_1)a_{21}=0\\
\end{smallmatrix}\right.
\text{\quad or \quad} [\sigma_1(1,2)\Lambda_2-\lambda_1I]
a^{(1)}=0,\text{\,where\,}
a^{(1)}=\left(\begin{smallmatrix} 0\\
0\\
a_{21}\\
\end{smallmatrix}\right).
$$
{\bf Let} $n=3$ and
$[A,\sigma_1^{\Lambda}]=[A,\sigma_2^{\Lambda}]=0$ where
$$
A=\left(\begin{smallmatrix}
a_{00}&a_{01}&a_{02}&a_{03}\\
a_{10}&a_{11}&a_{12}&a_{13}\\
a_{20}&a_{21}&a_{22}&a_{23}\\
a_{30}&a_{31}&a_{32}&a_{33}\\
\end{smallmatrix}\right),\quad
\sigma_1^{\Lambda}=\left(\begin{smallmatrix}
1&3&3&1\\
0&1&2&1\\
0&0&1&1\\
0&0&0&1\\
\end{smallmatrix}\right)
\left(\begin{smallmatrix}
\lambda_0&0&0&0\\
0&\lambda_1&0&0\\
0&0&\lambda_2&0\\
0&0&0&\lambda_3\\
\end{smallmatrix}\right)=\left(\begin{smallmatrix}
\lambda_0&3\lambda_1&3\lambda_2&\lambda_3\\
0&\lambda_1&2\lambda_2&\lambda_3\\
0&0&\lambda_2&\lambda_3\\
0&0&0&\lambda_3
\end{smallmatrix}\right).
$$
The relation $A\sigma_1^{\Lambda}=\sigma_1^{\Lambda}A$ gives us
$$
\left(\begin{smallmatrix}
a_{00}\lambda_0&(3a_{00}+a_{01})\lambda_1&(3a_{00}+2a_{01}+a_{02})\lambda_2
&(a_{00}+a_{01}+a_{02}+a_{03})\lambda_3\\
a_{10}\lambda_0&(3a_{10}+a_{11})\lambda_1&(3a_{10}+2a_{11}+a_{12})\lambda_2
&(a_{10}+a_{11}+a_{12}+a_{13})\lambda_3\\
a_{20}\lambda_0&(3a_{20}+a_{21})\lambda_1&(3a_{20}+2a_{21}+a_{22})\lambda_2
&(a_{20}+a_{21}+a_{22}+a_{23})\lambda_3\\
a_{30}\lambda_0&(3a_{30}+a_{31})\lambda_1&(3a_{30}+2a_{31}+a_{32})\lambda_2
&(a_{30}+a_{31}+a_{32}+a_{33})\lambda_3\\
\end{smallmatrix}\right)=
$$
$$
%
\left(\begin{smallmatrix}
\lambda_0a_{00}+3\lambda_1a_{10}+3\lambda_2a_{20}+\lambda_3a_{30}&
\lambda_0a_{01}+3\lambda_1a_{11}+3\lambda_2a_{21}+\lambda_3a_{31}\\
\lambda_1a_{10}+2\lambda_2a_{20}+\lambda_3a_{30}&
\lambda_1a_{11}+2\lambda_2a_{21}+\lambda_3a_{31}\\
\lambda_2a_{20}+\lambda_3a_{30}& \lambda_2a_{21}+\lambda_3a_{31}\\
\lambda_3a_{30}& \lambda_3a_{31}\\
\end{smallmatrix}\right.
$$
$$
\left.\begin{smallmatrix}
\lambda_0a_{02}+3\lambda_1a_{12}+3\lambda_2a_{22}+\lambda_3a_{32}&
\lambda_0a_{03}+3\lambda_1a_{13}+3\lambda_2a_{23}+\lambda_3a_{33}\\
\lambda_1a_{12}+2\lambda_2a_{22}+\lambda_3a_{32}&
\lambda_1a_{13}+2\lambda_2a_{23}+\lambda_3a_{33}\\
\lambda_2a_{22}+\lambda_3a_{32}&
\lambda_2a_{23}+\lambda_3a_{33}\\
\lambda_3a_{32}&\lambda_3a_{32}\\
\end{smallmatrix}\right).
$$
If we compare the first columns we get
$$
%
%
 \left\{\begin{smallmatrix}
3\lambda_1a_{10}+3\lambda_2a_{20}+\lambda_3a_{30}=0\\
(\lambda_1-\lambda_0)a_{10}+2\lambda_2a_{20}+\lambda_3a_{30}=0\\
(\lambda_2-\lambda_0)a_{20}+\lambda_3a_{30}= 0\\
(\lambda_3-\lambda_0)a_{30}=0\\
\end{smallmatrix}\right.\text{\quad or \quad}
[\sigma_1(1,3)\Lambda_3-\lambda_0I] a^{(0)}=0,\text{\,where\,}
a^{(0)}=\left(\begin{smallmatrix} 0\\
a_{10}\\
a_{20}\\
a_{30}\\
\end{smallmatrix}\right).
$$
Let $a^{(0)}=0$. If we compare the second columns we get the
system
$$
%
%
 \left\{\begin{smallmatrix}
2\lambda_1a_{21}+3\lambda_3a_{31}=0\\
(\lambda_2-\lambda_1)a_{21}+\lambda_3a_{31}=0\\
(\lambda_3-\lambda_1)a_{31}= 0\\
\end{smallmatrix}\right.
\text{\quad or \quad} [\sigma_1(1,3)\Lambda_3-\lambda_1I]
a^{(1)}=0,\text{\,where\,}
a^{(1)}=\left(\begin{smallmatrix} 0\\
0\\
a_{21}\\
a_{31}\\
\end{smallmatrix}\right).
$$
Let $a^{(0)}=a^{(1)}=0$. If we compare the third columns we get
the system
$$
 \left\{\begin{smallmatrix}
\lambda_3a_{32}=0\\
(\lambda_3-\lambda_2)a_{32}= 0\\
\end{smallmatrix}\right.
\text{\quad or \quad} [\sigma_1(1,3)\Lambda_3-\lambda_2I]
a^{(2)}=0,\text{\,where\,}
a^{(2)}=\left(\begin{smallmatrix} 0\\
0\\
0\\
a_{32}\\
\end{smallmatrix}\right).
$$
For the general case $n\in{\mathbb N}$ let us  consider the
following equations
\begin{equation}
\label{A-trian}
[\sigma_1(1,n)\Lambda_n-\lambda_kI]a^{(k)},\text{\,\,where\,\,}
a^{(k)}=(0,...,a_{k+1,k},a_{k+2,k},...,a_{n,k})^t,\,\,0\leq k\leq
n-1.
\end{equation}
To prove Lemma it is sufficient to show that all solutions of the
equations (\ref{A-trian}) are trivial. We rewrite the latter
equations in the following form:
\begin{equation}
\label{A-trian1} \sigma_1^{\Lambda,k}(1,n)b^{(k)}=0 ,
\text{\,\,where\,\,}
\sigma_1^{\Lambda,k}(1,n):=[\sigma_1(1,n)-\lambda_k\Lambda_n^{-1}],\,\,
b^{(k)}=\Lambda_na^{(k)},
\end{equation}
$0\leq k\leq n-1$. If we denote
\begin{equation}
\label{F_kn(1,lam)} F_{k,n}(1,\lambda)=
[\sigma_1(1,n)-\lambda_k(\Lambda_n)^{-1}]^s
\end{equation}
(for  notation $A^s$ see (\ref{sharp})) we get by Lemma
\ref{l.Pas.exp}
$$
F_{k,n}(1,\lambda)= [\sigma_1(1,n)-\lambda_k(\Lambda_n)^{-1}]^s
=\exp\left(\sum_{r=0}^{n-1} (r+1)E_{rr+1}\right)-
\lambda_k(\Lambda_n^\sharp)^{-1}=
$$
$$
\left(\begin{smallmatrix} 1-\nu_0&1&1&1&1&1&1&...\\
0&1-\nu_1^k&2&3&4&5&6&...\\
0&0&1-\nu_2^k&3&6&10&15&...\\
0&0&0&1-\nu_3^k&4&10&20&...\\
0&0&0&0&1-\nu_4^k&5&15&...\\
0&0&0&0&0&1-\nu_5^k&6&...\\
0&0&0&0&0&0&1-\nu_6^k&...\\
&&&&...&&&\\
\end{smallmatrix}\right),
$$
where $\lambda_k(\Lambda_n^\sharp)^{-1}={\rm
diag}(\nu_r^{k})_{r=0}^n$ and $\nu_r^k=\lambda_k/\lambda_{n-r}$.
Let us set $(k_n):=\sigma_1(1,n)-\lambda_k(\Lambda_n)^{-1}.$ {\bf
For} $n=2$ we get
$$
\sigma_1^{\Lambda,k}(1,2)= \sigma_1(1,2)-\lambda_k\Lambda_2^{-1}=
\left[ \left(\begin{smallmatrix}
1&2&1\\
0&1&1\\
0&0&1\\
\end{smallmatrix}\right)-\lambda_k\left(\begin{smallmatrix}
\lambda_0&0&0\\
0&\lambda_1&0\\
0&0&\lambda_2\\
\end{smallmatrix}\right)^{-1}\right]
= \left(\begin{smallmatrix}
1-\nu_0^k&2&1\\
0&1-\nu_1^k&1\\
0&0&1-\nu_2^k\\
\end{smallmatrix}\right).
$$
The equations (\ref{A-trian1}) gives us
$\sigma_1^{\Lambda,k}(1,2)b^{(k)}=0,\,\,k=0,1$ i.e.
$$
\left(\begin{smallmatrix}
0&2&1\\
0&1-\nu_1^0&1\\
0&0&1-\nu_2^0\\
\end{smallmatrix}\right)\left(\begin{smallmatrix}
0\\
b_{10}\\
b_{20}\\
\end{smallmatrix}\right)=0,\quad\left(\begin{smallmatrix}
1-\nu_0^1&2&1\\
0&0&1\\
0&0&1-\nu_2^1\\
\end{smallmatrix}\right)\left(\begin{smallmatrix}
0\\
0\\
b_{21}\\
\end{smallmatrix}\right)=0.
$$
We see that $b^{(0)}=0$ if some of minors $ M^{01}_{12}(0_2),$ $
M^{02}_{12}(0_2),\,\, M^{12}_{12}(0_2)$ are not $0$. Since
$M^{0}_{2}(0_2)=1$ we conclude that $b^{(1)}=0$. We have
$$
M^{01}_{12}(0_2)= \left|\begin{smallmatrix}
2&1\\
1-\nu_1^0&1\\
\end{smallmatrix}\right|,\,\,
 M^{02}_{12}(0_2)=\left|\begin{smallmatrix}
2&1\\
0&1-\nu_2^0\\
\end{smallmatrix}\right|,\,\,
  M^{12}_{12}(0_2)=\left|\begin{smallmatrix}
1-\nu_1^0&1\\
0&1-\nu_2^0\\
\end{smallmatrix}\right|,\,\, M^{0}_{2}(1_2)=1,
$$
hence
$$
M^{01}_{12}(0_2)= M^{01}_{12}(F_{0,2}^s(1,\lambda)),\,\,
M^{02}_{12}(0_2)= M^{02}_{12}(F_{0,2}^s(1,\lambda)),\,\,
$$
$$
M^{12}_{12}(0_2)=M^{12}_{12}(F_{0,2}^s(1,\lambda)),\,\,
M^{0}_{2}(0_2)=M^{0}_{2}(F_{1,2}^s(1,\lambda)),
$$
where
$\nu^{(r)}=(\nu^{(r)}_k)_{k=0}^2,\,\,\nu^{(r)}_k=\lambda_r/\lambda_{2-k}$
and $0\leq r\leq [\frac{2}{2}]=1$.

{\bf For} $n=3$ we have
$$
\left[ \sigma_1(1,3)-\lambda_k\Lambda_3^{-1}\right]= \left[
\left(\begin{smallmatrix}
1&3&3&1\\
0&1&2&1\\
0&0&1&1\\
0&0&0&1\\
\end{smallmatrix}\right)-\lambda_k\left(\begin{smallmatrix}
\lambda_0&0&0&0\\
0&\lambda_1&0&0\\
0&0&\lambda_2&0\\
0&0&0&\lambda_3\\
\end{smallmatrix}\right)^{-1}\right]
=
\left(\begin{smallmatrix} 1-\nu_0^k&3&3&1\\
0&1-\nu_1^k&      2&1\\
0&0      &1-\nu_2^k&1\\
0&0      &0      &1-\nu_3^k
\end{smallmatrix}\right),
$$
the equations (\ref{A-trian1}) are
$\sigma_1^{\Lambda,k}(1,3)b^{(k)}=0,\,\,k=0,1,2$ i.e.
$$
\left(\begin{smallmatrix}0&3&3&1\\
0&1-\nu_1^0&      2&1\\
0&0      &1-\nu_2^0&1\\
0&0      &0      &1-\nu_3^0
\end{smallmatrix}\right)\left(\begin{smallmatrix}
0\\
b_{10}\\
b_{20}\\
b_{30}\\
\end{smallmatrix}\right)=0,\,
\left(\begin{smallmatrix} 1-\nu_0^1&3&3&1\\
0&0&      2&1\\
0&0      &1-\nu_2^1&1\\
0&0      &0      &1-\nu_3^1
\end{smallmatrix}\right)\left(\begin{smallmatrix}
0\\
0\\
b_{21}\\
b_{31}\\
\end{smallmatrix}\right)=0,\,
\left(\begin{smallmatrix} 1-\nu_0^2&3&3&1\\
0&1-\nu_1^2&2&1\\
0&0        &0&1\\
0&0        &0&1-\nu_3^2
\end{smallmatrix}\right)\left(\begin{smallmatrix}
0\\
0\\
0\\
b_{32}\\
\end{smallmatrix}\right)=0.
$$
We see that $b^{(0)}=0$ if some of minors $
M^{i_0i_1i_2}_{123}(0_3),\,\,0\leq i_0<i_1<i_2\leq 3$ are not $0$.
Since $M^{01}_{23}(1_3)= \left|\begin{smallmatrix}
3&1\\
2&1\\
\end{smallmatrix}\right|=1\not=0$ (resp $M^{0}_{3}(2_3)=1\not=0$)
we conclude that $b^{(1)}=0$ (resp $b^{(2)}=0$). As before we
conclude that
$$
M^{i_0i_1i_2}_{123}(0_3)=M^{i_0i_1i_2}_{123}(F_{0,23}^s(1,\lambda)),\,\,
M^{01}_{23}(1_3)=M^{01}_{23}(F_{1,3}^s(1,\lambda)),
$$
where $\nu^{(r)}=(\nu^{(r)}_k)_{k=0}^3,\,\,
\nu^{(r)}_k=\lambda_r/\lambda_{3-k}$ and $0\leq r\leq
[\frac{3}{2}]=1$.

For $n=4$ and $n=5$ we have
$$
\sigma_1^\Lambda(4,\lambda_k) =
\left(\begin{smallmatrix}1-\nu_0^k&4&6&4&1\\
0&1-\nu_1^k&3&3&1\\
0&0&1-\nu_2^k&      2&1\\
0&0&0      &1-\nu_3^k&1\\
0&0&0      &0      &1-\nu_4^k
\end{smallmatrix}\right),\,\,
\sigma_1^\Lambda(5,\lambda_k) =
\left(\begin{smallmatrix}1-\nu_0^k&5&10&10&5&1\\
0&1-\nu_1^k&4&6&4&1\\
0&0&1-\nu_2^k&3&3&1\\
0&0&0&1-\nu_3^k&      2&1\\
0&0&0&0      &1-\nu_4^k&1\\
0&0&0&0      &0      &1-\nu_5^k
\end{smallmatrix}\right),
$$
the equations (\ref{A-trian1}) are
$\sigma_1^{\Lambda,k}(1,4)b^{(k)}=0,\,\,0\leq k\leq 3$ i.e.
$$
\left(\begin{smallmatrix}0&4&6&4&1\\
0&1-\nu_1^0&3&3&1\\
0&0&1-\nu_2^0&      2&1\\
0&0&0      &1-\nu_3^0&1\\
0&0&0      &0      &1-\nu_4^0
\end{smallmatrix}\right)\left(\begin{smallmatrix}
0\\
b_{10}\\
b_{20}\\
b_{30}\\
b_{40}\\
\end{smallmatrix}\right)=0,\,
\left(\begin{smallmatrix}1-\nu_0^1&4&6&4&1\\
0&0&3&3&1\\
0&0&1-\nu_2^1&      2&1\\
0&0&0      &1-\nu_3^1&1\\
0&0&0      &0      &1-\nu_4^1
\end{smallmatrix}\right)\left(\begin{smallmatrix}
0\\
0\\
b_{21}\\
b_{31}\\
b_{41}\\
\end{smallmatrix}\right)=0,
$$
$$
\left(\begin{smallmatrix}1-\nu_1^2&4&6&4&1\\
0&1-\nu_1^2&3&3&1\\
0&0&0&      2&1\\
0&0&0      &1-\nu_3^2&1\\
0&0&0      &0      &1-\nu_4^02
\end{smallmatrix}\right)\left(\begin{smallmatrix}
0\\
0\\
0\\
b_{32}\\
b_{42}\\
\end{smallmatrix}\right)=0,\,
\left(\begin{smallmatrix}1-\nu_0^3&4&6&4&1\\
0&1-\nu_2^3&3&3&1\\
0&0&1-\nu_2^3&2&1\\
0&0&0        &0&1\\
0&0&0        &0&1-\nu_4^3
\end{smallmatrix}\right)\left(\begin{smallmatrix}
0\\
0\\
0\\
0\\
b_{43}\\
\end{smallmatrix}\right)=0.
$$
 We see that $b^{(0)}=0$ if some of minors $
M^{i_0i_1i_2i_3}_{1234}(0_4),\,\,0\leq i_0<i_1<i_2<i_3\leq 4$ are
not equal to $0$. Similarly, we conclude that $b^{(1)}=0$ if some
of minors $ M^{i_0i_1i_2}_{234}(1_4),\,\,0\leq i_0<i_1<i_2\leq 4$
are not equal to $0$. Since $M^{01}_{34}(2_4)=
\left|\begin{smallmatrix}
4&1\\
3&1\\
\end{smallmatrix}\right|=1\not=0$ (resp $M^{0}_{4}(3_4)=1\not=0$)
we conclude that $b^{(2)}=0$ (resp $b^{(3)}=0$) .

In general we conclude that the system of equations
(\ref{A-trian1})
$$
\sigma_1^{\Lambda,k}(1,n)b^{(k)}=0 ,\quad 0\leq k\leq n-1
$$
has only trivial solutions $b^{(k)}=0$ if and only if for any
$0\leq r\leq \left[\frac{n}{2}\right]$ there exists $0\leq
i_0<i_i<...<i_{n-r-1}\leq n $ such that (see (\ref{M^i_j...=0}))
$$
M^{i_0i_i...i_{n-r-1}}_{r+1r+2...n}(r_n)=
M^{i_0i_i...i_{n-r-1}}_{r+1r+2...n}
(F_{r,n}^s(1,\lambda))\not=0\,\,
\text{\,\,where\,\,}\nu^{(r)}=(\nu^{(r)}_k)_{k=0}^n,\,\,
\nu^{(r)}_k=\frac{\lambda_r}{\lambda_{n-k}}.
$$
\qed\end{pf}
For $n=5$ the equations (\ref{A-trian1}) are
$\sigma_1^{\Lambda,k}(1,5)b^{(k)}=0,\,\,0\leq k\leq 4$ i.e.
$$
\left(\begin{smallmatrix}0&5&10&10&5&1\\
0&1-\nu_1^0&4&6&4&1\\
0&0&1-\nu_2^0&3&3&1\\
0&0&0&1-\nu_3^0&      2&1\\
0&0&0&0      &1-\nu_4^0&1\\
0&0&0&0      &0      &1-\nu_5^0
\end{smallmatrix}\right)\left(\begin{smallmatrix}
0\\
b_{10}\\
b_{20}\\
b_{30}\\
b_{40}\\
b_{50}\\
\end{smallmatrix}\right)=0,
\left(\begin{smallmatrix}1-\nu_0^1&5&10&10&5&1\\
0&0&4&6&4&1\\
0&0&1-\nu_2^1&3&3&1\\
0&0&0&1-\nu_3^1&      2&1\\
0&0&0&0      &1-\nu_4^1&1\\
0&0&0&0      &0      &1-\nu_5^1
\end{smallmatrix}\right)\left(\begin{smallmatrix}
0\\
0\\
b_{21}\\
b_{31}\\
b_{41}\\
b_{51}\\
\end{smallmatrix}\right)=0,
$$
$$
\left(\begin{smallmatrix}1-\nu_0^2&5&10&10&5&1\\
0&1-\nu_1^2&4&6&4&1\\
0&0&0&3&3&1\\
0&0&0&1-\nu_3^2&      2&1\\
0&0&0&0      &1-\nu_4^2&1\\
0&0&0&0      &0      &1-\nu_5^2
\end{smallmatrix}\right)\left(\begin{smallmatrix}
0\\
0\\
0\\
b_{32}\\
b_{42}\\
b_{52}\\
\end{smallmatrix}\right)=0,
\left(\begin{smallmatrix}1-\nu_0^3&5&10&10&5&1\\
0&1-\nu_1^3&4&6&4&1\\
0&0&1-\nu_2^3&3&3&1\\
0&0&0&0&      2&1\\
0&0&0&0      &1-\nu_4^3&1\\
0&0&0&0      &0      &1-\nu_5^3
\end{smallmatrix}\right)\left(\begin{smallmatrix}
0\\
0\\
0\\
0\\
b_{43}\\
b_{53}\\
\end{smallmatrix}\right)=0,
$$
$$
\left(\begin{smallmatrix}1-\nu_0^4&5&10&10&5&1\\
0&1-\nu_1^4&4&6&4&1\\
0&0&1-\nu_2^4&3&3&1\\
0&0&0&1-\nu_3^4&      2&1\\
0&0&0&0      &0&1\\
0&0&0&0      &0      &1-\nu_5^4
\end{smallmatrix}\right)\left(\begin{smallmatrix}
0\\
0\\
0\\
0\\
0\\
b_{54}\\
\end{smallmatrix}\right)=0.
$$

{\bf Definition 1.}{\it We say that the values of $\Lambda_n={\rm
diag} (\lambda_k)_{k=0}^n$ are suspected (for reducibility) if for
some  $0\leq r\leq \left[\frac{n}{2}\right]$ (see
(\ref{M^i_j...=0}))
\begin{equation}
\label{susp} M^{i_0i_i...i_{n-r-1}}_{r+1r+2...n}
(F_{rn}^s(1,\lambda))=0\,\, \text{\,\,for \,\,all\,\,\,} 0\leq
i_0<i_i<...<i_r\leq n.
\end{equation}
}
Our {\bf aim now is to describe shortly the suspected values} of
$\Lambda_n$, for example if $r=0$ we get that for all $0\leq
i_0<i_i<...<i_{n-1}\leq n$
\begin{equation}
\label{shortM...=0}
M^{i_0i_i...i_{n-1}}_{12...n}(F_{0n}^s(1,\lambda))=0\Leftrightarrow
M^{01...n-1}_{12...n} (F_{0n}^s(1,\lambda))=0,\text{\,\, and\,\,}
M_n^n(F_{0n}^s(1,\lambda))=0.
\end{equation}
{\bf To  complete the proof of  the Theorem \ref{t.IrrB_3} we
should show that representation is operator reducible for
suspected values of $\Lambda_n$ }.

Firstly we find the {\bf list} of the suspected values for
$q=1,\,\,0\leq r\leq \left[\frac{n}{2}\right]$. {\bf For}
$n=2,\,\,r=0$ we see that $ M^{01}_{12}(0_2), M^{02}_{12}(0_2),
M^{12}_{12}(0_2)$ all are zeros if and only if
$M^{01}_{12}(0_2)=0$ and $M^2_2(0_2)=0.$ Since
$$
D_2(\nu):=M^{01}_{12}(0_2)=1+\nu^0_1=(\lambda_0+\lambda_1)/\lambda_1,\text{\,\,and\,\,}
M^{2}_{2}(0)=1-\nu^0_2=(\lambda_2-\lambda_0)/\lambda_2
$$
we have the suspected value $\Lambda_2=\Lambda_2^{(2)}$ where
\begin{equation}
\label{red2}\Lambda_2^{(2)}={\rm
diag}(\lambda_0,-\lambda_0,\lambda_0) =\lambda_0{\rm
diag}(1,-1,1),\quad \text{rep. is reducible}.
\end{equation}
In this case we have
\begin{equation}
\label{reduc} \sigma_1^{\Lambda}=\left(\begin{smallmatrix}
1&-2&1\\
0&-1&1\\
0&0&1\\
\end{smallmatrix}\right),\quad
\sigma_2^{\Lambda}= \left(\begin{smallmatrix}
1&0&0\\
1&-1&0\\
1&-2&1
\end{smallmatrix}\right).
\end{equation}
Let us denote
$$
A=\left(\begin{smallmatrix}
0&-2&2\\
1&-3&1\\
2&-2&0\\
\end{smallmatrix}\right).
$$
\begin{rem} One may verify that the operator $A$ commute with
$\sigma_1^{\Lambda}$ and $\sigma_2^{\Lambda}$ defined by
(\ref{reduc}). The invariant subspace $V_2=\langle
e_2=(2,1,2)\rangle$ is generated by the eigenvector $e_2:=(2,1,2)$
for $\sigma_1^{\Lambda}$ and $\sigma_2^{\Lambda}$ i.e.
$\sigma_1^{\Lambda}e_2=e_2$ and $\sigma_2^{\Lambda}e_2=e_2$. The
representation is operator irreducible $\Leftrightarrow
\Lambda_2\not=\lambda_0{\rm diag}(1,-1,1)$.
\end{rem}


{\bf For} $n=3,\,\,r=0$ we see that all minors
$M^{i_0i_1i_2}_{123}(0_3),\,\,0\leq i_0<i_1<i_2\leq 3$ are zeros
if and only if $M^{012}_{123}(0_3)=0$ and $M^3_3(0_3)=0.$ By Lemma
\ref{l.D_n^0} we have
$$
D_3^{(0)}(\nu)^*=\left|\begin{smallmatrix} 3&3&1\\
1-\nu_1&2&1\\
0&1-\nu_2&1
\end{smallmatrix}\right|=\frac{2\lambda_0}{\lambda_1\lambda_2}
\left(\lambda_0+\lambda_1+\lambda_2\right)\text{\,\,and\,\,}
M^{3}_{3}(0_3)=1-\nu^0_3=(\lambda_3-\lambda_0)/\lambda_3
$$
hence the suspected  $\Lambda_3$ is as follows
$\Lambda_3^{(3)}:={\rm
diag}(\lambda_0,\lambda_1,\lambda_2,\lambda_0)$ with
$\lambda_0+\lambda_1+\lambda_2=0$ or
$$
\Lambda_3^{(3)}=\lambda_0{\rm
diag}(1,\alpha_1,\alpha_2,1)\text{\,\,\,with\,\,\,}1+\alpha_1+\alpha_2=0.
$$

{\bf For} $n=4,\,\,r=0$ we get that all minors
$M^{i_0i_1i_2i_3}_{1234}(0_4),\,\,0\leq i_0<i_1<i_2<i_3\leq 4$ are
zeros if and only if $M^{0123}_{1234}(0_4)=0$ and $M^4_4(0_4)=0.$
By Lemma \ref{l.D_n^0} we have
$$
\Lambda_4^{(4)}=\lambda_0{\rm
diag}(1,\alpha_1,\alpha_2,\alpha_3,1)\text{\,\,\,with\,\,\,}1+\alpha_1+\alpha_2+\alpha_3=0.
$$
{\bf For} $r=1$ we get that all minors
$M^{i_0i_1i_2}_{234}(1_4),\,\,0\leq i_0<i_1<i_2\leq 4$ are zeros
if and only if $ M^{i_0i_1i_2}_{234}(1_4),\,\,0\leq
i_0<i_1<i_2\leq 3\text{\,\,and\,\,}M^4_4(1_4)=0. $

The general rule is similar
$M^{i_0i_i...i_{n-r-1}}_{r+1r+2...n}(r_n)=0,\,\,0\leq
i_0<i_i<...<i_{n-1}\leq n\Leftrightarrow
M^{i_0i_i...i_{n-r-1}}_{r+1r+2...n}(r_n)=0,\,\,0\leq
i_0<i_i<...<i_{n-1}\leq n-1$ and $M^n_n(r_n)=0$.\\
For $r=0$ and  the general case $n\in {\mathbb N}$  we should
calculate the following determinants:
\begin{equation}
\label{det} D_n^{(0)}(\nu):= M^{01...n-1}_{12...n}\left[
\sigma_1(1,n)-\lambda_0\Lambda_n^{-1}\right],
\end{equation}
\begin{equation}
\label{det_k} D_n^{(k)}(\nu):= M^{01...k-2}_{k+1k+2...n}\left[
\sigma_1(1,n)-\lambda_k\Lambda_n^{-1}\right],\quad 0\leq k\leq
[n/2].
\end{equation}
Let us denote by (*) the conditions (see Remark 4.5)
$\lambda_r\lambda_{n-r}=c,\,\,1\leq r\leq n$ and by
$D_n^{(k)}(\nu)^*$ the value of $D_n^{(k)}(\nu)$ under these
conditions.
\begin{lem}
\label{l.D_n^0}
 We have
$$
D_n^{(0)}(\nu):= M^{01...n-1}_{23...n}\left[
\sigma_1(1,n)-\lambda_0\Lambda_n^{-1}\right]=1+
\sum_{r=1}^{n-1}\sum_{1\leq i_1<i_2<...<i_r\leq
n-1}a_{i_1i_2...i_r}\nu_{i_1}\nu_{i_2}...\nu_{i_r},
$$
\begin{equation}
\label{red_n} D_n^{(0)}(\nu)^*
=\frac{(n-2)!\lambda_0^{n-2}}{\prod_{k=1}^{n-1}\lambda_k}
\sum_{k=0}^{n-1}\lambda_k.
\end{equation}
\end{lem}
\begin{pf}
For the following notion and lemma below see \cite{KosAlb06J}. We
define $G_m(\lambda)$ {\it the generalization of the
characteristic polynomial}
 $p_C(t)={\rm
det}\,(tI-C),\, t\in {\mathbb C}$ of the matrix $C\in {\rm
Mat}(m,{\mathbb C})$:
\begin{equation}
\label{d.G_k(lambda)}
 G_m(\lambda)={\rm
det}\,C_m(\lambda) ,\,\,\lambda \in {\mathbb C}^m,\quad{\rm
where}\quad C_m(\lambda)=C+\sum_{k=1}^m\lambda_kE_{kk}.
\end{equation}
 We denote by
$$
M^{i_1i_2...i_r}_{j_1j_2...j_r}(C),\,\,({\rm resp.\,\,}
A^{i_1i_2...i_r}_{j_1j_2...j_r}(C)),\,\, 1\leq i_1<...<i_r\leq
m,\,\, 1\leq j_1<...<j_r\leq m
$$
the  minors (resp. the cofactors) of the matrix $C$ with
$i_1,i_2,...,i_r$ rows and $j_1,j_2,...,j_r$ columns. By
definition $A^{12...m}_{12...m}(C)=M^\emptyset_\emptyset(C)=1$ and
$M^{12...m}_{12...m}(C)=A^\emptyset_\emptyset(C)={\rm det}\,{C}$.
\begin{lem}
\label{l.detC-LI} For the generalized characteristic polynomial
 $G_m(\lambda)$ of  $C\in {\rm Mat}(m,{\mathbb C})$ and
 $\lambda=(\lambda_1,\lambda_2,...,\lambda_m)\in {\mathbb
C}^m$ we have:  $G_m(\lambda)=$
\begin{equation}
\label{detC-LI}
 {\rm
det}\,\left(C+\sum_{k=1}^m\lambda_kE_{kk}\right)\!=\!{\rm det}\,C+
\sum_{r=1}^m\sum_{1\leq i_1<i_2<...<i_r\leq
m}\lambda_{i_1}\lambda_{i_2}...\lambda_{i_r}A^{i_1i_2...i_r}_{i_1i_2...i_r}(C).
\end{equation}
\end{lem}
\begin {rem}
If we set
$\lambda_\alpha=\lambda_{i_1}\lambda_{i_2}...\lambda_{i_r}$ where
$\alpha=(i_1,i_2,...,i_r)$ and
$A^\alpha_\alpha(C)=A^{i_1i_2...i_r}_{i_1i_2...i_r}(C),\,\,
\lambda_\emptyset=1,\,\,A^\emptyset_\emptyset(C)={\rm det}\,C$ we
may write (\ref{detC-LI}) as follows:
\begin{equation}
\label{detC-LI.2} G_m(\lambda)={\rm det}\,C_m(\lambda)=
\sum_{\emptyset\subseteq\alpha\subseteq\{1,2,...,m\}}\lambda_\alpha
A^\alpha_\alpha(C).
\end{equation}
\end{rem}
Denote by $C_n$ the matrix corresponding to minor
$M^{01...n-1}_{23...n}\left[\sigma_1(1,n)\right]$. Using Lemma
\ref{l.detC-LI} we have  for $n=2,3,4$ and $n=5$ where
$\nu_k=\nu^{(0)}_k=\lambda_0/\lambda_k$
$$
D_2^{(0)}(\nu)= \left|\begin{smallmatrix}
2      &1\\
1-\nu_1&1\\
\end{smallmatrix}\right|={\rm det}C_2+\nu_1A^1_0(C_2)=
1+\nu_1,
$$
$$
 D_3^{(0)}(\nu)=
\left|\begin{smallmatrix} 3&3&1\\
1-\nu_1&2&1\\
0&1-\nu_2&1
\end{smallmatrix}\right|={\rm det}C_3+\nu_1A^1_0(C_3)+\nu_2A^2_1(C_3)+
\nu_1\nu_2A^{12}_{01}(C_3) =1+2\nu_1+2\nu_2+\nu_1\nu_2,
$$
$$
D_4^{(0)}(\nu)= \left|\begin{smallmatrix}
4&6&4&1\\
1-\nu_1&3&3&1\\
0&1-\nu_2&2&1\\
0&0&1-\nu_3&1
\end{smallmatrix}\right|={\rm det}C_4+\nu_1A^1_0(C_4)+\nu_2A^2_1(C_4)+
\nu_3A^3_2(C_4)+ \nu_1\nu_2A^{12}_{01}(C_4)$$
$$+
\nu_1\nu_3A^{13}_{02}(C_4) + \nu_2\nu_3A^{23}_{12}(C_4)+
\nu_1\nu_2\nu_3A^{123}_{012}(C_4)=
$$
$$
1+3\nu_1+5\nu_2+3\nu_3+3\nu_1\nu_2+5\nu_1\nu_3+3\nu_2\nu_3+
\nu_1\nu_2\nu_3,
$$
$$
D_5^{(0)}(\nu)= \left|\begin{smallmatrix}5&10&10&5&1\\
1-\nu_1&4&6&4&1\\
0&1-\nu_2&3&3&1\\
0&0&1-\nu_3&2&1\\
0&0&0&1-\nu_4&1
\end{smallmatrix}\right|=
1+4\nu_1+9\nu_2+9\nu_3+4\nu_4+
$$
$$
6\nu_1\nu_2
+16\nu_1\nu_3+11\nu_1\nu_4+11\nu_2\nu_3+16\nu_2\nu_4+6\nu_3\nu_4
$$
$$
+4\nu_2\nu_3\nu_4+9\nu_1\nu_3\nu_4+9\nu_1\nu_2\nu_4+4\nu_1\nu_2\nu_3+
\nu_1\nu_2\nu_3\nu_4.
$$
To prove the general formulas we use Lemma \ref{l.detC-LI}. We
have
$$
D_n^{(0)}(\nu)= \!{\rm det}\,C_n+ \sum_{r=1}^{n-1}\sum_{1\leq
i_1<i_2<...<i_r\leq
m}\nu_{i_1}\nu_{i_2}...\nu_{i_r}A^{i_1i_2...i_r}_{i_1-1i_2-1...i_r-1}(C_n).
$$

To prove (\ref{red_n}) we get for $n=2,3$
$$
D_2^{(0)}(\nu)= \left|\begin{smallmatrix}
2      &1\\
1-\nu_1&1\\
\end{smallmatrix}\right|=
1+\nu_1=1+\lambda_0/\lambda_1=(\lambda_0+\lambda_1)/\lambda_1,
$$
$$
D_3^{(0)}(\nu)=
\left|\begin{smallmatrix} 3&3&1\\
1-\nu_1&2&1\\
0&1-\nu_2&1
\end{smallmatrix}\right|
=1+2\nu_1+2\nu_2+\nu_1\nu_2=\frac{1}{\lambda_1\lambda_2}\left(\lambda_1\lambda_2+
2\lambda_0\lambda_2+2\lambda_0\lambda_1+\lambda_0^2\right),\quad
$$
$$
(*)\Rightarrow\lambda_1\lambda_2=\lambda_0^2\Rightarrow
D_3^{(0)}(\nu)^*=\frac{2\lambda_0}{\lambda_1\lambda_2}
\left(\lambda_0+\lambda_1+\lambda_2\right).
$$
$$
D_4^{(0)}(\nu)=
1+3\nu_1+5\nu_2+3\nu_3+3\nu_1\nu_2+5\nu_1\nu_3+3\nu_2\nu_3+
\nu_1\nu_2\nu_3,
$$
$$
D_4^{(0)}(\nu^{(0)})=1+3\frac{\lambda_0}{\lambda_1}+5\frac{\lambda_0}{\lambda_2}+
3\frac{\lambda_0}{\lambda_3}+3\frac{\lambda_0^2}{\lambda_1\lambda_2}+
5\frac{\lambda_0^2}{\lambda_1\lambda_3}+3\frac{\lambda_0^2}{\lambda_2\lambda_3}+
\frac{\lambda_0^3}{\lambda_1\lambda_2\lambda_3}=
$$
$$
\left(\lambda_1\lambda_2\lambda_3+3\lambda_0\lambda_2\lambda_3+5\lambda_0\lambda_1\lambda_3
+3\lambda_0\lambda_1\lambda_2+3\lambda_0^2\lambda_3+5\lambda_0^2\lambda_2+
3\lambda_0^2\lambda_1+3\lambda_0^3
\right)/\lambda_1\lambda_2\lambda_3=
$$
(since $\lambda_0\lambda_4=\lambda_1\lambda_3=\lambda_2^2$ and
$\lambda_0=\lambda_4$ we get
$\lambda_0^2=\lambda_1\lambda_3=\lambda_2^2$ so
$\lambda_0=\pm\lambda_2$. If $\lambda_0=\lambda_2$ we get)
$$
\left( \lambda_0^2(\lambda_2\pm 3\lambda_3+5\lambda_0 \pm
3\lambda_1)+\lambda_0^2(3\lambda_3+5\lambda_2+
3\lambda_1+\lambda_0 \right)/\lambda_1\lambda_2\lambda_3
$$
$$
\frac{6\lambda_0^2}{\lambda_1\lambda_2\lambda_3}( \lambda_0+(1\pm
1 )/2 \lambda_1+\lambda_2+(1\pm 1)/2\lambda_3 ).
$$
$$
D_4^{(0)}(\nu)^* =\frac{6\lambda_0^2}{\lambda_1\lambda_2\lambda_3}
\left\{\begin{smallmatrix}
\left(\lambda_0+\lambda_1+\lambda_2+\lambda_3\right)&\text{\,\,if\,\,}
\lambda_0=\lambda_2\\
\left(\lambda_0+\lambda_2\right)&\text{\,\,if\,\,}
\lambda_0=-\lambda_2\\
\end{smallmatrix}\right.
$$
$$
D_5^{(0)}(\nu)= 1+4\nu_1+9\nu_2+9\nu_3+4\nu_4+
$$
$$
6\nu_1\nu_2
+16\nu_1\nu_3+11\nu_1\nu_4+11\nu_2\nu_3+16\nu_2\nu_4+6\nu_3\nu_4
$$
$$
+4\nu_2\nu_3\nu_4+9\nu_1\nu_3\nu_4+9\nu_1\nu_2\nu_4+4\nu_1\nu_2\nu_3+
\nu_1\nu_2\nu_3\nu_4 ,
$$
$$
D_5^{(0)}(\nu)=
1+4\frac{\lambda_0}{\lambda_1}+9\frac{\lambda_0}{\lambda_2}+9\frac{\lambda_0}{\lambda_3}
+4\frac{\lambda_0}{\lambda_4}+
$$
$$
6\frac{\lambda_0^2}{\lambda_1\lambda_2}
+16\frac{\lambda_0^2}{\lambda_1\lambda_3}+11\frac{\lambda_0^2}{\lambda_1\lambda_4}
+11\frac{\lambda_0^2}{\lambda_2\lambda_3}+16\frac{\lambda_0^2}{\lambda_2\lambda_4}
+6\frac{\lambda_0^2}{\lambda_3\lambda_4}
$$
$$
+4\frac{\lambda_0^3}{\lambda_2\lambda_3\lambda_4}+9\frac{\lambda_0^3}
{\lambda_1\lambda_3\lambda_4}+9\frac{\lambda_0^3}{\lambda_1\lambda_2\lambda_4}
+4\frac{\lambda_0^3}{\lambda_1\lambda_2\lambda_3}+
\frac{\lambda_0^4}{\lambda_1\lambda_2\lambda_3\lambda_4} ,
$$
$$
(*)\Rightarrow D_5^{(0)}(\nu)^*
=\frac{24\lambda_0^3}{\lambda_1\lambda_2\lambda_3\lambda_4}
\left(\lambda_0+\lambda_1+\lambda_2+\lambda_3+\lambda_4\right).
$$

{\bf For} $n=3$ we get (we write $\nu_k$ for $\nu_k^0$)
$$
D_3^{(0)}(\nu)^*=\frac{2\lambda_0}{\lambda_1\lambda_2}
\left(\lambda_0+\lambda_1+\lambda_2\right)=
\frac{2\lambda_0^2}{\lambda_1\lambda_2}(1+\alpha_1+\alpha_2)=0
$$
$$
\Lambda_3={\rm
diag}(1,\alpha_1,\alpha_2,1),\,\,\alpha_1+\alpha_2=-1,\,\,
\alpha_1\alpha_2=1,
$$
\begin{equation}
\label{alpha3}
\text{equation\,\,}\alpha^2+\alpha+1=0,\,\,\alpha_{1,2}=-1/2\pm\sqrt{1/4-1}=(-1\pm
i\sqrt{3})/2.
\end{equation}
\begin{equation}
\label{red3} \Lambda_3^{(3)} =\lambda_0{\rm
diag}(1,\alpha_1,\alpha_2,1),\,\,\{\alpha_1,\alpha_2\}=\{\exp 2\pi
i/3,\exp 4\pi i/3 \}.
\end{equation}
{\bf For} $n=4$ we get $D_4^{(0)}(\nu)^* =6\lambda_0^2
\left(\lambda_0+\lambda_1+\lambda_2+\lambda_3\right)/(\lambda_1\lambda_2\lambda_3)$
$$
\text{Since\,\,}\Lambda_4={\rm
diag}(\lambda_0,\lambda_1,\lambda_2,\lambda_3,\lambda_4)\text{\,\,with\,\,}
\lambda_0\lambda_4=\lambda_1\lambda_3=\lambda_2^2\text{\,\,we
have}
$$
$$
\Lambda_4=\lambda_0{\rm
diag}(1,\alpha_1,\alpha_2,\alpha_3,1),\,\,\text{where\,\,}\alpha_k=\lambda_k/\lambda_0,
$$
$$
\alpha_1+\alpha_2+\alpha_3=-1\text{\,\,with\,\,}\alpha_1\alpha_3=\alpha_2^2=1.
$$
Indeed we have $\alpha_2=\pm 1$. a) let $\alpha_2=1$ then we have
$$
\alpha_1+\alpha_3=-2,\,\,\alpha_1\alpha_3=1,\text{\,\,equation\,\,}
\alpha^2+2\alpha+1=0,\,\,(\alpha+1)^2=0,\,\,\alpha_{1,3}=-1
$$
\begin{equation}
\label{red4a} \text{hence\,\,}\Lambda_4^{(2)} =\lambda_0{\rm
diag}(1,-1,1,-1,1),\quad
(\sigma_1^\Lambda)^2=(\sigma_2^\Lambda)^2=1,\,\,\, \text{\bf rep.
is reducible}.
\end{equation}
b) let $\alpha_2=-1$ then we have $
\alpha_1+\alpha_3=0,\,\,\alpha_1\alpha_3=1,\text{\,\,equation\,\,}
\alpha^2+1=0,\,\,\alpha_{1,3}=\pm i, $
\begin{equation}
\label{red4b} \Rightarrow\Lambda_4^{(4)} =\lambda_0{\rm
diag}(1,\pm i,-1,\mp i,1).
\end{equation}
{\bf For} $n=5$ we get $D_5^{(0)}(\nu)^*
=\frac{24\lambda_0^3}{\lambda_1\lambda_2\lambda_3\lambda_4}
\left(\lambda_0+\lambda_1+\lambda_2+\lambda_3+\lambda_4\right), $
$$
\text{since\,\,}\Lambda_5={\rm
diag}(\lambda_0,\lambda_1,\lambda_2,\lambda_3,\lambda_4,\lambda_5)\text{\,\,with\,\,}
\lambda_0\lambda_5=\lambda_1\lambda_4=\lambda_2\lambda_3.
$$
$$
\Lambda_5=\lambda_0{\rm
diag}(1,\alpha_1,\alpha_2,\alpha_3,\alpha_4,1),\,\,\alpha_k=\lambda_k/\lambda_0
$$
$$
\alpha_1+\alpha_2+\alpha_3+\alpha_4=-1\text{\,\,with\,\,}
\alpha_1\alpha_4=\alpha_2\alpha_3=1,\quad \text{reducible} .
$$
$$
\left\{\begin{smallmatrix} \alpha_1+\alpha_4=k,\\
 \alpha_1\alpha_4=1
\end{smallmatrix}\right.
\left\{\begin{smallmatrix} \alpha_2+\alpha_3=-(1+k),\\
 \alpha_2\alpha_3=1
\end{smallmatrix}\right.
$$
$$
\alpha^2-k\alpha+1=0,\,\,\alpha^2+(1+k)\alpha+1=0,
$$
$$
\alpha_{1,4}=k/2\pm\sqrt{(k/2)^2-1},\,\,\alpha_{2,3}=-(1+k)/2\pm\sqrt{[(1+k)/2]^2-1}.
$$
$$
\Lambda_n=\lambda_0{\rm
diag}(1,\alpha_1,\alpha_2,...,\alpha_{n-1},1),\,\,\sum_{k=1}^{n-1}\alpha_k=-1,\,\,
\alpha_k\alpha_{n-k}=1,\,\,1\leq k\leq n-1.
$$
\begin{equation}
\label{red,q=1} \Lambda_n={\rm
diag}(\lambda_0,\lambda_1,...,\lambda_{n-1},\lambda_n),\quad
\lambda_0=\lambda_n,\quad\text{and}\quad
\sum_{k=0}^{n-1}\lambda_k=0.
\end{equation}
The solution in the general case are (see
(\ref{red2})--(\ref{red4b}))
\begin{equation}
\label{red_2m+1}
 \Lambda_n=\lambda_0{\rm
diag}(\alpha_k)_{k=0}^n,\,\,\alpha_k
=\exp( \pm\frac{2\pi i
k}{n}),\,\,\text{for\,\,\,}n=2m+1.
\end{equation}
\begin{equation}
\label{red_2m}
 \Lambda_n=\lambda_0{\rm
diag}(\alpha_k)_{k=0}^n,\,\,\alpha_k^{(0)}=\exp( \pm\frac{2\pi i
k}{2m}),\,\,\alpha_k^{(1)}=\exp( \pm\frac{2\pi i
k}{m})\,\,\text{for\,\,\,}n=2m.
\end{equation}
{\bf For} $n=2$ we have
$$
\Lambda_2 =\lambda_0{\rm diag}(1,-1,1),\,\,
\alpha_1+\alpha_2=-1,\,\,\exp\left(\frac{2\pi i k}{2}\right),
$$
{\bf For} $n=3$ we have $ \alpha_1+\alpha_2=-1,\,\,
\alpha_1\alpha_2=1,\,\,\alpha^2+\alpha+1=0 $
$$
\Lambda_3^{(3)} =\lambda_0{\rm
diag}(\alpha^0,\alpha,\alpha^2,\alpha^3),\,\,\alpha^3=1,\,\,\alpha\not=1.
$$
{\bf For} $n=4$ we have
$$
\Lambda_4=\lambda_0{\rm diag}(1,\alpha_1,\alpha_2,\alpha_3,1),
$$
$$
\Lambda_4^{(2)} =\lambda_0{\rm diag}(1,-1,1,-1,1),\,\,\exp(
\pm\frac{2\pi i k}{2}),\,\,0\leq k\leq 4,
$$
$$
\alpha_1+\alpha_2+\alpha_3=-1,\,\,\alpha_1\alpha_3=\alpha_2^2=1,\,\,
$$
$$
a)\,\,\alpha_2=1,\,\,\alpha_1+\alpha_3=-2,\,\alpha_1\alpha_3=\alpha_2^2=1\Rightarrow\alpha^2+2\alpha+1=0
$$
$$
\quad \Lambda_4^{(4)} =\lambda_0{\rm diag}(1,\pm i,-1,\mp
i,1),\,\,\exp( \pm\frac{2\pi i k}{4}),\,\,0\leq k\leq 4.
$$
$$
b)\,\,\alpha_2=-1,\,\,\alpha_1+\alpha_3=0,
\,\alpha_1\alpha_3=1\Rightarrow\alpha^2+1=0.
$$
{\bf For} $n=5$ we have
$$
\Lambda_5=\lambda_0{\rm
diag}(1,\alpha_1,\alpha_2,\alpha_3,\alpha_4,1),
$$
$$
\alpha_1+\alpha_2+\alpha_3+\alpha_4=-1,\text
{\,\,with\,\,}\alpha_1\alpha_4=\alpha_2\alpha_3=1.
$$
$$
\Lambda_5^{(5)}(\alpha_1,\alpha_2)=\lambda_0{\rm
diag}(1,\alpha_1,\alpha_2,\alpha_2^{-1},\alpha_1^{-1},1),\,\,
\alpha_1+\alpha_1^{-1}+\alpha_2+\alpha_2^{-1}=0.
$$
In particular we have
$$
\Lambda_5^{(5)} =\lambda_0{\rm
diag}(\alpha^k)_{k=0}^5,\,\,\alpha^5=1,\,\,\alpha\not=1.
$$
{\bf For} $n=6$ we get
$$
\Lambda_6=\lambda_0{\rm
diag}(1,\alpha_1,\alpha_2,\alpha_3,\alpha_4,\alpha_5,1),
$$
$$
\alpha_1+\alpha_2+\alpha_3+\alpha_4+\alpha_5=-1
\text{\,\,with\,\,}\alpha_1\alpha_5=\alpha_2\alpha_4=\alpha_3^2=1,
$$
$$
a)
\,\,\alpha_3=1,\,\Lambda_6^{(6),1}(\alpha_1,\alpha_2)=\lambda_0{\rm
diag}(1,\alpha_1,\alpha_2,1,\alpha_2^{-1},\alpha_1^{-1},1),\,\,
\alpha_1+\alpha_1^{-1}+\alpha_2+\alpha_2^{-1}=0,
$$
$$
b)
\,\,\alpha_3=-1,\,\Lambda_6^{(6),-1}(\alpha_1,\alpha_2)=\lambda_0{\rm
diag}(1,\alpha_1,\alpha_2,-1,\alpha_2^{-1},\alpha_1^{-1},1),\,\,
\alpha_1+\alpha_1^{-1}+\alpha_2+\alpha_2^{-1}=-2,
$$
In particular we have
$$
\Lambda_6^{(6)} =\lambda_0{\rm
diag}(\alpha^k)_{k=0}^6,\,\,\alpha^6=1,\,\,\alpha\not=1.
$$
In the {\bf general case} we have for $n=2m+1$
$$
\Lambda_{2m+1}=\lambda_0{\rm diag}(1,\alpha_1,...,\alpha_{2m},1),
$$
$$
\sum_{k=1}^{2m}\alpha_k=-1,\text{\,\,with\,\,}\alpha_k\alpha_{2m-k}=
\alpha_1\alpha_{2m}=1,
$$
$$
\Lambda_{2m+1}^{(2m+1)}(\alpha_1,...,\alpha_m)=\lambda_0{\rm
diag}(1,\alpha_1,...,\alpha_m,\alpha_m^{-1},...,\alpha_1^{-1},1),\,\,
\sum_{k=1}^m(\alpha_k+\alpha_k^{-1})=0.
$$
{\bf For} $n=2m+2$ we have $ \Lambda_{2m+2}=\lambda_0{\rm
diag}(1,\alpha_1,...,\alpha_{2m+1},1), $
$$
\sum_{k=1}^{2m+1}\alpha_k=-1,\text{\,\,with\,\,}\alpha_k\alpha_{2m+1-k}=
\alpha_1\alpha_{2m+1}=1,
$$
$$
a) \,\,
\Lambda_{2m+2}^{(2m+2),1}(\alpha_1,...,\alpha_2)=\lambda_0{\rm
diag}(1,\alpha_1,...,\alpha_m,1,\alpha_m^{-1},...,\alpha_1^{-1},1),\,\,
\sum_{k=1}^m(\alpha_k+\alpha_k^{-1})=0,
$$
$$
b) \,\,
\Lambda_{2m+2}^{(2m+2),-1}(\alpha_1,\alpha_2)=\lambda_0{\rm
diag}(1,\alpha_1,...,\alpha_m,-1,\alpha_m^{-1},...,\alpha_1^{-1},1),\,\,
\sum_{k=1}^m(\alpha_k+\alpha_k^{-1})=-2.
$$
In particular we have in both cases:
$$
\Lambda_{n}^{(n)} =\lambda_0{\rm
diag}(\alpha^k)_{k=0}^{n},\,\,\alpha^{n}=1,\,\,\alpha\not=1.
$$
\qed\end{pf}
%

 {\bf Case 3).} We prove the {\bf irreducibility}  of the
representation
$$
\sigma_1\mapsto\sigma_1^D=\sigma_1(q,n)D_n^\sharp(q)\quad
\sigma_2\mapsto \sigma_2^D=D_n(q)\sigma_2(q,n),
$$ where (see
(\ref{D_n(q)})) $D_n(q)={\rm diag}(q_r)_{r=0}^n$. By Lemma
\ref{l.case3} below we show that the operator $A$, commuting with
$\sigma_1(q,n)D_n^\sharp(q)$ is upper triangular under certain
conditions. Further, by relation (\ref{Br_n(q)}) $A$ commute with
$S(q)\Lambda$ hence by Lemma \ref{l.(S)'} $A$ is diagonal: $A={\rm
diag}(a_{00},...,a_{nn}).$ Using again the commutation
$\sigma_1(q,n)\Lambda A\Lambda^{-1}=A\sigma_1(q,n)$ we get
$\sigma_1(q,n) A=A\sigma_1(q,n)$, since $\Lambda A\Lambda^{-1}=A$,
hence, by Lemma \ref{l.Pas.exp} $A$ commute with
$\beta(q)=\ln_{(q)}\sigma_1(q,n)=(\sigma_1(q,n)-I)_1$ where
$$
\ln_{(q)}\sigma_1(q,n)=\sum_{r=1}^n\frac{(-1)^r}{(r)_q}(\sigma_1(q,n)-I)^r
$$
and if
$\beta_{kk+1}(q,n):=(\sigma_1(q,n)-I)_{kk+1}=C_{n-k}^{n-k-1}(q)\not=0$
we conclude by Lemma \ref{l.(beta)'} that $A$  is trivial.

{\bf Definition 2.}{\it We say that the value of $q$ is suspected
(for reducibility of the representation $\sigma^D(q,n)$) if for
some $2\leq r\leq n$ holds $(r)_q=0$.}
\begin{lem}
\label{l.case3} Let an operator $A\in {\rm Mat(n+1,{\mathbb C})}$
commute  with $\sigma_1(q,n)D_n^\sharp(q)$.
then $A$ is also upper triangular, i.e.
\begin{equation}
\label{Acase3} A=\sum_{0\leq k\leq m\leq n}a_{km}E_{km}.
\end{equation}
if for any $0\leq r\leq \left[\frac{n}{2}\right]$ there exists
$0\leq i_0<i_i<...<i_r\leq n $ such that
$$
 M^{i_0i_i...i_{n-r-1}}_{r+1r+2...n} (F_{rn}^s(q,1))\not=0.\,\,
$$
\end{lem}
\begin{pf} For $n=1$ we get
$\sigma_1(q,1)D_1^\sharp(q)=\sigma_1(1,1)$ hence we are in the
case 1) i.e. $q=1$. For $n=2$ we have (see (\ref{si_1(q)}) and
(\ref{D_n(q)}))
$$
\sigma_1(q,2) =\left(\begin{smallmatrix}
1&1+q&1\\
0&1&1\\
0&0&1
\end{smallmatrix}\right),\quad
D_2^\sharp(q)=\left(\begin{smallmatrix}
q&0&0\\
0&1&0\\
0&0&1
\end{smallmatrix}\right),\quad
 \sigma_1(q,2)D_2^\sharp(q)= \left(\begin{smallmatrix}
q&1+q&1\\
0&1&1\\
0&0&1
\end{smallmatrix}\right),
$$
$$
\left(\begin{smallmatrix}
q&1+q&1\\
0&1&1\\
0&0&1
\end{smallmatrix}\right)
\left(\begin{smallmatrix}
a_{00}&a_{01}&a_{02}\\
a_{10}&a_{11}&a_{12}\\
a_{20}&a_{21}&a_{22}
\end{smallmatrix}\right)=
\left(\begin{smallmatrix}
qa_{00}+(1+q)a_{10}+a_{20}&qa_{01}+(1+q)a_{11}+a_{21}&qa_{21}+(1+q)a_{12}+a_{22}\\
a_{10}+a_{20}&a_{11}+a_{21}&a_{12}+a_{22}\\
a_{20}&a_{21}&a_{22}
\end{smallmatrix}\right),
$$
$$
\left(\begin{smallmatrix}
a_{00}&a_{01}&a_{02}\\
a_{10}&a_{11}&a_{12}\\
a_{20}&a_{21}&a_{22}
\end{smallmatrix}\right)\left(\begin{smallmatrix}
q&1+q&1\\
0&1&1\\
0&0&1
\end{smallmatrix}\right)=
\left(\begin{smallmatrix}
a_{00}q\,&\,a_{00}(1+q)+a_{01}&a_{00}+a_{01}+a_{02}\\
a_{10}q\,&\,a_{10}(1+q)+a_{11}&a_{10}+a_{11}+a_{12}\\
a_{20}q\,&\,a_{20}(1+q)+a_{21}&a_{20}+a_{21}+a_{22}
\end{smallmatrix}\right).
$$
If we compare the first columns we get
$$
\left\{\begin{smallmatrix}
(1+q)a_{10}+a_{20}=0\\
(1-q)a_{10}+a_{20}=0\\
(1-q)a_{20}=0.\\
\end{smallmatrix}\right.
 \text{\quad or \quad}
[\sigma_1(q,2)D_2^\sharp(q)-qI] a^{(0)}=0,\text{\,where\,}
a^{(0)}=\left(\begin{smallmatrix} 0\\
a_{10}\\
a_{20}\\
\end{smallmatrix}\right).
$$
Let $ a^{(0)}=0$. If we compare the second columns we get
$$
\left\{\begin{smallmatrix}
a_{21}=0\\
a_{21}=0\\
\end{smallmatrix}\right.
\text{\quad or \quad} [\sigma_1(q,2)D_2^\sharp(q)-I]
a^{(1)}=0,\text{\,where\,}
a^{(1)}=\left(\begin{smallmatrix} 0\\
0\\
a_{21}\\
\end{smallmatrix}\right).
$$

By analogy for $n=3$ we have
$$
[\sigma_1(q,3)D_3^\sharp(q)-q_3I] a^{(0)}=0,\quad
[\sigma_1(q,3)D_3^\sharp(q)-q_2I] a^{(1)}=0,
$$
$$
[\sigma_1(q,3)D_3^\sharp(q)-q_1I] a^{(2)}=0,
$$
where $a^{(0)}=(0,a_{10},a_{20},a_{30})^t,\,
a^{(1)}=(0,0,a_{21},a_{31})^t,\,a^{(2)}=(0,0,0,a_{32})^t.$
For general $n$ we get
$$
[\sigma_1(q,n)D_n^\sharp(q)-q_{n-k}I] a^{(k)}=0,\quad
a^{(k)}=(0,0,...,a_{k+1,k},...,a_{nk})^t,\quad 0\leq k<n.
$$
To prove Lemma it is sufficient to show that all solutions of the
latter equations
are trivial. We rewrite the
latter equations in the following forms:
\begin{equation}
\label{case3} [\sigma_1(q,n)-q_{n-k}(D_n^\sharp(q))^{-1}]
b^{(k)}=0,\quad 0\leq
k<n,\text{\,\,where\,\,}b^{(k)}=D_n^\sharp(q)a^{(k)}.
\end{equation}
Set $(k_n):=\sigma_1(q,n)-q_{n-k}(D_n^\sharp(q))^{-1}.$ The
equations (\ref{case3}) {\bf for} $n=2$ gives us
$$
\sigma_1(q,2)=\left(\begin{smallmatrix}
q&1+q&1\\
0&1&1\\
0&0&1
\end{smallmatrix}\right),\,\,
\left(\begin{smallmatrix}
0&1+q&1\\
0&1-q&1\\
0&0&1-q\\
\end{smallmatrix}\right)\left(\begin{smallmatrix} 0\\
b_{10}\\
b_{20}\\
\end{smallmatrix}\right)=0,\quad
\left(\begin{smallmatrix}
q-1&1+q&1\\
0&0&1\\
0&0&0\\
\end{smallmatrix}\right)\left(\begin{smallmatrix} 0\\
0\\
b_{21}\\
\end{smallmatrix}\right)=0.
$$
Hence $b^{(0)}=0$ if some of minors
$M^{12}_{12}(0_2),\,\,M^{01}_{12}(0_2)$ or $M^{02}_{12}(0_2)$ are
not zero. Further we get $b^{(1)}=0$ since $M^{0}_{2}(1_2)=1$. We
have $$M^{01}_{12}(0_2)=\left|\begin{smallmatrix}
1+q&1\\
1-q&1\\
\end{smallmatrix}\right|=M^{01}_{12}(F_{02}^s(q,1))=2q
$$
$$
M^{02}_{12}(0_2)=\left|\begin{smallmatrix}
1+q&1\\
0&1-q\\
\end{smallmatrix}\right|=M^{12}_{12}(F_{02}^s(q,1)),\,\,
M^{12}_{12}(0_2)=\left|\begin{smallmatrix}
1-q&1\\
0&1-q\\
\end{smallmatrix}\right|=M^{12}_{12}(F_{02}^s(q,1)).
$$
The {\bf suspected case} (see Remark) is
$\beta_{01}(q,2)=C_{2}^{1}(q)=1+q=0$  i.e. $q=-1$. We show later
that the representation is reducible in this case.
{\bf For} $n=3$ we have $D_3^\sharp(q)={\rm diag}(q^3,q,1,1)$ and
$$
(k_3)=\sigma_1(q,3)-q_k(D_3^\sharp(q))^{-1}$$
$$
=\left(\begin{smallmatrix}
1&1+q+q^2&1+q+q^2&1\\
0&1&1+q&1\\
0&0&1&1\\
0&0&0&1\\
\end{smallmatrix}\right)-q_k\left(\begin{smallmatrix}
q^3&0&0&0\\
0&q&0&0\\
0&0&1&0\\
0&0&0&1\\
\end{smallmatrix}\right)^{-1}=\left(\begin{smallmatrix}
1-q_k/q_3&1+q+q^2&1+q+q^2&1\\
0&1-q_k/q_2&1+q&1\\
0&0&1-q_k/q_1&1\\
0&0&0&1-q_k/q_0\\
\end{smallmatrix}\right),
$$
so the equations (\ref{case3}) for $n=3$ gives us
$$
\left(\begin{smallmatrix}
0&1+q+q^2&1+q+q^2&1\\
0&1-q^2&1+q&1\\
0&0&1-q^3&1\\
0&0&0&1-q^3\\
\end{smallmatrix}\right)\left(\begin{smallmatrix} 0\\
b_{10}\\
b_{20}\\
b_{30}\\
\end{smallmatrix}\right)=0,\,\,
\left(\begin{smallmatrix}
1-q^{-2}&1+q+q^2&1+q+q^2&1\\
0&0&1+q&1\\
0&0&1-q&1\\
0&0&0&1-q\\
\end{smallmatrix}\right)\left(\begin{smallmatrix} 0\\
0\\
b_{21}\\
b_{31}\\
\end{smallmatrix}\right)=0,
$$
$$
\left(\begin{smallmatrix}
1-q^{-3}&1+q+q^2&1+q+q^2&1\\
0&1-q^{-1}&1+q&1\\
0&0&0&1\\
0&0&0&0\\
\end{smallmatrix}\right)\left(\begin{smallmatrix} 0\\
0\\
0\\
b_{32}\\
\end{smallmatrix}\right)=0.
$$
We conclude that $b^{(0)}=0$ if
$M^{i_0i_1i_2}_{123}(0_3)=M^{i_0i_1i_2}_{123}(F_{03}^s(q,1))\not=0$
for some $0\leq i_0<i_1<i_2\leq 3$, further $b^{(1)}=0$ since
$M^{01}_{23}(1_3)=M^{01}_{23}(F_{13}^s(q,1))=q^2\not=0$, and
$b^{(2)}=0$ since
$M^{0}_{3}(2_3)=M^{0}_{3}(F_{23}^s(q,1))=1\not=0$. We have
\begin{equation}
M^{123}_{123}(0_3)=\left|\begin{smallmatrix}
1-q^2&1+q&1\\
0&1-q^3&1\\
0&0&1-q^3\\
\end{smallmatrix}\right|=(1-q^2)(1-q^3)^2,\quad
M^0_3(2)=1,
\end{equation}
\begin{equation}
M^{012}_{012}(0_3)=\left|\begin{smallmatrix}
1+q+q^2&1+q+q^2&1\\
1-q^2&1+q&1\\
0&1-q^3&1\\
\end{smallmatrix}\right|=2q^3(1+q+q^2),\,\,
M^{01}_{23}(1)=\left|\begin{smallmatrix}
1+q+q^2&1\\
1+q&1\\
\end{smallmatrix}\right|=q^2.
\end{equation}
The {\bf suspected case} (see Remark) are
$\beta_{01}(q,3)=C_{3}^{1}(q)=1+q+q^2=0$ and
$\beta_{12}(q,3)=C_{2}^{1}(q)=1+q=0$ i.e. $q^3=1$ and
$q^2=1,\,\,q\not=1$. Finally the suspected values are
$q=\alpha_1^{(3)}=\exp(2\pi i /3),\,\,q=\alpha_2^{(3)}=\exp(2\pi i
2/3)$ and $q=\alpha_1^{(2)}=\exp(2\pi i /2)=-1$ where
\begin{equation}
\alpha_k^{(s)}=\exp(2\pi i k/s),\,\,0\leq k \leq s,\,\,s=1,2,...
\end{equation}
We show later that the representation is reducible in this case.

Since $D_4^\sharp(q)={\rm diag}(q^6,q^3,q,1,1)$ and
$$
(k_4)=\sigma_1(q,4)-q_k(D_4^\sharp(q))^{-1}$$
$$
=\left(\begin{smallmatrix}
1&(1+q)(1+q^2)&(1+q^2)(1+q+q^2)&(1+q)(1+q^2)&1\\
0&1&1+q+q^2&1+q+q^2&1\\
0&0&1&1+q&1\\
0&0&0&1&1\\
0&0&0&0&1\\
\end{smallmatrix}\right)-q_k\left(\begin{smallmatrix}
q^6&0&0&0&0\\
0&q^3&0&0&0\\
0&0&q&0&0\\
0&0&0&1&0\\
0&0&0&0&1\\
\end{smallmatrix}\right)^{-1}=
$$
$$
\left(\begin{smallmatrix}
1-q_k/q_4&(1+q)(1+q^2)&(1+q^2)(1+q+q^2)&(1+q)(1+q^2)&1\\
0&1-q_k/q_3&1+q+q^2&1+q+q^2&1\\
0&0&1-q_k/q_2&1+q&1\\
0&0&0&1-q_k/q_1&1\\
0&0&0&0&1-q_k/q_0\\
\end{smallmatrix}\right),
$$
the equations (\ref{case3}) {\bf for} $n=4$ gives us
$$
\left(\begin{smallmatrix}
0&(1+q)(1+q^2)&(1+q^2)(1+q+q^2)&(1+q)(1+q^2)&1\\
0&1-q^3&1+q+q^2&1+q+q^2&1\\
0&0&1-q^5&1+q&1\\
0&0&0&1-q^6&1\\
0&0&0&0&1-q^6\\
\end{smallmatrix}\right)\left(\begin{smallmatrix} 0\\
b_{10}\\
b_{20}\\
b_{30}\\
b_{40}\\
\end{smallmatrix}\right)=0,
$$
$$
\left(\begin{smallmatrix}
1-q^{-3}&(1+q)(1+q^2)&(1+q^2)(1+q+q^2)&(1+q)(1+q^2)&1\\
0&0&1+q+q^2&1+q+q^2&1\\
0&0&1-q^2&1+q&1\\
0&0&0&1-q^3&1\\
0&0&0&0&1-q^3\\
\end{smallmatrix}\right)\left(\begin{smallmatrix} 0\\
0\\
b_{21}\\
b_{31}\\
b_{41}\\
\end{smallmatrix}\right)=0,
$$
$$
\left(\begin{smallmatrix}
1-q^{-5}&(1+q)(1+q^2)&(1+q^2)(1+q+q^2)&(1+q)(1+q^2)&1\\
0&1-q^{-2}&1+q+q^2&1+q+q^2&1\\
0&0&0&1+q&1\\
0&0&0&1-q&1\\
0&0&0&0&1-q\\
\end{smallmatrix}\right)\left(\begin{smallmatrix} 0\\
0\\
0\\
b_{32}\\
b_{42}\\
\end{smallmatrix}\right)=0,
$$
$$
\left(\begin{smallmatrix}
1-q^{-6}&(1+q)(1+q^2)&(1+q^2)(1+q+q^2)&(1+q)(1+q^2)&1\\
0&1-q^{-3}&1+q+q^2&1+q+q^2&1\\
0&0&1-q^{-1}&1+q&1\\
0&0&0&0&1\\
0&0&0&0&0\\
\end{smallmatrix}\right)\left(\begin{smallmatrix} 0\\
0\\
0\\
0\\
b_{43}\\
\end{smallmatrix}\right)=0.
$$
The {\bf suspected case}  are
$\beta_{01}(q,4)=C_{4}^{1}(q)=1+q+q^2+q^3=0,\,\,\beta_{12}(q,4)=C_{3}^{1}(q)=1+q+q^2=0
$ and $\beta_{23}(q,4)=C_{2}^{1}(q)=1+q=0$ i.e. $q^s=1,\,\,2\leq
s\leq 4$ and $q\not=1$. Finally the suspected values are
$q=\alpha_k^{(s)},\,\,1\leq k<s\leq 4$.

For the general case $n\in{\mathbb N}$ the {\bf suspected case}
are $q^k=1,\,\,2\leq k\leq n$ and $q\not=1$ i.e.
$q=\alpha_k^{(s)},\,\,1\leq k<s\leq n$.

 \qed\end{pf}
\begin{lem}
\label{3)red} The representation
$$
\sigma_1\mapsto\sigma_1^D(q,n):=\sigma_1(q,n)D_n^\sharp(q)\quad
\sigma_2\mapsto \sigma_2^D(q,n):=D_n(q)\sigma_2(q,n),
$$
is irreducible if and only if $(n)_q=1+q+...+q^{n-1}\not=0.$
\end{lem}
\begin{pf}
 For $n=2$ and $(2)_q=1+q=0$ we have (see (\ref{si_1(q)}) and
(\ref{D_n(q)}))
 $$
\sigma_1^D(q,2) =\left(\begin{smallmatrix}
q&1+q&1\\
0&1&1\\
0&0&1
\end{smallmatrix}\right)=
\left(\begin{smallmatrix}
-1&0&1\\
0&1&1\\
0&0&1
\end{smallmatrix}\right)
,\,\, 
\sigma_2^D(q,2)=\left(\begin{smallmatrix}
1&0&0\\
-1&1&0\\
1&-(1+q)&q\\
\end{smallmatrix}\right)=\left(\begin{smallmatrix}
1&0&0\\
-1&1&0\\
1&0&-1\\
\end{smallmatrix}\right).
 $$
The vector $e_2=(0,1,0)$ is the eigenvector for $\sigma_1^D(q,2)
$ and $
\sigma_2^D(q,2)$ with eigenvalue equal to $1$:
$$
\sigma_1^D(q,2)
e_2=e_2,\quad
\sigma_2^D(q,2)e_2=e_2.
$$
hence the subspace $V_2=\{te_2=(0,t,0)\mid t\in {\mathbb C}\}$ is
nontrivial invariant subspace for $\sigma_1^D(q,2)$ and
$\sigma_2^D(q,2)$.

Let $1+q\not=0$. If we set $\sigma_k=\sigma_k^D(q,2),\,\,k=1,2$ we
get
 $$
\sigma_1-I=\left(\begin{smallmatrix}
q-1&1+q&1\\
0&0&1\\
0&0&0
\end{smallmatrix}\right),\quad
\sigma_2-I=\left(\begin{smallmatrix}
0&0&0\\
-1&0&0\\
1&-(1+q)&q-1\\
\end{smallmatrix}\right).
$$
Since
$$
(\sigma_1-I)^2=\left(\begin{smallmatrix}
(q-1)^2&(q-1)(1+q)&2q\\
0&0&0\\
0&0&0
\end{smallmatrix}\right),\quad
(\sigma_2-I)^2=\left(\begin{smallmatrix}
0&0&0\\
0&0&0\\
2q&-(1+q)(q-1)&(q-1)^2\\
\end{smallmatrix}\right)
$$
we conclude that
\begin{align}
\label{char-1}
 A_1:=(\sigma_1-I)^2-(q-1)(\sigma_1-I)=
\left(\begin{smallmatrix}
0&0&1+q\\
0&0&1-q\\
0&0&0
\end{smallmatrix}\right),\\
A_2:= (\sigma_2-I)^2-(q-1)(\sigma_2-I)=
\left(\begin{smallmatrix}
0&0&0\\
1-q&0&0\\
1+q&0&0
\end{smallmatrix}\right).
\end{align}

Further we get
$$
(1+q)^{-1}A_1A_2=\left(\begin{smallmatrix}
1+q&0&0\\
1-q&0&0\\
0&0&0
\end{smallmatrix}\right),\quad
(1+q)^{-1}A_2A_1=\left(\begin{smallmatrix}
0&0&0\\
0&0&1-q\\
0&0&1+q
\end{smallmatrix}\right).
$$
Finally we have 5 matrix
$$
a=\left(\begin{smallmatrix}
0&0&1+q\\
0&0&1-q\\
0&0&0
\end{smallmatrix}\right),\,b=\left(\begin{smallmatrix}
0&0&0\\
1-q&0&0\\
1+q&0&0
\end{smallmatrix}\right),\,
c=\left(\begin{smallmatrix}
1+q&0&0\\
1-q&0&0\\
0&0&0
\end{smallmatrix}\right),\,
d=\left(\begin{smallmatrix}
0&0&0\\
0&0&1-q\\
0&0&1+q
\end{smallmatrix}\right),\,\,
qS(q)=\left(\begin{smallmatrix}
0&0 &q\\
0&-q&0\\
1&0 &0
\end{smallmatrix}\right).
$$
We get
$$
(1+q)^{-1}(a-d)=\left(\begin{smallmatrix}
0&0&1\\
0&0&0\\
0&0&-1
\end{smallmatrix}\right),\,
(1+q)^{-1}(c-b)=\left(\begin{smallmatrix}
1&0&0\\
0&0&0\\
-1&0&0
\end{smallmatrix}\right),
$$
hence we have
$$
qS(q)(1+q)^{-1}(c-b)=\left(\begin{smallmatrix}
0&0 &q\\
0&-q&0\\
1&0 &0
\end{smallmatrix}\right)
\left(\begin{smallmatrix}
1&0&0\\
0&0&0\\
-1&0&0
\end{smallmatrix}\right)=
\left(\begin{smallmatrix}
-q&0&0\\
0&0&0\\
1&0&0
\end{smallmatrix}\right),
$$
and
$$
qS(q)(1+q)^{-1}(a-d) =\left(\begin{smallmatrix}
0&0 &q\\
0&-q&0\\
1&0 &0
\end{smallmatrix}\right)\left(\begin{smallmatrix}
0&0&1\\
0&0&0\\
0&0&-1
\end{smallmatrix}\right)
= \left(\begin{smallmatrix}
0&0&-q\\
0&0&0\\
0&0&1
\end{smallmatrix}\right).
$$
Since $q\not=1$ we can obtain
$$
\left(\begin{smallmatrix}
1&0&0\\
0&0&0\\
0&0&0
\end{smallmatrix}\right),\,\left(\begin{smallmatrix}
0&0&0\\
0&0&0\\
1&0&0
\end{smallmatrix}\right),\,\,
\left(\begin{smallmatrix}
0&0&1\\
0&0&0\\
0&0&0
\end{smallmatrix}\right),\,\,
\left(\begin{smallmatrix}
0&0&0\\
0&0&0\\
0&0&1
\end{smallmatrix}\right),\,\,
\left(\begin{smallmatrix}
0&0&0\\
0&1&0\\
0&0&0
\end{smallmatrix}\right)=I-\left(\begin{smallmatrix}
1&0&0\\
0&0&0\\
0&0&1
\end{smallmatrix}\right).
$$
Using again
 $$
\sigma_1-I=\left(\begin{smallmatrix}
q-1&1+q&1\\
0&0&1\\
0&0&0
\end{smallmatrix}\right),\quad
\sigma_2-I=\left(\begin{smallmatrix}
0&0&0\\
-1&0&0\\
1&-(1+q)&q-1\\
\end{smallmatrix}\right),
$$
we conclude that we can obtain the following matrices
$$
\beta:=\left(\begin{smallmatrix}
0&1+q&0\\
0&0&1\\
0&0&0
\end{smallmatrix}\right),\quad
-\beta^\sharp:=\left(\begin{smallmatrix}
0&0&0\\
-1&0&0\\
0&-(1+q)&0\\
\end{smallmatrix}\right).
$$
By Remark \ref{(T,T*)=Mat} two latter matrices generate ${\rm
Mat}(3,{\mathbb C})$ if $(2)_q=1+q\not=0$. Using the latter and
the previous matrices we can obtain
$$
\left(\begin{smallmatrix}
0&1&0\\
0&0&0\\
0&0&0
\end{smallmatrix}\right),\,\left(\begin{smallmatrix}
0&0&0\\
0&0&1\\
0&0&0
\end{smallmatrix}\right),\,\,
\left(\begin{smallmatrix}
0&0&0\\
1&0&0\\
0&0&0
\end{smallmatrix}\right),\,\,
\left(\begin{smallmatrix}
0&0&0\\
0&0&0\\
0&1&0
\end{smallmatrix}\right).
$$
Indeed we have
$$
\left(\begin{smallmatrix}
1&0&0\\
0&0&0\\
0&0&0
\end{smallmatrix}\right)\left(\begin{smallmatrix}
0&1+q&0\\
0&0&1\\
0&0&0
\end{smallmatrix}\right)=\left(\begin{smallmatrix}
0&1+q&0\\
0&0&0\\
0&0&0
\end{smallmatrix}\right),\,\,
\left(\begin{smallmatrix}
0&0&1\\
0&0&0\\
0&0&0
\end{smallmatrix}\right)\left(\begin{smallmatrix}
0&0&0\\
1&0&0\\
0&1+q&0
\end{smallmatrix}\right)=\left(\begin{smallmatrix}
0&0&0\\
0&0&0\\
0&1+q&0
\end{smallmatrix}\right).
$$

 Finally we conclude that we can obtain all matrix
units $E_{kn},\,0\leq k\leq 2$ so the algebra, generated by two
matrices $\sigma_1^D(q,2)$ and $\sigma_2^D(q,2)$ coincides with
the algebra ${\rm Mat}(3,{\mathbb C})$. So our representation is
irreducible for $n=2$ when $(2)_q=1+q\not=0$ by the Remark
\ref{mat()-irr}.

Let $n=3$ and $(3)_q=1+q+q^2=0$. Then $q=\alpha_k^{(3)}:=\exp(2\pi
i k/3),\,\,k=1,2$ and we have
$$
\sigma_1^D(q,3)=\left(\begin{smallmatrix}
q^3&q(1+q+q^2)&1+q+q^2&1\\
0&q&1+q&1\\
0&0&1&1\\
0&0&0&1\\
\end{smallmatrix}\right)=
\left(\begin{smallmatrix}
1&0&0&1\\
0&\alpha_3^{(k)}&1+\alpha_3^{(k)}&1\\
0&0&1&1\\
0&0&0&1\\
\end{smallmatrix}\right),
$$
$$
\sigma_2^D(q,3)=\left(\begin{smallmatrix}
1&0&0&0\\
-1&1&0&0\\
1&-(1+q)&q&0\\
-1&(1+q+q^2)&-q(1+q+q^2)&q^3\\
\end{smallmatrix}\right)=
\left(\begin{smallmatrix}
1&0&0&0\\
-1&1&0&0\\
1&-(1+\alpha_3^{(k)})&\alpha_3^{(k)}&0\\
-1&0&0&1\\
\end{smallmatrix}\right).
$$
Obviously, the subspace $V_3=\{(0,t_1,t_2,0)\mid (t_1,t_2)\in
{\mathbb C}^2\}$ is invariant subspace for $\sigma_1^D(q,3)$ and
$\sigma_2^D(q,3)$.

Let $(3)_q=1+q+q^2\not=0$. If we set
$\sigma_k=\sigma_k^D(q,3),\,\,k=1,2$ we get
$$
\sigma_1-1=\left(\begin{smallmatrix}
q^3-1&q(1+q+q^2)&1+q+q^2&1\\
0&q-1&1+q&1\\
0&0&0&1\\
0&0&0&0\\
\end{smallmatrix}\right),\,\,\,
\sigma_2-1=\left(\begin{smallmatrix}
0&0&0&0\\
-1&0&0&0\\
1&-(1+q)&q-1&0\\
-1&(1+q+q^2)&-q(1+q+q^2)&q^3-1\\
\end{smallmatrix}\right).
$$
To generalize expressions (\ref{char-1}) we note the following
\begin{rem} Let $P_A$ be the characteristic polynomial of the matrix
$A$, in the space ${\mathbb C}^{n+1}$ with the spectra
$(\lambda_k)_{k=0}^n$ i.e.
\begin{equation}\label{char-1.n}
P_A(x)=
\prod_{k=0}^n(x-\lambda_k),\text{\quad then\quad}
A_1=P_{\sigma_1}(\sigma_1)(\sigma_1-I)^{-1}.
\end{equation}
\end{rem}
Indeed we have ${\rm Sp}\,{\sigma_1}=\{q,1,1\}$, hence
$$
P_{\sigma_1}(\sigma_1)=(\sigma_1-qI)(\sigma_1-I)(\sigma_1-I)
$$
and
$$
A_1:=(\sigma_1-I)^2-(q-1)(\sigma_1-I)=(\sigma_1-qI)(\sigma_1-I)=
P_{\sigma_1}(\sigma_1)(\sigma_1-I)^{-1}.
$$
We would like to find the expression for
$P_{\sigma_1}(\sigma_1)(\sigma_1-I)^{-1}$ (when  $n=3$) in the
following form (see (\ref{char-1}))
$$
P_{\sigma_1}(\sigma_1)(\sigma_1-I)^{-1}=\left(\begin{smallmatrix} 0&0&0&x\\
0&0&0&(1+q)t\\
0&0&0&(1-q)t\\
0&0&0&0\\
\end{smallmatrix}\right).
$$
To find $x$ and $t$ we use the identity
$(\sigma_1-I)P_{\sigma_1}(\sigma_1)(\sigma_1-I)^{-1}=0$ i.e.
$$
\left(\begin{smallmatrix}
q^3-1&q(1+q+q^2)&1+q+q^2&1\\
0&q-1&1+q&1\\
0&0&0&1\\
0&0&0&0\\
\end{smallmatrix}\right)
\left(\begin{smallmatrix} 0&0&0&x\\
0&0&0&(1+q)t\\
0&0&0&(1-q)t\\
0&0&0&0\\
\end{smallmatrix}\right)=0.
$$
We have
$$
0=(q^3-1)x+[q(1+q+q^2)(1+q)+
(1+q+q^2)(1-q)]t=(1+q+q^2)\times
$$
$$
\{(q-1)x+[q(1+q)+(1-q)]t\}= (1+q+q^2)[(q-1)x+(1+q^2)t]
$$
hence $x_3:=x=1+q^2,\,\,t_3:=t=1-q$. Before we have calculated
$x_2=1+q,\,\,t_2=1-q$. Finally we have
$$
P_{\sigma_1}(\sigma_1)(\sigma_1-I)^{-1}=\left(\begin{smallmatrix}
0&0&0&1+q^2\\
0&0&0&1-q^2\\
0&0&0&(1-q)^2\\
0&0&0&0\\
\end{smallmatrix}\right).
$$
Let $n=4$ and $(4)_q=1+q+q^2+q^3=0$. Then
$q=\alpha_k^{(4)}:=\exp(2\pi i k/4),\,\,k=1,2,3$ and we have
$$
\left(\begin{smallmatrix}
q^6&q^3(1+q)(1+q^2)&q(1+q^2)(1+q+q^2)&(1+q)(1+q^2)&1\\
0&q^3&q(1+q+q^2)&1+q+q^2&1\\
0&0&q&1+q&1\\
0&0&0&1&1\\
0&0&0&0&1\\
\end{smallmatrix}\right)=
\left(\begin{smallmatrix}
1&0&0&0&1\\
0&q^3&q(1+q+q^2)&1+q+q^2&1\\
0&0&q&1+q&1\\
0&0&0&1&1\\
0&0&0&0&1\\
\end{smallmatrix}\right),
$$
$$
\left(\begin{smallmatrix}
1&0&0&0&0\\
-1&1&0&0&0\\
1&-(1+q)&q&0&0\\
-1&1+q+q^2&-q(1+q+q^2)&q^3&0\\
1&-(1+q)(1+q^2)&q(1+q)(1+q^2)&-q^3(1+q^2)(1+q+q^2)&q^6\\
\end{smallmatrix}\right)=
\left(\begin{smallmatrix}
1&0&0&0&0\\
-1&1&0&0&0\\
1&-(1+q)&q&0&0\\
-1&1+q+q^2&-q(1+q+q^2)&q^3&0\\
1&0&0&0&1\\
\end{smallmatrix}\right).
$$
Obviously, the subspace $V_4=\{(0,t_1,t_2,t_3,0)\mid
(t_1,t_2,t_3)\in {\mathbb C}^3\}$ is invariant subspace for
$\sigma_1^D(q,4)$ and $\sigma_2^D(q,4)$. Let
$(4)_q=1+q+q^2+q^3\not=0$. We would like to find the expression
for $P_{\sigma_1}(\sigma_1)(\sigma_1-I)^{-1}$ and $n=4$ in the
following form
$$
P_{\sigma_1}(\sigma_1)(\sigma_1-I)^{-1}=
\left(\begin{smallmatrix} 0&0&0&0&x\\
0&0&0&0&(1+q^2)t\\
0&0&0&0&(1-q^2)t\\
0&0&0&0&(1-q)^2t\\
0&0&0&0&0\\
\end{smallmatrix}\right).
$$
As before we get
$$
\left(\begin{smallmatrix}
q^6-1&q^3(1+q)(1+q^2)&q(1+q^2)(1+q+q^2)&(1+q)(1+q^2)&1\\
0&q^3-1&q(1+q+q^2)&1+q+q^2&1\\
0&0&q-1&1+q&1\\
0&0&0&0&1\\
0&0&0&0&0\\
\end{smallmatrix}\right)\left(\begin{smallmatrix} 0&0&0&0&x\\
0&0&0&0&(1+q^2)t\\
0&0&0&0&(1-q^2)t\\
0&0&0&0&(1-q)^2t\\
0&0&0&0&0\\
\end{smallmatrix}\right)=0,
$$
hence
$$
(q^6-1)x+(1+q^2)[q^3(1+q)(1+q^2)+q(1+q+q^2)(1-q^2)+
(1+q)(1-q)^2]t=
$$
$$
(q^6-1)x+(1+q^2)[1+2q^3+q^6]t=(q^3+1)(q^3-1)x+(1+q^2)(q^3+1)^2.
$$
Finally we conclude that
$$
x_4:=x=(1+q^2)(1+q^3),\,\,t_4:=t=(1-q^3).
$$
$$
\left(\begin{smallmatrix}
0&0&1+q\\
0&0&1-q\\
0&0&0
\end{smallmatrix}\right),\,\,\left(\begin{smallmatrix}
0&0&0&1+q^2\\
0&0&0&1-q^2\\
0&0&0&(1-q)^2\\
0&0&0&0\\
\end{smallmatrix}\right),\,\,\left(\begin{smallmatrix} 0&0&0&0&(1+q^2)(1+q^3)\\
0&0&0&0&(1+q^2)(1-q^3)\\
0&0&0&0&(1-q^2)(1-q^3)\\
0&0&0&0&(1-q)^2(1-q^3)\\
0&0&0&0&0\\
\end{smallmatrix}\right).
$$

 For $n=3$ and $q=-1$
$$
\sigma_1^D(q,3)= \left(\begin{smallmatrix}
-1&-1&1&1\\
0&-1&0&1\\
0&0&1&1\\
0&0&0&1\\
\end{smallmatrix}\right),\,\,
\sigma_2^D(q,3)= \left(\begin{smallmatrix}
1&0&0&0\\
-1&1&0&0\\
1&0&-1&0\\
-1&1&-1&-1\\
\end{smallmatrix}\right).
$$
We can prove as before that representation is irreducible. For
general $n$  the proof is similar. \qed\end{pf}

{\bf Case 4).} We prove the following lemma (see case 3)).
\begin{lem}
\label{l.case4} Let an operator $A\in {\rm Mat(n+1,{\mathbb C})}$
commute  with $\sigma_1(q)D_n^\sharp(q)\Lambda_n$ where
$\Lambda_n={\rm diag}(\lambda_0,\lambda_1,...,\lambda_n)$ with
$\lambda_r\lambda_{n-r}=c,\,\,0\leq r\leq n$ then $A$ is also
upper triangular, i.e.
\begin{equation}
\label{Acase4} A=\sum_{0\leq k\leq m\leq n}a_{km}E_{km}.
\end{equation}
if for any $0\leq r\leq \left[\frac{n}{2}\right]$ there exists
$0\leq i_0<i_i<...<i_r\leq n $ such that 
$$
 M^{i_0i_i...i_{n-r-1}}_{r+1r+2...n} (F_{rn}^s(q,\lambda))\not=0\,\,
\text{\,\,where\,\,}\nu^{(r)}=(\nu^{(r)}_k),\,\,
\nu^{(r)}_k=\lambda_r/\lambda_{n-k}.\\
$$
\end{lem}
\begin{pf}
For $n=2$ we have (see (\ref{si_1(q)}) and (\ref{D_n(q)}))
$$
 \sigma_1(q,2)D_2^\sharp(q)\Lambda= \left(\begin{smallmatrix}
q&1+q&1\\
0&1&1\\
0&0&1
\end{smallmatrix}\right)\left(\begin{smallmatrix}
\lambda_0&0&0\\
0&\lambda_1&0\\
0&0&\lambda_2
\end{smallmatrix}\right)=
\left(\begin{smallmatrix}
q\lambda_0&(1+q)\lambda_1&\lambda_2\\
0&\lambda_1&\lambda_2\\
0&0&\lambda_1
\end{smallmatrix}\right),
$$
$$
\left(\begin{smallmatrix}
q\lambda_0&(1+q)\lambda_1&\lambda_2\\
0&\lambda_1&\lambda_2\\
0&0&\lambda_1
\end{smallmatrix}\right)
\left(\begin{smallmatrix}
a_{00}&a_{01}&a_{02}\\
a_{10}&a_{11}&a_{12}\\
a_{20}&a_{21}&a_{22}
\end{smallmatrix}\right)=
$$
$$
\left(\begin{smallmatrix}
q\lambda_0a_{00}+(1+q)\lambda_1a_{10}+\lambda_2a_{20}&
q\lambda_0a_{01}+(1+q)\lambda_1a_{11}+\lambda_2a_{21}
&q\lambda_0a_{21}+(1+q)\lambda_1a_{12}+\lambda_2a_{22}\\
\lambda_1a_{10}+\lambda_2a_{20}&\lambda_1a_{11}+\lambda_2a_{21}&\lambda_1a_{12}+\lambda_2a_{22}\\
\lambda_2a_{20}&\lambda_2a_{21}&\lambda_2a_{22}
\end{smallmatrix}\right),
$$
$$
\left(\begin{smallmatrix}
a_{00}&a_{01}&a_{02}\\
a_{10}&a_{11}&a_{12}\\
a_{20}&a_{21}&a_{22}
\end{smallmatrix}\right)\left(\begin{smallmatrix}
q&1+q&1\\
0&1&1\\
0&0&1
\end{smallmatrix}\right)
\left(\begin{smallmatrix}
\lambda_0&0&0\\
0&\lambda_1&0\\
0&0&\lambda_2
\end{smallmatrix}\right)
= \left(\begin{smallmatrix}
a_{00}q\lambda_0\,&\,[a_{00}(1+q)+a_{01}]\lambda_1&[a_{00}+a_{01}+a_{02}]\lambda_2\\
a_{10}q\lambda_0\,&\,[a_{10}(1+q)+a_{11}]\lambda_1&[a_{10}+a_{11}+a_{12}]\lambda_2\\
a_{20}q\lambda_0\,&\,[a_{20}(1+q)+a_{21}]\lambda_1&[a_{20}+a_{21}+a_{22}]\lambda_2
\end{smallmatrix}\right).
$$
If we compare the first columns we get
$$
\left\{\begin{smallmatrix}
(1+q)\lambda_1a_{10}+\lambda_2a_{20}=0\\
(\lambda_1-q\lambda_0)a_{10}+\lambda_2a_{20}=0\\
(\lambda_2-q\lambda_0)a_{20}=0.\\
\end{smallmatrix}\right.
 \text{\quad or \quad}
[\sigma_1(q,2)D_2^\sharp(q)\Lambda_2-q_2\lambda_0I]
a^{(0)}=0,\text{\,where\,}
a^{(0)}=\left(\begin{smallmatrix} 0\\
a_{10}\\
a_{20}\\
\end{smallmatrix}\right).
$$
Let $ a^{(0)}=0$. If we compare the second columns we get
$$
\left\{\begin{smallmatrix}
\lambda_2a_{21}=0\\
(\lambda_2-\lambda_1)a_{21}=0\\
\end{smallmatrix}\right.
\text{\quad or \quad}
[\sigma_1(q,2)D_2^\sharp(q)\Lambda_2-q_1\lambda_1I]
a^{(1)}=0,\text{\,where\,}
a^{(1)}=\left(\begin{smallmatrix} 0\\
0\\
a_{21}\\
\end{smallmatrix}\right).
$$
By analogy for $n=3$ we have
$$
[\sigma_1(q,3)D_3^\sharp(q)\Lambda_3-q_3\lambda_0I]
a^{(0)}=0,\quad
[\sigma_1(q,3)D_3^\sharp(q)\Lambda_3-q_2\lambda_1I] a^{(0)}=0,
$$
$$
[\sigma_1(q,3)D_3^\sharp(q)\Lambda_3-q_1\lambda_2I] a^{(0)}=0,
$$
where $a^{(0)}=(0,a_{10},a_{20},a_{30})^t,\,
a^{(1)}=(0,0,a_{21},a_{31})^t,\,a^{(2)}=(0,0,0,a_{32})^t.$

For general $n$ we get
$$
[\sigma_1(q,n)D_n^\sharp(q)\Lambda_n-q_{n-k}\lambda_kI]
a^{(k)}=0,\quad a^{(k)}=(0,0,...,a_{k+1,k},...,a_{nk})^t,\quad
0\leq k<n.
$$
To prove Lemma it is sufficient to show that all solutions of the
latter equations are trivial.

{\bf Let us set}
$(k_n):=\sigma_1^{\Lambda,k}(q,n):=\sigma_1(q,n)-q_{n-k}\lambda_k(D_n^\sharp(q)
\Lambda_n)^{-1}.$ We rewrite the latter equations in the following
forms:
\begin{equation}
\label{case4}
\sigma_1^{\Lambda,k}(q,n)b^{(k)}=0,\quad 0\leq
k<n,\text{\,\,where\,\,}b^{(k)}=D_n^\sharp(q)\Lambda_na^{(k)}.
\end{equation}
If we denote
\begin{equation}
\label{F_kn(q,lam)} F_{k,n}(q,\lambda)=
[\sigma_1(q,n)-q_{n-k}\lambda_k(D_n^\sharp(q)\Lambda_n)^{-1}]^s
\end{equation}
we get by Lemma \ref{l.Pas.exp}
$$
F_{k,n}(q,\lambda)=
[\sigma_1(q,n)-q_{n-k}\lambda_k(D_n^\sharp(q)\Lambda_n)^{-1}]^s=
$$
$$
=\exp_{(q)}\left(\sum_{r=0}^{n-1} (r+1)_qE_{rr+1}\right)-
q_{n-k}\lambda_k(D_n(q)\Lambda_n^\sharp)^{-1}
$$
The equations (\ref{case3}) for $n=2$ gives us (we set
$\nu^k_m=\lambda_k/\lambda_m$) $(0_2)b^{(0)}=0$ and
$(1_2)b^{(1)}=0$ or
$$
\left(\begin{smallmatrix}
q&1+q&1\\
0&1&1\\
0&0&1
\end{smallmatrix}\right),\,\,
\left(\begin{smallmatrix}
0&1+q&1\\
0&1-q\nu_1^0&1\\
0&0&1-q\nu_2^0
\end{smallmatrix}\right)\left(\begin{smallmatrix} 0\\
b_{10}\\
b_{20}\\
\end{smallmatrix}\right)=0,\quad
\left(\begin{smallmatrix}
1-q\nu_0^1&1+q&1\\
0&0&1\\
0&0&1-\nu_2^1\\
\end{smallmatrix}\right)\left(\begin{smallmatrix} 0\\
0\\
b_{21}\\
\end{smallmatrix}\right)=0.
$$
Hence $b^{(0)}=0$ if $M^{i_0i_1}_{12}(0_2) \not=0$ for some $0\leq
i_0<i_1\leq 2$ and $b^{(1)}=0$ since $M^{0}_{2}(1_2)=1$.
We have
$$
M^{01}_{12}(0_2)= \left|\begin{smallmatrix}
1+q&1\\
1-q\nu_1^0&1\\
\end{smallmatrix}\right|,\,\,
 M^{02}_{12}(0_2)=\left|\begin{smallmatrix}
1+q&1\\
0&1-q\nu_2^0\\
\end{smallmatrix}\right|,\,\,
  M^{12}_{12}(0_2)=\left|\begin{smallmatrix}
1-q\nu_1^0&1\\
0&1-q\nu_2^0\\
\end{smallmatrix}\right|,\,\, M^{0}_{2}(1_2)=1,
$$
hence
$$
M^{01}_{12}(0_2)= M^{01}_{12}(F_{0,2}^s(q,\lambda)),\,\,
M^{02}_{12}(0_2)= M^{02}_{12}(F_{0,2}^s(q,\lambda)),\,\,
$$
$$
M^{12}_{12}(0_2)=M^{12}_{12}(F_{0,2}^s(q,\lambda)),\,\,
M^{0}_{2}(0_2)=M^{0}_{2}(F_{1,2}^s(q,\lambda)),
$$
where
$\nu^{(r)}=(\nu^{(r)}_k)_{k=0}^2,\,\,\nu^{(r)}_k=\lambda_r/\lambda_{2-k}$
and $0\leq r\leq [\frac{2}{2}]=1$. We have
$$
D_2(q,\nu):=M^{01}_{12}(0_2)=\left|\begin{smallmatrix}
1+q&1\\
1-q\nu_1&1
\end{smallmatrix}\right|=q(1+\nu_1)=q(\lambda_1+\lambda_0)/\lambda_1.
$$
$$
M^{2}_{2}(0_2)=1-q\nu_2^0=(\lambda_2-q\lambda_0)\lambda_2.
$$


We see that $M^{i_0i_1}_{12}(0_2) =0$ for all $0\leq i_0<i_1\leq
2$ if and only if $M^{01}_{12}(0_2)=0$ and $M^{2}_{2}(0_2)$ i.e.
$\lambda_0+\lambda_1=0$ and $\lambda_2-q\lambda_0=0$.

\begin{rem} We note that conditions $\lambda_0+\lambda_1=0$ and
$\lambda_2-q\lambda_0=0$ contradicts with conditions
$\lambda_r\lambda_{n-r}=c,\,\,0\leq r\leq n$ (see (\ref{Rep(q)'}),
for $n=2$ we have $\lambda_0\lambda_2=\lambda_1^2$). Indeed
otherwise we have $\lambda_2=q\lambda_0$ and
$\lambda_1=-\lambda_0$ hence $\lambda_1^2=\lambda_0^2$ and
$\lambda_0\lambda_2=q\lambda_0^2$ so
$\lambda_1^2\not=\lambda_0\lambda_2$ if $q\not=1.$
\end{rem}
In the general case $n\in {\mathbb N}$ we should calculate the
following determinant:
\begin{equation}
\label{det-q} D_n(q,\nu):= M^{01...n-1}_{23...n}\left[
\sigma_1(q,n)-q_n\lambda_0(D_n^\sharp(q)\Lambda_n)^{-1}\right].
\end{equation}
Since $D_3^\sharp(q)={\rm diag}(q^3,q,1,1)$ and
$$(k_3)=\sigma_1(q,3)-q_k(D_3^\sharp(q)\Lambda)^{-1}$$
$$
=\left(\begin{smallmatrix}
1&1+q+q^2&1+q+q^2&1\\
0&1&1+q&1\\
0&0&1&1\\
0&0&0&1\\
\end{smallmatrix}\right)-q_{3-k}\lambda_k\left(\begin{smallmatrix}
q^3\lambda_0&0&0&0\\
0&q\lambda_1&0&0\\
0&0&\lambda_2&0\\
0&0&0&\lambda_3\\
\end{smallmatrix}\right)^{-1}=
$$
$$
\left(\begin{smallmatrix}
1-q_{3-k}/q_3\nu^k_{0}&1+q+q^2&1+q+q^2&1\\
0&1-q_{3-k}/q_2\nu^k_{1}&1+q&1\\
0&0&1-q_{3-k}/q_1\nu^k_{2}&1\\
0&0&0&1-q_{3-k}/q_0\nu^k_{3}\\
\end{smallmatrix}\right),
$$
the equations (\ref{case3}) {\bf for} $n=3$ gives us
$$
\left(\begin{smallmatrix}
0&1+q+q^2&1+q+q^2&1\\
0&1-q^2\nu^{0}_{1}&1+q&1\\
0&0&1-q^3\nu^{0}_{2}&1\\
0&0&0&1-q^3\nu^{0}_{3}\\
\end{smallmatrix}\right)\left(\begin{smallmatrix} 0\\
b_{10}\\
b_{20}\\
b_{30}\\
\end{smallmatrix}\right)=0,\,\,
\left(\begin{smallmatrix}
1-q^{-2}\nu^{1}_{0}&1+q+q^2&1+q+q^2&1\\
0&0&1+q&1\\
0&0&1-q\nu^{1}_{2}&1\\
0&0&0&1-q\nu^{1}_{3}\\
\end{smallmatrix}\right)\left(\begin{smallmatrix} 0\\
0\\
b_{21}\\
b_{31}\\
\end{smallmatrix}\right)=0,
$$
$$
\left(\begin{smallmatrix}
1-q^{-3}\nu^{2}_{0}&1+q+q^2&1+q+q^2&1\\
0&1-q^{-1}\nu^{2}_{1}&1+q&1\\
0&0&0&1\\
0&0&0&1-\nu^{2}_{3}\\
\end{smallmatrix}\right)\left(\begin{smallmatrix} 0\\
0\\
0\\
b_{32}\\
\end{smallmatrix}\right)=0.
$$
We have
$$
M^{123}_{123}(0_3)=\left|\begin{smallmatrix}
1-q^2\nu^{0}_{1}&1+q&1\\
0&1-q^3\nu^{0}_{2}&1\\
0&0&1-q^3\nu^{0}_{3}\\
\end{smallmatrix}\right|,\,
M^{012}_{012}(0_3)=\left|\begin{smallmatrix}
1+q+q^2&1+q+q^2&1\\
1-q^2\nu^{0}_{1}&1+q&1\\
0&1-q^3\nu^{0}_{2}&1\\
\end{smallmatrix}\right|
$$
$$
M^{01}_{23}(1_3)=\left|\begin{smallmatrix}
1+q+q^2&1\\
1+q&1\\
\end{smallmatrix}\right|=q^2,\quad
M^0_3(2_3)=1,
$$
hence
$$
M^{i_0i_1i_2}_{123}(0_3)=
M^{i_0i_1i_2}_{123}(F_{0,3}^s(q,\lambda)),\,\,
M^{i_0i_1}_{23}(1_3)= M^{i_0i_1}_{23}(F_{1,3}^s(q,\lambda)),\,\,
$$
$$
M^{i_0}_{3}(2_3)=M^{i_0}_{3}(F_{2,3}^s(q,\lambda)).
$$
Since $D_4^\sharp(q)={\rm diag}(q^6,q^3,q,1,1)$ and
$$
(k_4)=\sigma_1(q,4)-q_{4-k}\lambda_k(D_4^\sharp(q)\Lambda_4)^{-1}$$
$$
=\left(\begin{smallmatrix}
1&(1+q)(1+q^2)&(1+q^2)(1+q+q^2)&(1+q)(1+q^2)&1\\
0&1&1+q+q^2&1+q+q^2&1\\
0&0&1&1+q&1\\
0&0&0&1&1\\
0&0&0&0&1\\
\end{smallmatrix}\right)-q_{4-k}\lambda_k\left(\begin{smallmatrix}
q^6\lambda_0&0&0&0&0\\
0&q^3\lambda_1&0&0&0\\
0&0&q\lambda_2&0&0\\
0&0&0&\lambda_3&0\\
0&0&0&0&\lambda_4\\
\end{smallmatrix}\right)^{-1}=
$$
$$
\left(\begin{smallmatrix}
1-\nu^{k}_{0}q_k/q_4&(1+q)(1+q^2)&(1+q^2)(1+q+q^2)&(1+q)(1+q^2)&1\\
0&1-\nu^{k}_{1}q_k/q_3&1+q+q^2&1+q+q^2&1\\
0&0&1-\nu^{k}_{2}q_k/q_2&1+q&1\\
0&0&0&1-\nu^{k}_{3}q_k/q_1&1\\
0&0&0&0&1-\nu^{k}_{4}q_k/q_0\\
\end{smallmatrix}\right),
$$
the equations (\ref{case3}) {\bf for} $n=4$ gives us
$$
\left(\begin{smallmatrix}
0&(1+q)(1+q^2)&(1+q^2)(1+q+q^2)&(1+q)(1+q^2)&1\\
0&1-q^3\nu^{0}_{1}&1+q+q^2&1+q+q^2&1\\
0&0&1-q^5\nu^{0}_{2}&1+q&1\\
0&0&0&1-q^6\nu^{0}_{3}&1\\
0&0&0&0&1-q^6\nu^{0}_{4}\\
\end{smallmatrix}\right)\left(\begin{smallmatrix} 0\\
b_{10}\\
b_{20}\\
b_{30}\\
b_{40}\\
\end{smallmatrix}\right)=0,
$$
$$
\left(\begin{smallmatrix}
1-q^{-3}\nu^{1}_{0}&(1+q)(1+q^2)&(1+q^2)(1+q+q^2)&(1+q)(1+q^2)&1\\
0&0&1+q+q^2&1+q+q^2&1\\
0&0&1-q^2\nu^{1}_{2}&1+q&1\\
0&0&0&1-q^3\nu^{1}_{3}&1\\
0&0&0&0&1-q^3\nu^{1}_{4}\\
\end{smallmatrix}\right)\left(\begin{smallmatrix} 0\\
0\\
b_{21}\\
b_{31}\\
b_{41}\\
\end{smallmatrix}\right)=0,
$$
$$
\left(\begin{smallmatrix}
1-q^{-5}\nu^{2}_{0}&(1+q)(1+q^2)&(1+q^2)(1+q+q^2)&(1+q)(1+q^2)&1\\
0&1-q^{-2}\nu^{2}_{1}&1+q+q^2&1+q+q^2&1\\
0&0&0&1+q&1\\
0&0&0&1-q\nu^{2}_{3}&1\\
0&0&0&0&1-q\nu^{2}_{4}\\
\end{smallmatrix}\right)\left(\begin{smallmatrix} 0\\
0\\
0\\
b_{32}\\
b_{42}\\
\end{smallmatrix}\right)=0,
$$
$$
\left(\begin{smallmatrix}
1-q^{-6}\nu^{3}_{0}&(1+q)(1+q^2)&(1+q^2)(1+q+q^2)&(1+q)(1+q^2)&1\\
0&1-q^{-3}\nu^{3}_{1}&1+q+q^2&1+q+q^2&1\\
0&0&1-q^{-1}\nu^{3}_{3}&1+q&1\\
0&0&0&0&1\\
0&0&0&0&1-\nu^{3}_{4}\\
\end{smallmatrix}\right)\left(\begin{smallmatrix} 0\\
0\\
0\\
0\\
b_{43}\\
\end{smallmatrix}\right)=0,
$$
hence for $n=4$ we have
$$
M^{i_0i_1i_2i_3}_{1234}(0_4)=
M^{i_0i_1i_2i_3}_{1234}(F_{1,4}^s(q,\lambda)),\,\,
M^{i_0i_1i_2}_{234}(1_4)=M^{i_0i_1i_2}_{234}(F_{1,4}^s(q,\lambda)),
$$
$$
M^{i_0i_1}_{34}(2_4)= M^{i_0i_1}_{23}(F_{2,4}^s(q,\lambda)),\,\,
M^{i_0}_{4}(3_4)=M^{i_0}_{4}(F_{3,4}^s(q,\lambda)).
$$
In general we conclude that the system of equations (\ref{case4})
$$
\sigma_1^{\Lambda,k}(q,n)b^{(k)}=0 ,\quad 0\leq k\leq n-1
$$
has only trivial solutions $b^{(k)}=0$ if and only if for any
$0\leq r\leq \left[\frac{n}{2}\right]$ there exists $0\leq
i_0<i_i<...<i_{n-r-1}\leq n $ such that (see (\ref{M^i_j...=0}))
$$
M^{i_0i_i...i_{n-r-1}}_{r+1r+2...n}(r_n)=
M^{i_0i_i...i_{n-r-1}}_{r+1r+2...n}
(F_{r,n}^s(q,\lambda))\not=0.
$$\qed\end{pf}
 {\bf Definition 3.}{\it We say that the values of
$\Lambda_n={\rm diag} (\lambda_k)_{k=0}^n$ are suspected (for
reducibility) if for some  $0\leq r\leq \left[\frac{n}{2}\right]$
(see (\ref{M^i_j...=0}))
$$
M^{i_0i_i...i_{n-r-1}}_{r+1r+2...n} (F_{rn}^s(q,\lambda))=0\,\,
\text{\,\,for \,\,all\,\,} 0\leq i_0<i_i<...<i_r\leq n.
$$}
Our {\bf aim now is to describe shortly the suspected values} of
$\Lambda_n$, (see definition (\ref{susp})). For example if $r=0$
we get that for all $0\leq i_0<i_i<...<i_{n-1}\leq n$
\begin{equation}
\label{shortM(q)=0}
M^{i_0i_i...i_{n-1}}_{12...n}(F_{0n}^s(q,\lambda))=0\Leftrightarrow
M^{01...n-1}_{12...n} (F_{0n}^s(q,\lambda))=0,\text{\,\, and\,\,}
M_n^n(F_{0n}^s(q,\lambda))=0.
\end{equation}
{\bf To  complete the proof of  the Theorem \ref{t.IrrB_3} we
should show that representation is reducible for suspected values
of $\Lambda_n$ }.
We should calculate the determinant:
$$D_3(q,\nu):=M^{012}_{123}\left[
\sigma_1(q,3)-q^3\lambda_0(D_3^\sharp(q)\Lambda_3)^{-1}\right]
$$
$$
=M^{012}_{123}\left[ \left(\begin{smallmatrix}
1&1+q+q^2&1+q+q^2&1\\
0&1&1+q&1\\
0&0&1&1\\
0&0&0&1\\
\end{smallmatrix}\right)-q^3\lambda_0\left(\begin{smallmatrix}
q^3\lambda_0&0&0&0\\
0&q\lambda_1&0&0\\
0&0&\lambda_2&0\\
0&0&0&\lambda_3\\
\end{smallmatrix}\right)^{-1}\right].
$$
{\bf Let us denote by (*) the conditions} (see (\ref{Rep(q)'}))
$\lambda_r\lambda_{n-r}=c,\,\,1\leq r\leq n$ and by $D_n(q,\nu)^*$
the value of $D_n(q,\nu)$ under these conditions.

We have
$$
D_3(q,\nu)=
\left|\begin{smallmatrix} 1+q+q^2&1+q+q^2&1\\
1-q^2\nu_1&1+q&1\\
0&1-q^3\nu_2&1
\end{smallmatrix}\right|
=\left|\begin{smallmatrix} q+q^2(1+\nu_1)&q^2\\
1-q^2\nu_1&q+q^3\nu_2\\
\end{smallmatrix}\right|
$$
$$
q^2\left|\begin{smallmatrix} 1+q(1+\nu_1)&1\\
1-q^2\nu_1&1+q^2\nu_2\\
\end{smallmatrix}\right|=q^3\left[
1+(1+q)\nu_1+ q(1+q)\nu_2+q^2\nu_1\nu_2\right].
$$
$$
\lambda_3=q\lambda_0\text{\,\,hence\,\,}(*)\Rightarrow\,\,\lambda_0\lambda_3=q\lambda_0^2=
\lambda_1\lambda_2,\,\, \text{\,\,so\,\,}\nu_1\nu_2=q^{-1},
$$
$$
 D_3^*(q,\nu)= q^3\left[ 1+q+(1+q)\nu_1+
q(1+q)\nu_2\right]
$$
$$
=q^3(1+q)\left[1+\nu_1+q\nu_2\right]=\frac{q^3(1+q)}{\lambda_1\lambda_2}
\left[\lambda_1\lambda_2+\lambda_0\lambda_2+q\lambda_0\lambda_1\right]
$$
$$
=\frac{q^3(1+q)}{\lambda_1\lambda_2}
\left[q\lambda_0^2+\lambda_0\lambda_2+q\lambda_0\lambda_1\right]=
\frac{q^3(1+q)\lambda_0}{\lambda_1\lambda_2}
\left[q\lambda_0+q\lambda_1+\lambda_2\right]=0.
$$
Finally we get
$$
 D_3^*(q,\nu)=\frac{q^4(1+q)\lambda_0^2}{\lambda_1\lambda_2}
\left[1+\alpha_1+q^{-1}\alpha_2\right]=0\text{\,\,with\,\,}
\alpha_1\alpha_2=q,
\text{\,\,where\,\,}\alpha_k=\lambda_0/\lambda_k.
$$
If we replace
$\tilde\alpha_1=\alpha_1,\,\,\tilde\alpha_2=q^{-1}\alpha_2$ we get
$\tilde\alpha_1+\tilde\alpha_2=-1,\,\,\tilde\alpha_1\tilde\alpha_2=1.$
Using (\ref{alpha3}) and (\ref{red3}) we conclude that
$$
\tilde\alpha_{k}=(-1\pm i\sqrt{3})/2=\exp( \pm\frac{2\pi i
k}{3}),\,\,k=1,2,
$$
$$
\Lambda_3 =\lambda_0{\rm diag}(1,(-1\pm i\sqrt{3})/2,q(-1\mp
i\sqrt{3})/2,q^3)=
$$
$$
\lambda_0{\rm
diag}(\alpha_0,\alpha_1,q\alpha_2,q^3\alpha_3)=\lambda_0D_3(q){\rm
diag}(\alpha_k)_{k=0}^3 ,\,\,\text{\,\,where\,\,}\alpha_k=\exp(
\pm\frac{2\pi i k}{3}).
$$
The latter values of the matrix $\Lambda_3$ contradicts with
conditions $\lambda_r\lambda_{n-r}=c$. Indeed,
$\lambda_1\lambda_2/\lambda_0^2=\alpha_1q\alpha_2=q$ but
$\lambda_0\lambda_3/\lambda_0^2=\alpha_0q^3\alpha_3=q^3$. They
coincide when $q=q^3$ or $q=\pm 1.$

{\bf For} $n=4$ we have the following determinant:
$$D_4(q,\nu):=M^{0123}_{1234}\left[
\sigma_1(q,4)-q^6\lambda_0(D_4^\sharp(q)\Lambda_4)^{-1}\right]
$$
$$
=M^{012}_{123}\left[ \left(\begin{smallmatrix}
1&(1+q)(1+q^2)&(1+q^2)(1+q+q^2)&(1+q)(1+q^2)&1\\
0&1&1+q+q^2&1+q+q^2&1\\
0&0&1&1+q&1\\
0&0&0&1&1\\
0&0&0&0&1\\
\end{smallmatrix}\right)-q^6\lambda_0\left(\begin{smallmatrix}
q^6\lambda_0&0&0&0&0\\
0&q^3\lambda_1&0&0&0\\
0&0&q\lambda_2&0&0\\
0&0&0&\lambda_3&0\\
0&0&0&0&\lambda_4\\
\end{smallmatrix}\right)^{-1}\right].
$$
We have
$$
D_4(q,\nu)= \left|\begin{smallmatrix}
(1+q)(1+q^2)&(1+q^2)(1+q+q^2)&(1+q)(1+q^2)&1\\
1-q^3\nu_1&1+q+q^2&1+q+q^2&1\\
0&1-q^5\nu_2&1+q&1\\
0&0&1-q^6\nu_3&1
\end{smallmatrix}\right|=
\left|\begin{smallmatrix} q+q^2+q^3(1+\nu_1)&q^2(1+q+q^2)&q^3\\
1-q^3\nu_1&q+q^2+q^5\nu_2&q^2\\
0&1-q^5\nu_2&q+q^6\nu_3\\
\end{smallmatrix}\right|
$$
In the general case $n\in{\mathbb N}$ we should find
$$D_n(q,\nu) \text{\,\,and\,\,} D_n(q,\nu)^*$$ and
prove the reducibility.
\subsection{Subspace irreducibility}

Let us denote by $\sigma^\Lambda(q,n)$ the representation of $B_3$
defined by (\ref{Rep(q)'}).

{\bf Problem.} {\it To find a criteria of the subspace
irreducibility for all representations $\sigma^\Lambda(q,n)$}.

We study here only some particular cases.

{\bf Theorem 4.}{\it The representation $\sigma^\Lambda(q,n)$ is
{\rm subspace irreducible} for $n=1$ if and only if
$\Lambda_1\not=\lambda_0(1,\alpha)$ where $\alpha^2-\alpha+1=0$.}
\begin{pf}
{\it Reducibility} for $n=1$ . The eigenvalues of
$\sigma_1^\Lambda(1,1)=\left(\begin{smallmatrix}
\lambda_0&\lambda_1\\
0&\lambda_1
\end{smallmatrix}\right)$ and $\sigma_2^\Lambda(1,1)=
\left(\begin{smallmatrix}
\lambda_1&0\\
-\lambda_0&\lambda_0
\end{smallmatrix}\right)$ are the following
$$
e_{\lambda_0}^{(1)}=\left(\begin{smallmatrix} 1\\
0\\
\end{smallmatrix}\right),\,\,e_{\lambda_1}^{(1)}=\left(\begin{smallmatrix} \lambda_1\\
\lambda_1-\lambda_0\\
\end{smallmatrix}\right),\,\,
e_{\lambda_0}^{(2)}=\left(\begin{smallmatrix} 0\\
1\\
\end{smallmatrix}\right),\,\,e_{\lambda_1}^{(1)}=\left(\begin{smallmatrix}
\lambda_1^{-1}-\lambda_0^{-1}\\
\lambda_1^{-1}\\
\end{smallmatrix}\right).
$$
or in a short form
$$ (e_{\lambda_0}^{(1)},e_{\lambda_1}^{(1)}):=
\left(\begin{smallmatrix} 1&\lambda_1\\
0&\lambda_1-\lambda_0\\
\end{smallmatrix}\right),\quad
(e_{\lambda_0}^{(2)},e_{\lambda_1}^{(2)}):=
\left(\begin{smallmatrix}
0&\lambda_1^{-1}-\lambda_0^{-1}\\
1&\lambda_1^{-1}\\
\end{smallmatrix}\right).
$$
We see that $e_{\lambda_1}^{(1)}$ and $e_{\lambda_1}^{(2)}$ are
linearly independent if and only if
\begin{equation}
\label{irr1} {\rm det
}_1(\lambda_0,\lambda_1):=\frac{1}{\lambda_0\lambda_1}\left|
\begin{smallmatrix}
\lambda_1&\lambda_0-\lambda_1\\
\lambda_1-\lambda_0&\lambda_0\\
\end{smallmatrix}
\right|=
 \left|
\begin{smallmatrix}
\lambda_1&\lambda_1^{-1}-\lambda_0^{-1}\\
\lambda_1-\lambda_0&\lambda_1^{-1}\\
\end{smallmatrix}
\right|=\frac{\lambda_0}{\lambda_1}+\frac{\lambda_1}{\lambda_0}-1\not=0.
\end{equation}
If we set $\alpha=\frac{\lambda_1}{\lambda_0}$ we get
$\alpha+\alpha^{-1}-1=0$ or $\alpha^2-\alpha+1=0$ hence if
\begin{equation}
\label{Irr,n=2} \Lambda_2=\lambda_0{\rm
diag}(1,\alpha):\,\,\alpha^2-\alpha+1=0,\,\,
\alpha_{1,2}=\frac{1}{2}\pm\frac{i\sqrt{3}}{2},\quad
\alpha_{1,2}=\exp(\pm 2\pi i/6),
\end{equation}
the representation is reducible.

{\it Irreducibility}. We have
$$
A_1:=\sigma_1^\Lambda(1,1)-\lambda_0I=\left(\begin{smallmatrix}
0&\lambda_1\\
0&\lambda_1-\lambda_0\\
\end{smallmatrix}\right),\quad
A_2:=\sigma_2^\Lambda(1,1)-\lambda_0I= \left(\begin{smallmatrix}
\lambda_1-\lambda_0&0\\
-\lambda_0&0\\
\end{smallmatrix}\right)
$$
If $\lambda_1-\lambda_0=0$ we get
$E_{01}=\left(\begin{smallmatrix}
0&1\\
0&0\\
\end{smallmatrix}\right)$ and $E_{10}=\left(\begin{smallmatrix}
0&0\\
1&0\\
\end{smallmatrix}\right)$. They generate the matrix algebra ${\rm Mat}(2,{\mathbb
C})$ (see Remark) hence the representation is irreducible.

If $\lambda_1-\lambda_0\not=0$ we have
$$
A_1A_2=\left(\begin{smallmatrix}
0&\lambda_1\\
0&\lambda_1-\lambda_0\\
\end{smallmatrix}\right)\left(\begin{smallmatrix}
\lambda_1-\lambda_0&0\\
-\lambda_0&0\\
\end{smallmatrix}\right)=-\lambda_0\left(\begin{smallmatrix}
\lambda_1&0\\
\lambda_1-\lambda_0&0\\
\end{smallmatrix}\right),\quad A_2=\left(\begin{smallmatrix}
\lambda_1-\lambda_0&0\\
-\lambda_0&0\\
\end{smallmatrix}\right),
$$
$$
A_2A_1=\left(\begin{smallmatrix}
\lambda_1-\lambda_0&0\\
-\lambda_0&0\\
\end{smallmatrix}\right)\left(\begin{smallmatrix}
0&\lambda_1\\
0&\lambda_1-\lambda_0\\
\end{smallmatrix}\right)=\lambda_1\left(\begin{smallmatrix}
0&\lambda_1-\lambda_0\\
0&-\lambda_0\\
\end{smallmatrix}\right),\quad
A_1=\left(\begin{smallmatrix}
0&\lambda_1\\
0&\lambda_1-\lambda_0\\
\end{smallmatrix}\right),
$$
hence $A_1A_2$ and $A_2$ (resp. $A_2A_1$ and $A_1$) generate
elements $E_{00}$ and  $E_{10}$ (resp. $E_{01}$ and  $E_{11}$) iff
${\rm det }_1(\lambda_0,\lambda_1)\not=0$, the representation is
irreducible. \qed\end{pf}
{\bf For} $n=2$ {\it reducibility}. The eigenvalues of
$$
\sigma_1^\Lambda(1,2)=\left(\begin{smallmatrix}
\lambda_0&2\lambda_1&\lambda_2\\
0&\lambda_1&\lambda_2\\
0&0&\lambda_2
\end{smallmatrix}\right),\,\,\sigma_2^\Lambda(1,2)=\left(\begin{smallmatrix}
\lambda_2&0&0\\
\lambda_1&-\lambda_1&0\\
\lambda_0&-2\lambda_0&\lambda_0\\
\end{smallmatrix}\right)
$$
are as follows
$$
e_{\lambda_0}^{(1)}=\left(\begin{smallmatrix} 1\\
0\\
0\\
\end{smallmatrix}\right),\,\,e_{\lambda_1}^{(1)}=
\left(\begin{smallmatrix} 2\lambda_1\\
\lambda_1-\lambda_0\\
0\\
\end{smallmatrix}\right),\,\,
e_{\lambda_2}^{(1)}=\left(\begin{smallmatrix} \lambda_2(\lambda_2+\lambda_1)\\
\lambda_2(\lambda_2-\lambda_0)\\
(\lambda_2-\lambda_1)(\lambda_2-\lambda_0)\\
\end{smallmatrix}\right),
$$
$$
e_{\lambda_0}^{(2)}=\left(\begin{smallmatrix} 0\\
0\\
1\\
\end{smallmatrix}\right),\,\,e_{\lambda_1}^{(2)}=\left(\begin{smallmatrix}
0\\
\lambda_1^{-1}-\lambda_0^{-1}\\
\lambda_1^{-1}\\
\end{smallmatrix}\right),\,\,
e_{\lambda_2}^{(2)}=\left(\begin{smallmatrix}
(\lambda_2^{-1}-\lambda_1^{-1})(\lambda_2^{-1}-\lambda_0^{-1})\\
\lambda_2^{-1}(\lambda_2^{-1}-\lambda_0^{-1})\\
\lambda_2^{-1}(\lambda_2^{-1}+\lambda_1^{-1})\\
\end{smallmatrix}\right)
$$
or in the short form
$$
(e_{\lambda_0}^{(1)},e_{\lambda_1}^{(1)},e_{\lambda_2}^{(1)})
=\left(\begin{smallmatrix} 1&2\lambda_1&\lambda_2(\lambda_2+\lambda_1)\\
0&\lambda_1-\lambda_0&\lambda_2(\lambda_2-\lambda_0)\\
0&0&(\lambda_2-\lambda_0)(\lambda_2-\lambda_1)\\
\end{smallmatrix}\right),
$$
$$
(e_{\lambda_0}^{(2)},e_{\lambda_1}^{(2)},e_{\lambda_2}^{(2) })=
\left(\begin{smallmatrix} 0&0&
(\lambda_2^{-1}-\lambda_1^{-1})(\lambda_2^{-1}-\lambda_0^{-1})\\
0&\lambda_1^{-1}-\lambda_0^{-1}&\lambda_2^{-1}(\lambda_2^{-1}-\lambda_0^{-1})\\
1&\lambda_1^{-1}&\lambda_2^{-1}(\lambda_2^{-1}+\lambda_1^{-1})\\
\end{smallmatrix}\right),
$$
$$
(e_{\lambda_0}^{(2)},e_{\lambda_1}^{(2)},e_{\lambda_2}^{(2)})_{ij}=
(e_{\lambda_0}^{(1)},e_{\lambda_1}^{(1)},e_{\lambda_2}^{(1)})_{2-i,j}.
$$


Let now $\lambda_2-\lambda_0=0$. To find $e_{\lambda_0}^{(1)}$ we
write $(\sigma_1^\Lambda(1,2)-\lambda_0I)a=0$ or
$$
\left(\begin{smallmatrix}
0&2\lambda_1&\lambda_0\\
0&\lambda_1-\lambda_0&\lambda_0\\
0&0&0\\
\end{smallmatrix}\right)\left(\begin{smallmatrix}
a_0\\
a_1\\
a_2\\
\end{smallmatrix}\right)=0,\,\,\left\{\begin{smallmatrix}
2\lambda_1a_1+\lambda_0a_2=0\\
(\lambda_1-\lambda_0)a_1+\lambda_0a_2=0\\
\end{smallmatrix}\right.
$$
We see that $(a_1,a_2)=(0,0)$ i.e. $e_{\lambda_0}^{(1)}=(1,0,0)$
iff $\lambda_1+\lambda_0\not=0$. Indeed
$$\left|
\begin{smallmatrix}
2\lambda_1&\lambda_0\\
\lambda_1-\lambda_0&\lambda_0\\
\end{smallmatrix}
\right|=\lambda_0(\lambda_1+\lambda_0)=0\Leftrightarrow
\lambda_1+\lambda_0=0.
$$
In this case we have
$$
(e_{\lambda_0}^{(1)},e_{\lambda_1}^{(1)})
=\left(\begin{smallmatrix} 1&2\lambda_1\\
0&\lambda_1-\lambda_0\\
0&0\\
\end{smallmatrix}\right),\quad
(e_{\lambda_0}^{(2)},e_{\lambda_1}^{(2)})=
\left(\begin{smallmatrix} 0&0\\
0&\lambda_1^{-1}-\lambda_0^{-1}\\
1&\lambda_1^{-1}\\
\end{smallmatrix}\right),
$$
so representation is irreducible. If $\lambda_1+\lambda_0=0$ we
conclude that $(a_1,a_2)=(\lambda_0,-2\lambda_1)$ and
$e_{\lambda_0}^{(1),t}=(t,\lambda_0,-2\lambda_1).$ To find
$e_{\lambda_0}^{(2)}$ we write
$(\sigma_2^\Lambda(1,2)-\lambda_0I)a=0$ or
$$
\left(\begin{smallmatrix}
0&0&0\\
\lambda_1&\lambda_1-\lambda_0&0\\
\lambda_0&-2\lambda_0&0\\
\end{smallmatrix}\right)\left(\begin{smallmatrix}
a_0\\
a_1\\
a_2\\
\end{smallmatrix}\right)=0,\,\,\left\{\begin{smallmatrix}
\lambda_1a_0+(\lambda_1-\lambda_0)a_1=0\\
\lambda_0a_0-2\lambda_0a_1=0\\
\end{smallmatrix}\right.
$$
$e_{\lambda_0}^{(2),\tau}=(2\lambda_0,\lambda_0,\tau).$ We have
$$
(e_{\lambda_0}^{(1),t},e_{\lambda_1}^{(1)})
=\left(\begin{smallmatrix} t&2\lambda_1\\
\lambda_0&\lambda_1-\lambda_0\\
-2\lambda_1&0&\\
\end{smallmatrix}\right),\quad
(e_{\lambda_0}^{(2),\tau},e_{\lambda_1}^{(2)})=
\left(\begin{smallmatrix} 2\lambda_0&0\\
\lambda_0&\lambda_1^{-1}-\lambda_0^{-1}\\
\tau&\lambda_1^{-1}\\
\end{smallmatrix}\right).
$$
Two vectors $(a,b,t)$ and $(\tau,c,d)$ are proportional if and
only if
$$
\left|
\begin{smallmatrix}
a&b\\
\tau&c\\
\end{smallmatrix}
\right|=0,\,\,\left|
\begin{smallmatrix}
b&t\\
c&d\\
\end{smallmatrix}
\right|=0,\,\, \left|\begin{smallmatrix}
a&t\\
\tau&d\\
\end{smallmatrix}
\right|=0.
$$
If $b\not=0$ and $c\not=0$ we have
$\tau=\frac{ac}{b},\,\,t=\frac{bd}{c}$ and hence
$\left|\begin{smallmatrix}
a&t\\
\tau&d\\
\end{smallmatrix}
\right|=0$. Two vectors
$e_{\lambda_0}^{(1),t}=(t,\lambda_0,-2\lambda_1)$ and
$e_{\lambda_0}^{(2),\tau}=(2\lambda_0,\lambda_0,\tau)$ are
proportional for $t=2\lambda_0$ and $\tau=-2\lambda_1$. In general
(without condition $\lambda_0\lambda_2=\lambda_1^2$) the family of
two matrix $\sigma_1^\Lambda(1,2)$ and $\sigma_2^\Lambda(1,2)$ is
irreducible if and only if $\lambda_2-\lambda_0=0$ and
$\lambda_0+\lambda_1\not=0$.

In our case we have $\lambda_0\lambda_2=\lambda_1^2$ hence
$\lambda_0^2=\lambda_1^2$  and $\lambda_0=\pm\lambda_1$. If
$\lambda_0=\lambda_1$  we get $\Lambda_2=\lambda_0{\rm
diag}(1,1,1)$, the representation is irreducible, if
$\lambda_0=-\lambda_1$ we get $\Lambda_2=\lambda_0{\rm
diag}(1,-1,1)$, the representation is reducible.

{\it Irreducibility}. Let us denote
$A_i=(\sigma_i^\Lambda(1,2)-\lambda_0I)(\sigma_i^\Lambda(1,2)-\lambda_1I),\,\,i=1,2$.
We have
$$
A_1
=\left(\begin{smallmatrix} 0&0&\lambda_2(\lambda_2+\lambda_1)\\
0&0&\lambda_2(\lambda_2-\lambda_0)\\
0&0&(\lambda_2-\lambda_1)(\lambda_2-\lambda_0)\\
\end{smallmatrix}\right),\quad
A_2=
\left(\begin{smallmatrix} (\lambda_2-\lambda_1)(\lambda_2-\lambda_0)&0&0\\
-\lambda_1(\lambda_2-\lambda_0)&0&0&\\
\lambda_0(\lambda_2+\lambda_1)&0&0\\
\end{smallmatrix}\right)
$$
$$
A_1A_2=\lambda_2(\lambda_2+\lambda_1)
\left(\begin{smallmatrix} \lambda_2(\lambda_2+\lambda_1)&0&0\\
\lambda_2(\lambda_2-\lambda_0)&0&0\\
(\lambda_2-\lambda_1)(\lambda_2-\lambda_0)&0&0\\
\end{smallmatrix}\right),\,\,
A_2A_1= \lambda_2(\lambda_2+\lambda_1)
\left(\begin{smallmatrix}0&0& (\lambda_2-\lambda_1)(\lambda_2-\lambda_0)\\
0&0&-\lambda_1(\lambda_2-\lambda_0)\\
0&0&\lambda_0(\lambda_2+\lambda_1)\\
\end{smallmatrix}\right).
$$
\begin{rem}
For $n=2$ the representation is subspace irreducible if and only
if $\Lambda_2=\lambda_0{\rm diag}(1,1,1).$
\end{rem}
 {\bf For} $n=3$ the eigenvalues of $\sigma_1^\Lambda(1,3)$ and
$\sigma_2^\Lambda(1,3)$ are
$$
(e_{\lambda_0}^{(1)},e_{\lambda_1}^{(1)},e_{\lambda_2}^{(1)},e_{\lambda_3}^{(1)})=
\left(\begin{smallmatrix}
1&3\lambda_1&3\lambda_2(\lambda_2+\lambda_1)&
\lambda_3[(\lambda_3+\lambda_2)(\lambda_3+\lambda_1)+\lambda_3(\lambda_2+\lambda_1)]\\
0&\lambda_1-\lambda_0&2\lambda_2(\lambda_2-\lambda_0)&
\lambda_3(\lambda_3+\lambda_2)(\lambda_3-\lambda_0)\\
0&0&(\lambda_2-\lambda_0)(\lambda_2-\lambda_1)&
\lambda_3(\lambda_3-\lambda_1)(\lambda_3-\lambda_0)\\
0&0&0&(\lambda_3-\lambda_2)(\lambda_3-\lambda_1)(\lambda_3-\lambda_0)\\
\end{smallmatrix}\right),
$$
$$
(e_{\lambda_0}^{(2)},e_{\lambda_1}^{(2)},e_{\lambda_2}^{(2)},e_{\lambda_3}^{(2)})_{ij}=
(e_{\lambda_0^{-1}}^{(1)},e_{\lambda_1^{-1}}^{(1)},e_{\lambda_2^{-1}}^{(1)},
e_{\lambda_3^{-1}}^{(1)})_{3-i,j}.
$$
{\bf For} $n=4$ the eigenvalues of $\sigma_1^\Lambda(1,4)$ and
$\sigma_2^\Lambda(1,4)$ are
$(e_{\lambda_0}^{(1)},e_{\lambda_1}^{(1)},e_{\lambda_2}^{(1)},e_{\lambda_2}^{(1)})=$
$$
 \left(\begin{smallmatrix}
1&4\lambda_1&6\lambda_2(\lambda_2+\lambda_1)&
4\lambda_3[(\lambda_3+\lambda_2)(\lambda_3+\lambda_1)+\lambda_3(\lambda_2+\lambda_1)]&
\lambda_4[(\lambda_4+\lambda_3)(\lambda_4+\lambda_2)(\lambda_4+\lambda_1)+
]
\\
0&\lambda_1-\lambda_0&3\lambda_2(\lambda_2-\lambda_0)&
3\lambda_3(\lambda_3+\lambda_2)(\lambda_3-\lambda_0)&
\lambda_4(\lambda_4+\lambda_3)(\lambda_4+\lambda_2)(\lambda_4-\lambda_0)
\\
0&0&(\lambda_2-\lambda_0)(\lambda_2-\lambda_1)&
\lambda_3(\lambda_3-\lambda_1)(\lambda_3-\lambda_0)&
\lambda_4(\lambda_4+\lambda_3)(\lambda_4-\lambda_1)(\lambda_4-\lambda_0)
\\
0&0&0&(\lambda_3-\lambda_2)(\lambda_3-\lambda_1)(\lambda_3-\lambda_0)&
\lambda_4(\lambda_4-\lambda_2)(\lambda_4-\lambda_1)(\lambda_4-\lambda_0)
\\
0&0&0&0&(\lambda_4-\lambda_3)(\lambda_4-\lambda_2)(\lambda_4-\lambda_1)(\lambda_4-\lambda_0)\\
\end{smallmatrix}\right),
$$
$$
(e_{\lambda_0}^{(2)},e_{\lambda_1}^{(2)},e_{\lambda_2}^{(2)},
e_{\lambda_3}^{(2)}, e_{\lambda_4}^{(2)})_{ij}=
(e_{\lambda_0^{-1}}^{(1)},e_{\lambda_1^{-1}}^{(1)},e_{\lambda_2^{-1}}^{(1)},
e_{\lambda_3^{-1}}^{(1)}, e_{\lambda_4^{-1}}^{(1)})_{4-i,j}.
$$
In general for different $\lambda_k,\,\,0\leq k\leq n$ we get
$$
(e_{\lambda_0}^{(2)},e_{\lambda_1}^{(2)},...,e_{\lambda_n}^{(2)})_{ij}=
(e_{\lambda_0^{-1}}^{(1)},e_{\lambda_1^{-1}}^{(1)},...,
e_{\lambda_n^{-1}}^{(1)})_{n-i,j},\,\,0\leq i,j\leq n.
$$
\begin{rem}
To study the subspace irreducibility it is necessary to compare
all possible subspaces generated by the eigenvalues of
$\sigma_1^\Lambda (1,n)$ and by the eigenvalues of
$\sigma_2^\Lambda (1,n)$.
\end{rem}
\subsection{Subspace reducibility}
Representation is irreducible in the case 1). We know only {\bf
some particular cases} in the case 2). We show that
representations are reducible for suspected values of $\Lambda_n$
for small $n$ and $r=0$ (see definition 1, p.24). Let $n=2$ and
$\Lambda_2^{(2)}={\rm diag }(1,-1,1)$. We have
$$
\sigma_1^\Lambda(1,2)=\left(\begin{smallmatrix}
1&-2&1\\
0&-1&1\\
0&0&1\\
\end{smallmatrix}\right),\quad
\sigma_1^\Lambda(1,2)=\left(\begin{smallmatrix}
1&0&0\\
1&-1&0\\
1&-2&1\\
\end{smallmatrix}\right).
$$
 The eigenvectors $e_0$ for $\sigma_1^\Lambda(1,2)$ and
$f_0$ for  $\sigma_2^\Lambda(1,2)$ corresponding to the eigenvalue
$\lambda_0=1$ are the following
$$
e_0^t=(t,1,2),\,\,f_0^\tau=(2,1,\tau).
$$
{\small Indeed we have
$$
\left(\begin{smallmatrix}
0&-2&1\\
0&-1&1\\
0&0&0\\
\end{smallmatrix}\right)
\left(\begin{smallmatrix}
t\\
1\\
2\\
\end{smallmatrix}\right)=0,
\quad \left(\begin{smallmatrix}
0&0&0\\
1&-2&0\\
1&-2&0\\
\end{smallmatrix}\right)
\left(\begin{smallmatrix}
2\\
1\\
\tau\\
\end{smallmatrix}\right)=0.
$$ }
We conclude that
$$
\langle e_0^t=(t,1,2)\rangle =\langle
f_0^\tau=(2,1,\tau)\rangle\Leftrightarrow t=\tau=2.
$$
Let $n=3$. We have two suspected cases:
$\Lambda_3^{(s)}=\lambda_0{\rm diag }(\alpha^{(s)}_k)_{k=0}^3,\,\,
s=2,3$ where $\alpha^{(s)}_k=\exp(2\pi i k/s )$. We get for
$\Lambda_2^{(2)}={\rm diag }(1,-1,1-1)$
$$
\sigma_1^\Lambda(1,2)=\left(\begin{smallmatrix} 1&-3&3&-1\\
0&-1&2&-1\\
0&0&1&-1\\
0&0&0&-1\\
\end{smallmatrix}\right),\quad
\sigma_2^\Lambda(1,2)=\left(\begin{smallmatrix}
-1&0&0 &0\\
-1&1&0 &0\\
-1&2&-1&0\\
-1&3&-3&1\\
\end{smallmatrix}\right)
$$
The eigenvectors $e_0$ for $\sigma_1^\Lambda(1,3)$ and $f_0$ for
$\sigma_2^\Lambda(1,3)$ corresponding to the eigenvalue
$\lambda_0=1$ are the following
$$
e_0^t=(t,1,1,0),\quad f_0^\tau=(0,1,1,\tau).
$$
{\small Indeed we have
$$
\left(\begin{smallmatrix} 0&-3&3&-1\\
0&-2&2&-1\\
0&0&0&-1\\
0&0&0&-2\\
\end{smallmatrix}\right)
\left(\begin{smallmatrix}
t\\
1\\
1\\
0\\
\end{smallmatrix}\right)=0,\quad
\left(\begin{smallmatrix}
-2&0& 0&0\\
-1&0& 0&0\\
-1&2&-2&0\\
-1&3&-3&0\\
\end{smallmatrix}\right)
\left(\begin{smallmatrix}
0\\
1\\
1\\
\tau\\
\end{smallmatrix}\right)=0.
$$ }
We conclude that
$$
\langle e_0^t=(t,1,1,0)\rangle =\langle
f_0^\tau=(0,1,1,\tau)\rangle\Leftrightarrow t=\tau=0.
$$
\begin{lem}
\label{l.s-red.n=3.4} In the suspected cases for $r=0$ and $n=3,4$
the representations $\sigma^\Lambda(1,n)$ is subspace irreducible.
\end{lem}
\begin{pf} For $n=3$ and $r=0$ the suspected values of $\Lambda_3$
are as follows
$$
\Lambda_3=\Lambda_3^{(2)}:=\lambda_0{\rm diag
}(1,-1,1,-1),\,\,\text{\,\,and\,\,\,\,}
\Lambda_3=\Lambda_3^{(3)}:=\lambda_0{\rm diag }(\exp 2\pi
ik/3)_{k=0}^3.
$$
If we set $\alpha_k=\exp(2\pi ik/3)$ and take $\lambda_0=1$ we get
$\Lambda_3^{(3)}:={\rm diag }(1,\alpha_1,\alpha_2,1)$ and
$$
\sigma_1^\Lambda(1,3)=\sigma_1(1,3)\Lambda_3=\left(\begin{smallmatrix}
1&3\alpha_1&3\alpha_2&1\\
0&\alpha_1&2\alpha_2&1\\
0&0&\alpha_2&1\\
0&0&0&1\\
\end{smallmatrix}\right),\quad
\sigma_2^\Lambda(1,3)=\Lambda_3^\sharp\sigma_2(1,3)=\left(\begin{smallmatrix}
1&0&0&0\\
-\alpha_2&\alpha_2&0&0\\
\alpha_1&-2\alpha_1&\alpha_1&0\\
-1&3&-3&1\\
\end{smallmatrix}\right).
$$
To find the eigenvectors $e_k$ for $\sigma_1^\Lambda(1,3)$ and
$f_k$ for  $\sigma_2^\Lambda(1,3)$ we have
$$ (\sigma_1^\Lambda(1,3)-I)e_0=0,\quad
(\sigma_1^\Lambda(1,3)-\alpha_1I)e_1=0,\quad
(\sigma_1^\Lambda(1,3)-\alpha_2I)e_2=0,
$$
$$ (\sigma_2^\Lambda(1,3)-I)f_0=0,\quad
(\sigma_2^\Lambda(1,3)-\alpha_1I)f_1=0,\quad
(\sigma_2^\Lambda(1,3)-\alpha_2I)f_2=0
$$
or
$$
\left(\begin{smallmatrix}
0&3\alpha_1&3\alpha_2&1\\
0&\alpha_1-1&2\alpha_2&1\\
0&0&\alpha_2-1&1\\
0&0&0&0\\
\end{smallmatrix}\right)\left(\begin{smallmatrix}
a_0\\
a_1\\
a_2\\
a_3\\
\end{smallmatrix}\right)=0,\,\,
\left(\begin{smallmatrix}
1-\alpha_1&3\alpha_1&3\alpha_2&1\\
0&0&2\alpha_2&1\\
0&0&\alpha_2-\alpha_1&1\\
0&0&0&1-\alpha_1\\
\end{smallmatrix}\right)\left(\begin{smallmatrix}
a_0\\
a_1\\
a_2\\
a_3\\
\end{smallmatrix}\right)=0,
$$
$$
\left(\begin{smallmatrix}
1-\alpha_2&3\alpha_1&3\alpha_2&1\\
0&\alpha_1-\alpha_2&2\alpha_2&1\\
0&0&0&1\\
0&0&0&1-\alpha_2\\
\end{smallmatrix}\right)\left(\begin{smallmatrix}
a_0\\
a_1\\
a_2\\
a_3\\
\end{smallmatrix}\right)=0.
$$
The solutions are the following
$$
\begin{smallmatrix}
e_0^t=(&t,      &1+\alpha_2,   &1-\alpha_1,&(1-\alpha_1)(1-\alpha_2)&),\\
e_1=(&3\alpha_1,&-(1-\alpha_1),&0,          &0&),\\
e_2=(&3(\alpha_1-1)^{-1},&2\alpha_2,&\alpha_2-\alpha_1,          &0&),\\
f_0^s=(&3(1-\alpha_1),&3,&2+\alpha_1,&s&),\\
f_1=(&0,&0,&1-\alpha_1,&3&),\\
f_2=(&0,&1-\alpha_1,&2,&1-\alpha_2&).\\
\end{smallmatrix}\,
$$
or
$$
\begin{smallmatrix}
e_0^t=(&t,      &\exp( 2\pi i 10/12),   &\sqrt{3}\exp( 2\pi i11 /12),&3&),\\
e_1=(&3\exp( 2\pi i4 /12),&\sqrt{3}\exp( 2\pi i5 /12),&0,          &0&),\\
e_2=(&\sqrt{3}\exp( 2\pi i7 /12),&2\exp( 2\pi i8 /12),&\sqrt{3}\exp( 2\pi i9 /12), &0&),\\
f_0^s=(&3\sqrt{3}\exp( 2\pi i11 /12),&3,&\sqrt{3}\exp( 2\pi i /12),&s&),\\
f_1=(&0,&0,&\sqrt{3}\exp( 2\pi i 11/12),&3&),\\
f_2=(&0,&\sqrt{3}\exp( 2\pi i11 /12),&2,&\sqrt{3}\exp( 2\pi i /12)&).\\
\end{smallmatrix}
$$
If we set $e(k):=\exp (2\pi i k/12)$
$$
\begin{smallmatrix}
e_0^t=(&t,      &e(10),   &\sqrt{3}e(11),&3&),\\
e_1=(&3e(4),&\sqrt{3}e(5),&0,          &0&),\\
e_2=(&\sqrt{3}e(7),&2e(8),&\sqrt{3}e(9), &0&),\\
f_0^\tau=(&3\sqrt{3}e(11),&3,&\sqrt{3}e(k),&\tau&),\\
f_1=(&0,&0,&\sqrt{3}e(11),&3&),\\
f_2=(&0,&\sqrt{3}e(11),&2,&\sqrt{3}e(1)&),\\
\end{smallmatrix}\text{\quad or\quad}
\begin{smallmatrix}
e_0^t=(&t,      &1,   &\sqrt{3}e(1),&3e(2)&),\\
e_1=(&\sqrt{3},&e(1),&0,          &0&),\\
e_2=(&\sqrt{3},&2e(1),&\sqrt{3}e(2), &0&),\\
f_0^\tau=(&3,&\sqrt{3}e(1),&e(2),&\tau&),\\
f_1=(&0,&0,&1,&\sqrt{3}e(1)&),\\
f_2=(&0,&\sqrt{3},&2e(1),&\sqrt{3}e(2)&).\\
\end{smallmatrix}
$$
We see that there no  one-dimensional invariant subspaces for
$\sigma_1^\Lambda(1,3)$ and $\sigma_2^\Lambda(1,3)$. We find the
two-dimensional subspace in the following form
$$
V_2^t:=\langle e_0^t,e_1\rangle.
$$
We can verify that $f_1\in V_2^t$ if and only if
$t=t_0:=\sqrt{3}e(-1)$. Further we get that $f_0^s\in V_2^{t_0}$
if and only if $s=s_0=\sqrt{3}e(3).$ Finally we conclude that two
dimensional subspace $V_2^{t_0}$ is invariant for
$\sigma_1^\Lambda(1,3)$ and $\sigma_2^\Lambda(1,3)$ since it is
generated by eigenvectors $e_0^{t_0},e_1$ for
$\sigma_1^\Lambda(1,3)$ and by eigenvectors $f_0^{s_0},f_1$ for
$\sigma_2^\Lambda(1,3)$.\\
Let $n=4$. We have three suspected cases:
$\Lambda_3^{(s)}=\lambda_0{\rm diag }(\alpha^{(s)}_k)_{k=0}^4,\,\,
s=2,3,4$. We get for $\Lambda_4^{(2)}={\rm diag }(1,-1,1-1,1)$
(resp. for  $\Lambda_4^{(3)}$ and $\Lambda_4^{(4)}$)
$$
\sigma_1^\Lambda(1,3)=\left(\begin{smallmatrix} 1&-4&6&-4&1\\
0&-1&3&-3&1\\
0&0&1&-2 &1\\
0&0&0&-1 &1\\
0&0&0&0  &1\\
\end{smallmatrix}\right),\quad
\sigma_2^\Lambda(1,3)=\left(\begin{smallmatrix}
1& 0&0&0 &0\\
1&-1&0&0 &0\\
1&-2&1&0 &0\\
1&-3&3&-1&0\\
1&-4&6&-4&1\\
\end{smallmatrix}\right),
$$
$$
\sigma_1^\Lambda(1,4)=\left(\begin{smallmatrix} 1&4\alpha_1&6\alpha_2&4&\alpha_1\\
0&\alpha_1&3\alpha_2&3&\alpha_1\\
0&0&\alpha_2&2 &\alpha_1\\
0&0&0&1 &\alpha_1\\
0&0&0&0  &\alpha_1\\
\end{smallmatrix}\right),\quad
\sigma_2^\Lambda(1,4)=\left(\begin{smallmatrix}
\alpha_1& 0&0&0 &0\\
-1&1&0&0 &0\\
\alpha_2&-2\alpha_2&\alpha_2&0 &0\\
-\alpha_1&3\alpha_1&-3\alpha_1&\alpha_1&0\\
1&-4&6&-4&1\\
\end{smallmatrix}\right),
$$
$$
\sigma_1^\Lambda(1,4)=\left(\begin{smallmatrix} 1&4i&-6&-4i&1\\
0&i&-3&-3i&1\\
0&0&-1&-2i &1\\
0&0&0&-i &1\\
0&0&0&0  &1\\
\end{smallmatrix}\right),\quad
\sigma_2^\Lambda(1,2)=\left(\begin{smallmatrix}
1& 0&0&0 &0\\
i&-i&0&0 &0\\
-1&2&-1&0 &0\\
-i&3i&-3i&i&0\\
1&-4&6&-4&1\\
\end{smallmatrix}\right).
$$
{\bf Problem, $n=4$.} To find the eigenvectors
$e_0^{(s),t_0},\,\,e_1^{(s),t_1},\,\,e_2^{(s),t_2}$
 for $\sigma_1^\Lambda(1,4)$ and
$f_0^{(s),\tau_0},\,\,f_1^{(s),\tau_1},\,\,f_2^{(s),\tau_2}$ for
$\sigma_2^\Lambda(1,4)$ in three different cases
$\Lambda_3^{(s)}=\lambda_0{\rm diag }(\alpha^{(s)}_k)_{k=0}^4,\,\,
s=2,3,4$.

Define in the case $\Lambda_3^{(2)}={\rm diag }(1,-1,1-1,1)$
$$
V_4^{(2),t}:=\langle e_0^{(2),t}=(t,1,1,1,2)\rangle,\quad
W_4^{(2),\tau}:=\langle f_0^{(2),\tau}=(2,1,1,1,\tau)\rangle.
$$
$V_4^{(2),t}=W_4^{(3),\tau}\Leftrightarrow t=\tau=2.$

In the case $\Lambda_3^{(3)}={\rm diag
}(1,\alpha_1,\alpha_2,1,\alpha_1)$ we have
$$ e_0^{(3),t}=(t,\sqrt{3}e(0),2e(1),\sqrt{3}e(2),0),\quad
f_0^{(3),\tau}=(0,\sqrt{3}e(0),2e(1),\sqrt{3}e(2),\tau),
$$
hence
$$
V_4^{(3),t}:=\langle e_0^{(3),t} \rangle=W_4^{(3),\tau }:=\langle
f_0^{(3),\tau} \rangle \Leftrightarrow t=\tau=0.
$$
Define in the case $\Lambda_4^{(4)}={\rm diag }(1,i,-1-i,1)$
$$
V_4^{(4),t}:=\langle e_0^{(4),t},e_1^{(4)},e_2^{(4)}\rangle,\quad
W_4^{(4),\tau}:=\langle f_0^{(4),\tau},f_1^{(4)},f_2^{(4)}\rangle.
$$
If we set $e(k):=\exp (2\pi i k/8)$ we get
$$
\begin{smallmatrix}
e_0^{(4),t}=(&t,            &e(0),  &\sqrt{2}e(1),&2e(2),&2\sqrt{2}e(3)&),\\
e_1^{(4)}=  (&2\sqrt{2},    &e(1),  &0,&0,     &0,            &),\\
e_2^{(4)}=  (&3\sqrt{2}e(-1),&3e(0),&\sqrt{2}e(1),&0,&0,   &),\\
e_3^{(4)}=  (&2\sqrt{2}e(0),&3e(0),&2\sqrt{2}e(2),&2e(3),&0   &),\\
f_0^{(4),\tau}=(&\sqrt{2}e(-1),&e(0),&\frac{1}{\sqrt{2}}e(1),&\frac{1}{2}e(2),&\tau&),\\
f_1^{(4)}=(&0,&0,&0,&e(0),&2\sqrt{2}e(1)&),\\
f_2^{(4)}=(&0,&0,    & \sqrt{2}e(1),&3e(2),&3\sqrt{2}e(3)&),\\
f_3^{(4)}=(&0,&2e(0),&2\sqrt{2}e(1),&3e(2),&2\sqrt{2}e(3)&).\\
\end{smallmatrix}
$$
When $V_4^{(4),t}=W_4^{(4),\tau}?$

{\bf Problem, $n=5$.} To find eigenvectors
$e_{0,5}^{(s),t_0},\,\,e_{1,k}^{(s),t_1},\,\,e_{2,5}^{(s),t_2},\,\,e_{3,5}^{(s),t_3},$
for $\sigma_1^\Lambda(1,5)$ and
$f_{0,5}^{(s),\tau_0},\,\,f_{1,5}^{(s),\tau_1},\,\,f_{2,5}^{(s),\tau_2},\,\,f_{3,5}^{(s),\tau_3}$
for $\sigma_2^\Lambda(1,5)$ corresponding to eigenvalues
$\lambda_{k,5}^{(s)}=\alpha_k^{(s)}$ in four different cases
$\Lambda_5^{(s)}={\rm diag }(\alpha^{(s)}_k)_{k=0}^5,\,\,
s=2,3,4,5$.

$$
\Lambda_5^{(2)}={\rm diag }(1,-1,1-1,1,1),\quad
V_5^{(2),t}:=\langle e_0^{(2),t} \rangle=W_5^{(2),\tau }:=\langle
f_0^{(2),\tau}\rangle,
$$
$$
\Lambda_5^{(3)}={\rm diag
}(1,\alpha_1,\alpha_2,1,\alpha_1,\alpha_2),\quad
V_5^{(3),t}:=\langle e_0^{(3),t} \rangle=W_5^{(3),\tau }:=\langle
f_0^{(3),\tau}\rangle,
$$
$$
\Lambda_5^{(4)}={\rm diag }(1,i,-1-i,1,i),\quad
V_5^{(4),t}:=\langle e_0^{(4),t} \rangle=W_5^{(4),\tau }:=\langle
f_0^{(4),\tau}\rangle,
$$
$$
\Lambda_5^{(5)}=\lambda_0{\rm diag }(\alpha^{(5)}_k)_{k=0}^5,\,
$$
$$
V_5^{(5),t}:= \langle
e_0^{(5),t},e_1^{(5)},e_2^{(5)},e_3^{(5)}\rangle,\quad
W_5^{(5),\tau}:=\langle
f_0^{(5),\tau},f_1^{(5)},f_2^{(5)},f_3^{(5)}\rangle.
$$

{\bf Problem, $n$.} To find the eigenvectors
$(e_k^{(s),t_k})_{k=0}^{n-2}$  for $\sigma_1^\Lambda(1,n)$ and
$(f_k^{(s),\tau_k})_{k=0}^{n-2}$ for $\sigma_2^\Lambda(1,n)$
corresponding to eigenvalues
$\lambda_{k,n}^{(s)}=\alpha^{(s)}_k,\,\,\,0\leq k\leq n-2$ in
$n-2$ different cases $\Lambda_n^{(s)}=\lambda_0{\rm diag
}(\alpha^{(s)}_k)_{k=0}^n,\,\, 2\leq s\leq n$.

Define for $2\leq s\leq n-3 $ the subspaces
$$
V_n^{(s),t}: =\langle e_0^{(s),t}\rangle,\quad
W_n^{(s),\tau}:=\langle f_0^{(s),\tau} \rangle,
$$
and
$$
V_n^{(n),t
}:=\langle e_0^{(n),t},\, e_k^{(n)}\mid 1\leq k\leq n-2
\rangle,\quad W_n^{(n),\tau
}:=\langle f_0^{(n),\tau},\,\,f_k^{(n)}\mid 1\leq k\leq n-2
\rangle.
$$
When $V_n^{(s),t}=W_n^{(s),\tau}$ for $0\leq s\leq n-2$ ?
\qed\end{pf}

{\bf Some general formulas.} Let $\Lambda_n(\alpha)={\rm diag}(
\alpha^{n-k})_{k=0}^n$ for some $\alpha\in {\mathbb C}$. We find
the eigenvectors $e_0(\alpha)$ (resp. $f_0(\alpha)$) of the
operator $\sigma_1(1,n)\Lambda_n(\alpha)$ (resp.
$\Lambda_n^\sharp(\alpha)\sigma_1(2,n)$ ) corresponding to the
eigenvalue $\lambda=1$. We would like to find the vectors
$e_0(\alpha)$ (resp. $f_0(\alpha)$) in the following form
$e=(\mu^{n-k})_{k=0}^n$. We have for $n=4$
$$
\left(
\begin{smallmatrix}
1&4&6&4&1\\
0&1&3&3&1\\
0&0&1&2&1\\
0&0&0&1&1\\
0&0&0&0&1\\
\end{smallmatrix}
\right)\Lambda_4(\alpha)
 \left(
\begin{smallmatrix}
\mu^4\\
\mu^3\\
\mu^2\\
\mu^1\\
\mu^0\\
\end{smallmatrix}
\right) = \left(
\begin{smallmatrix}
\mu^4\\
\mu^3\\
\mu^2\\
\mu^1\\
\mu^0\\
\end{smallmatrix}
\right),\,\,\text{or\,\,\,}
\begin{smallmatrix}
(\alpha\mu+1)^4=\mu^4\\
(\alpha\mu+1)^3=\mu^3\\
(\alpha\mu+1)^2=\mu^2\\
(\alpha\mu+1)=\mu\\
(\alpha\mu+1)^0=\mu^0\\
\end{smallmatrix}
$$
so we have $\left(\frac{\alpha\mu+1}{\mu}\right)^k=1$ for
$k=0,...,4$ or   $\alpha\mu+1=\mu$. Finally $\mu=(1-\alpha)^{-1}$
and we get for $n=4$ and for general $n$
$$
e_0(\alpha)=((1-\alpha)^{-(n-k)})_{k=0}^n.
$$
For $f_0(\alpha)$ we get
$$
\Lambda_4^\sharp(\alpha)\left(
\begin{smallmatrix}
1&0&0&0&0\\
-1&1&0&0&0\\
1&-2&1&0&0\\
-1&3&-3&1&0\\
1&-4&6&-4&1\\
\end{smallmatrix}
\right)
 \left(
\begin{smallmatrix}
\mu^4\\
\mu^3\\
\mu^2\\
\mu^1\\
\mu^0\\
\end{smallmatrix}
\right) = \left(
\begin{smallmatrix}
\mu^4\\
\mu^3\\
\mu^2\\
\mu^1\\
\mu^0\\
\end{smallmatrix}
\right),\,\,\text{or\,\,\,}
\begin{smallmatrix}
\alpha^0\mu^4=\mu^4\\
\alpha(-\mu^4+\mu^3)=\mu^3\\
\alpha^2(\mu^4-2\mu^3+\mu^2)=\mu^2\\
\alpha^3(-\mu^4+3\mu^3+3\mu^2+\mu)=\mu\\
\alpha^4(\mu^4-4\mu^3+6\mu^2-4\mu+\mu^0)=\mu^0\\
\end{smallmatrix},\,\,\,
\begin{smallmatrix}
\alpha^0\mu^4=\mu^4\\\
\alpha\mu^3(1-\mu)=\mu^3\\
\alpha^2\mu^2(1-\mu)^2=\mu^2\\
\alpha^3\mu(1-\mu)^3=\mu\\
\alpha^4(1-\mu)^4=\mu^0
\end{smallmatrix}
$$
or $ \alpha^k(1-\mu)^{k}=1,\,\,0\leq k\leq
n,\,\,\mu=1-\alpha^{-1}.$ Finally we get
$$
f_0(\alpha)=((1-\alpha^{-1})^{n-k})_{k=0}^n.
$$
For $n=2,\,\,\Lambda_2=(1,-1,1)$ and $e_2=(2,1,2)$ we have
$\sigma_1^{\Lambda_2}e_2=\sigma_2^{\Lambda_2}e_2=e_2$. Indeed (see
(\ref{reduc}))
\begin{equation}
 \sigma_1^{\Lambda_2}e_2=\left(\begin{smallmatrix}
1&-2&1\\
0&-1&1\\
0&0&1\\
\end{smallmatrix}\right)
\left(\begin{smallmatrix} 2\\
1\\
2\\
\end{smallmatrix}\right)=
\left(\begin{smallmatrix} 2\\
1\\
2\\
\end{smallmatrix}\right),
\quad \sigma_2^{\Lambda_2}e_2= \left(\begin{smallmatrix}
1&0&0\\
1&-1&0\\
1&-2&1
\end{smallmatrix}\right)\left(\begin{smallmatrix} 2\\
1\\
2\\
\end{smallmatrix}\right)=\left(\begin{smallmatrix} 2\\
1\\
2\\
\end{smallmatrix}\right).
\end{equation}
For $n=4,\,\,\Lambda_4=(1,-1,1,-1,1)$ and $e_4=(2,1,1,1,2)$ we
have $\sigma_1^{\Lambda_4}e_4=\sigma_2^{\Lambda_4}e_4=e_4$. Indeed
\begin{equation}
 \sigma_1^{\Lambda_n}e_4=\left(\begin{smallmatrix}
1&-4&6&-4&1\\
0&-1&3&-3&1\\
0&0 &1&-2&1\\
0&0 &0&-1&1\\
0&0 &0& 0&1\\
\end{smallmatrix}\right)
\left(\begin{smallmatrix} 2\\
1\\
1\\
1\\
2\\
\end{smallmatrix}\right)=
\left(\begin{smallmatrix} 2\\
1\\
1\\
1\\
2\\
\end{smallmatrix}\right)
,\quad \sigma_2^{\Lambda_4}e_4= \left(\begin{smallmatrix}
1&0 &0&0 &0\\
1&-1&0&0 &0\\
1&-2&1&0 &0\\
1&-3&3&-1&0\\
1&-4&6&-4&1\\
\end{smallmatrix}\right)\left(\begin{smallmatrix} 2\\
1\\
1\\
1\\
2\\
\end{smallmatrix}\right)=\left(\begin{smallmatrix} 2\\
1\\
1\\
1\\
2\\
\end{smallmatrix}\right).
\end{equation}
In the general case for $n=2m,\,\,\Lambda_{n}={\rm
diag}((-1)^k)_{k=0}^{n}$ and $e_n=(2,1,1,...,1,2)$ we have
$\sigma_1^{\Lambda_n}e_n=\sigma_2^{\Lambda_n}e_n=e_n$. Indeed we
have (see (\ref{si_1(1)})) if $k\not=0$ and $k\not=n$
$$
(\sigma_1^{\Lambda_n}e_n)_k=\sum_{m=k}^{n-r}\sigma_1(1,n)_{km}(e_n)_m=
\sum_{m=k}^{n-r}C_{n-k}^{n-m}(-1)^m+C_{n-k}^0=1,
$$
$$
(\sigma_1^{\Lambda_n}e_n)_0=\sum_{m=k}^{n}\sigma_1(1,n)_{0m}(e_n)_m=
\sum_{m=0}^{n}C_{n}^{n-m}(-1)^m+C_{n}^n+C_{n}^0=2,
$$
$$
(\sigma_1^{\Lambda_n}e_n)_n=C_{0}^02=2.
$$
Since $\sigma_2^{\Lambda_n}=(\sigma_1^{\Lambda_n})^\sharp$ and
$e_n$ is symmetric i.e. $(e_n)_k=(e_n)_{n-k}$ we also conclude
that $\sigma_2^{\Lambda_n}e_n=e_n$.

{\bf Reducibility.} We use the following notations
$$\Lambda_n^{(2)}={\rm
diag}(\alpha_k^{(2)})_{k=0}^n,\quad\alpha_k^{(2)}=\exp (2\pi i
k/2),\quad\sigma_1^{\Lambda_n^{(2)}}=\sigma_1(1,n)\Lambda_n^{(2)},\,\,
$$
$$
\sigma_2^{\Lambda_n^{(2)}}=(\Lambda_n^{(2)})^\sharp\sigma_2(1,n),\,\,\,
{\bf
1}=(1,1,...,1),\,\,\,\delta_0=(1,0,...,0),\,\,\,\delta_n=(0,...,0,1),
$$
$$
e_n^{(2)}={\bf
1}+\delta_0+\delta_n,\,\,\text{for\,\,}n=2m,\text{\,\,and\,\,}e_n^{(2)}={\bf
1}-\delta_0-\delta_n\text{\,\,for\,\,}n=2m+1.
$$
\begin{lem}  For any $n\geq 2$ holds
\begin{equation}
\sigma_1^{\Lambda_n^{(2)}}e_n^{(2)}=\sigma_2^{\Lambda_n^{(2)}}e_n^{(2)}=e_n^{(2)}.
\end{equation}
\end{lem}
\begin{pf} It is sufficient to prove that for $n=2m+1$ operator
$\sigma_1^{\Lambda_n^{(2)}}$ (resp. $\sigma_2^{\Lambda_n^{(2)}}$)
acts as follows (the vectors in the second line are the images of
the corresponding vectors in the first line for example
$\sigma_2^{\Lambda_n^{(2)}}\delta_0=-{\bf 1}$):
\begin{equation}
\label{e_n^1}
\left(\begin{smallmatrix} {\bf 1}&\delta_0&\delta_n\\
-\delta_n&\delta_0&-{\bf 1}
\end{smallmatrix}\right),\,\,
\left(\begin{smallmatrix} {\bf 1}&\delta_0&\delta_n\\
-\delta_0&-{\bf 1}&\delta_0\\
\end{smallmatrix}\right),
\end{equation}
and for $n=2m$ as follows
\begin{equation}
\label{e_n^0} \left(\begin{smallmatrix}
{\bf 1}&\delta_0&\delta_n\\
\delta_n&\delta_0&{\bf 1}
\end{smallmatrix}\right),\,\,
\left(\begin{smallmatrix}
{\bf 1}&\delta_0&\delta_n\\
\delta_0&{\bf 1}&\delta_n\\
\end{smallmatrix}\right).
\end{equation}
Indeed, in this case we get for $n=2m+1$
$$
\sigma_1^{\Lambda_n^{(2)}}e_n^{(2)}=\sigma_1^{\Lambda_n^{(2)}}({\bf
1}-\delta_0-\delta_n)=(-\delta_n-\delta_0+{\bf 1})=e_n^{(2)},
$$
$$
\sigma_2^{\Lambda_n^{(2)}}e_n^{(2)}=\sigma_2^{\Lambda_n^{(2)}}({\bf
1}-\delta_0-\delta_n)=(-\delta_0+{\bf 1}-\delta_n)=e_n^{(2)}.
$$
For $n=2m$ we get
$$
\sigma_1^{\Lambda_n^{(2)}}e_n^{(2)}=\sigma_1^{\Lambda_n^{(2)}}({\bf
1}+\delta_0+\delta_n)=(\delta_n+\delta_0+{\bf 1})=e_n^{(2)}
$$
$$
\sigma_2^{\Lambda_n^{(2)}}e_n^{(2)}=\sigma_2^{\Lambda_n^{(2)}}({\bf
1}+\delta_0+\delta_n)=(\delta_0+{\bf 1}+\delta_n)=e_n^{(2)}.
$$
The proof of (\ref{e_n^1}) and (\ref{e_n^0}) is based on the
identity $\sum_{r=0}^k(-1)^rC_k^r=0$. \qed\end{pf}

\subsection{Counterexamples}
By Theorem 3 and 4 we conclude that for  $n=1$ two matrices
$$
\sigma_1^\Lambda(1,1)=\left(\begin{smallmatrix}
\lambda_0&\lambda_1\\
0&\lambda_1
\end{smallmatrix}\right)\text{\,\,and\,\,}
\sigma_2^\Lambda(1,1)=\left(\begin{smallmatrix}
\lambda_1&0\\
-\lambda_0&\lambda_0
\end{smallmatrix}\right)
$$
with
$\lambda_1/\lambda_0=\alpha,\,\,\alpha^2-\alpha+1=0
$
are operator irreducible but they are not {\bf subspace
irreducible}.

For  $q\not= 1,\,\,\Lambda_2=I,\,\,n=2$ two matrices
$\sigma_1^D(q,2)$ and $\sigma_1^D(q,2)$ for $q=-1$
$$
\sigma_1^D(q,2) =\left(\begin{smallmatrix}
q&1+q&1\\
0&1&1\\
0&0&1
\end{smallmatrix}\right)=
\left(\begin{smallmatrix}
-1&0&1\\
0&1&1\\
0&0&1
\end{smallmatrix}\right)
,\,\, 
\sigma_2^D(q,2)=\left(\begin{smallmatrix}
1&0&0\\
-1&1&0\\
1&-(1+q)&q\\
\end{smallmatrix}\right)=\left(\begin{smallmatrix}
1&0&0\\
-1&1&0\\
1&0&-1\\
\end{smallmatrix}\right).
$$
are {\bf operator irreducible} since the minors
$$
M^{01}_{12}(0_2)=\left|\begin{smallmatrix}
1+q&1\\
1-q&1\\
\end{smallmatrix}\right|=2q,\quad
M^{2}_{2}(0_2)=1-q
$$
can not be zero simultaneously (see section 9.1, case 3), proof of
the Lemma \ref{l.case3}). By Lemma \ref{3)red} they are not {\bf
subspace irreducible} since $(2)_q=1+q=0$ etc.

\subsection{Equivalence}
 {\bf Theorem 5} {\it If two represenations $\sigma^\Lambda(q,n)$ and
$\sigma^{\Lambda'}(q',n)$  are equivalent i.e.
$$
\sigma_i^\Lambda(q,n)C=C\sigma_i^{\Lambda'}(q',n),\,\,i=1,2
$$
for some $C\in {\rm GL}(n+1,{\mathbb C})$ then $q/q'=1$ for $n=2m$
and $(q/q')^2=1$ for $n=2m-1.$ }
\begin{pf} To obtain a criteria of the equivalence it is necessary to study four
cases as in the proof of the Theorem 3 (see Section 9.1)
separately.

To prove the theorem it is sufficient to consider the commutation
relation for some $C\in {\rm GL}(n+1,{\mathbb C})$
$$
\lambda_0\lambda_nS(q)\Lambda_nC=C\lambda_0'\lambda_n'S(q')\Lambda_n'.
$$
\qed\end{pf}
\section{$q$-Pascal's triangle and  Tuba--Wenzl representations}
\begin{pf*}{Proof of the Remark 6.3 (the equivalence of the representations).}
In this section  we index row and columns of the matrix $A\in {\rm
Mat}(n,{\mathbb C})$ starting from $1$: $A=(a_{km})_{1\leq k,m\leq
n}$.  It is easy to see that {\bf for} $n=2$ the equivalence
\begin{equation}
\label{equiv2}
\Lambda^{-1}\sigma_1^{\lambda}\Lambda=\sigma_1^{\Lambda},\text{\,\,and\,\,}
\Lambda^{-1}\sigma_2^{\lambda}\Lambda=\sigma_2^{\Lambda}
\end{equation}
holds. Indeed, we have
$$
\sigma_1^{\lambda}=\left(\begin{smallmatrix}
\lambda_1&\lambda_1\\
0&\lambda_2&
\end{smallmatrix}\right)=\Lambda\sigma_1(1),\,\,
\sigma_1^{\Lambda}=\sigma_1(1)\Lambda=\left(\begin{smallmatrix}
\lambda_1&\lambda_2\\
0&\lambda_2&
\end{smallmatrix}\right),
\text{\,where\,} \sigma_1(1)=\left(\begin{smallmatrix}
1&1\\
0&1
\end{smallmatrix}\right),\,\,\Lambda=\left(\begin{smallmatrix}
\lambda_1&0\\
0&\lambda_2
\end{smallmatrix}\right)
$$
and
$$
\sigma_2^{\lambda}=\left(\begin{smallmatrix}
\lambda_2&0\\
-\lambda_2&\lambda_1
\end{smallmatrix}\right)=\sigma_2(1)\Lambda^\sharp,\,\,\sigma_2^{\Lambda}=
\Lambda^\sharp
\sigma_2(1)=\left(\begin{smallmatrix}
\lambda_2&0\\
-\lambda_1&\lambda_1
\end{smallmatrix}\right),\text{\,where\,} \sigma_2(1)=\left(\begin{smallmatrix}
1&0\\
-1&1
\end{smallmatrix}\right).
$$
Hence
$$
\Lambda^{-1}\sigma_1^{\lambda}\Lambda=\sigma_1^{\Lambda}\text{\,\,\,and\,\,\,}
\Lambda^{\sharp}\sigma_2^{\lambda}(\Lambda^\sharp)^{-1}=\sigma_2^{\Lambda}.
$$
But $\Lambda^{\sharp}=\Lambda^{-1}{\rm det}\Lambda$, which yields
(\ref{equiv2}).

We shall show that {\bf for} $n=3$ and $C={\rm
diag}(1,1,\lambda_3/\lambda_2)$ the equivalence
\begin{equation}
\label{equiv3}
\sigma_1^{\lambda}=C\sigma_1^{\Lambda}C^{-1},\text{\,\,and\,\,}
\sigma_2^{\lambda}=C\sigma_2^{\Lambda}C^{-1}
\end{equation}
holds. Indeed, we have
$$
\sigma_1^{\Lambda}=\sigma_1(q)\Lambda\quad\text{\,and\,}\quad
\sigma_2^{\Lambda}=\Lambda^\sharp(\sigma_1^{-1}(q^{-1}))^\sharp,
$$
$$
\text{where}\quad
q=\frac{\lambda_1\lambda_3}{\lambda_2^2},\quad\sigma_1(q)
=\left(\begin{smallmatrix}
1&1+q&1\\
0&1&1\\
0&0&1
\end{smallmatrix}\right),\quad
\Lambda=\left(\begin{smallmatrix}
\lambda_1&0&0\\
0&\lambda_2&0\\
0&0&\lambda_3
\end{smallmatrix}\right).
$$
To find $\sigma_1^{-1}(q),\,\,\sigma_2(q)$ and $\sigma_2^{-1}(q)$
we use the following formulas. Let $X$ be an upper triangular
matrix of infinite order with units on the principal diagonal
$X=I+x=I+\sum_{k<n}x_{kn}E_{kn}$, where $E_{kn}$ are matrix units
of infinite order. Let us denote by $x_{kn}^{-1}$ the matrix
element of the inverse matrix $X^{-1}$
$$
X^{-1}=(I+x)^{-1}=I+\sum_{k<n}x_{kn}^{-1}E_{kn}.
$$
 Since
$XX^{-1}=X^{-1}X=I$ we have
\begin{equation}
\label{x{kn}(-1).2}
\sum_{r=k}^{n}x_{kr}^{-1}x_{rn}=\sum_{r=k}^{n}x_{kr}x_{rn}^{-1}=\delta_{kn}.
\end{equation}
The following  explicit formula for $x_{kn}^{-1}$ holds (see
\cite{Kos88} formula (4.4)) $$x_{kk+1}^{-1}=-x_{kk+1},$$
\begin{equation}
\label{x{kn}(-1).1}
x_{kn}^{-1}=-x_{kn}+\sum_{r=1}^{n-k-1}(-1)^{r+1}\sum_{k\leq
i_1<i_2<...<i_r\leq n }x_{ki_1}x_{i_1i_2}...x_{i_rn},\quad k<n-1.
\end{equation}
We have
$$
\sigma_1^{-1}(q) =\left(\begin{smallmatrix}
1&-(1+q)&q\\
0&1&-1\\
0&0&1
\end{smallmatrix}\right),
$$
 hence
$$
\sigma_2(q)=(\sigma_1^{-1}(q^{-1}))^\sharp
=\left(\begin{smallmatrix}
      1&0&0\\
-1&1&0\\
 q^{-1}&-(1+q^{-1})&1
\end{smallmatrix}\right),\,\,
\sigma_1^{\Lambda}=\sigma_1(q)\Lambda=\left(\begin{smallmatrix}
\lambda_1&\lambda_1\lambda_3\lambda_2^{-1}+\lambda_2&\lambda_3\\
0&\lambda_2&\lambda_3\\
0&0&\lambda_3
\end{smallmatrix}\right),
$$
$$
\sigma_2^{\Lambda}=\Lambda^\sharp(\sigma_1^{-1}(q^{-1}))^\sharp=
\left(\begin{smallmatrix}
\lambda_3&0&0\\
-\lambda_2&\lambda_2&0\\
\lambda_2^2\lambda_3^{-1}&-\lambda_1-\lambda_2^2\lambda_3^{-1}&\lambda_1
\end{smallmatrix}\right).
$$
We compare $\sigma_1^{\Lambda}$ and $\sigma_2^{\Lambda}$ with the
following expressions (see (\ref{TW3}))
$$
\sigma_1^{ \lambda}=\left(\begin{smallmatrix}
\lambda_1&\lambda_1\lambda_3\lambda_2^{-1}+\lambda_2&\lambda_2\\
0&\lambda_2&\lambda_2\\
0&0&\lambda_3
\end{smallmatrix}\right),\quad
\sigma_2^{\lambda}=\left(\begin{smallmatrix}
\lambda_3&0&0\\
-\lambda_2&\lambda_2&0\\
\lambda_2&-\lambda_1\lambda_3\lambda_2^{-1}-\lambda_2&\lambda_1
\end{smallmatrix}\right).
$$
We have
$$
\sigma_1^{\lambda}\Lambda^{-1}=\left(\begin{smallmatrix}
1&1+q&\frac{\lambda_2}{\lambda_3}\\
0&1&\frac{\lambda_2}{\lambda_3}\\
0&0&1
\end{smallmatrix}\right)=:\sigma_1'(q)\text{\,\,and\,\,}
\sigma_1^{\Lambda}\Lambda^{-1}=\left(\begin{smallmatrix}
1&1+q&1\\
0&1&1\\
0&0&1
\end{smallmatrix}\right)=\sigma_1(q).
$$
We see that $\sigma_1'(q)=C\sigma_1(q)C^{-1}$  where $C={\rm
diag}(1,1,\lambda_3/\lambda_2)$, hence
$$
\sigma_1^{\lambda}\Lambda^{-1}=C\sigma_1^{\Lambda}\Lambda^{-1}C^{-1},\quad
\text{so}\quad
\sigma_1^{\lambda}=C\sigma_1^{\Lambda}\Lambda^{-1}C^{-1}\Lambda=
C\sigma_1^{\Lambda}C^{-1}
$$
since $\Lambda^{-1}C^{-1}\Lambda=C^{-1}$ (both $C$ and $\Lambda$
are diagonal). Further
$$
(\Lambda^\sharp)^{-1}\sigma_1^{\Lambda}=(\sigma_1^{-1}(q^{-1}))^\sharp
=\left(\begin{smallmatrix}
      1&0&0\\
-1&1&0\\
 q^{-1}&-(1+q^{-1})&1
\end{smallmatrix}\right),
$$
and
$$
(\Lambda^\sharp)^{-1}\sigma_1^{\lambda}= \left(\begin{smallmatrix}
      1&0&0\\
-1&1&0\\
\frac{\lambda_2}{\lambda_1}&-\frac{\lambda_3}{\lambda_2}-\frac{\lambda_2}
{\lambda_1}&1
\end{smallmatrix}\right)
=C\left(\begin{smallmatrix}
      1&0&0\\
-1&1&0\\
 q^{-1}&-(1+q^{-1})&1
\end{smallmatrix}\right)C^{-1}.
$$
Thus (\ref{equiv3}) holds.

{\bf For} $n=4$ we get if we put
$q=\left(\frac{\lambda_1\lambda_4}{\lambda_2\lambda_3}\right)^{1/2}=
D^{-1}$
$$
\sigma_1^{\lambda}=\sigma_1^{\Lambda}=\sigma_1(q)\Lambda\quad\text{\,and\,}\quad
\sigma_2^{\lambda}=\sigma_2^{\Lambda}=\Lambda^\sharp(\sigma_1^{-1}(q^{-1}))^\sharp.
$$
Indeed, we have
$$
\sigma_1(q)=\left(\begin{smallmatrix}
1&1+q+q^2&1+q+q^2&1\\
0&1&1+q&1\\
0&0&1&1\\
0&0&0&1
\end{smallmatrix}\right),\quad
\Lambda=\left(\begin{smallmatrix}
\lambda_1&0&0&0\\
0&\lambda_2&0&0\\
0&0&\lambda_3&0\\
0&0&0&\lambda_4
\end{smallmatrix}\right),
$$
hence $\sigma_1^{\Lambda}=\sigma_1(q)\Lambda=\sigma_1^{\lambda}$
(see (\ref{TW41}) and (\ref{TW42})). To find $\sigma_1^{-1}(q)$ we
have by (\ref{x{kn}(-1).2})
$$
 x_{kn}+\sum_{r=k+1}^{n-1}x_{kr}^{-1}x_{rn}+x_{kn}^{-1}=0,\,\,
 x_{kn}^{-1}+\sum_{r=k+1}^{n-1}x_{kr}x_{rn}^{-1}+x_{kn}=0,
$$
$$
x_{kk+1}=-x_{kk+1},\quad
x_{12}^{-1}=-(1+q+q^2),\,x_{23}^{-1}=-(1+q),\quad x_{34}^{-1}=-1,
$$
$$
x_{24}^{-1}= -x_{24}-x_{23}^{-1}x_{34}=-1+(1+q)=q,\,\,
$$
$$
x_{13}^{-1}=-x_{13}-x_{12}^{-1}x_{23}=-(1+q+q^2)+(1+q)(1+q+q^2)
=q(1+q+q^2),
$$
$$
x_{14}^{-1}=-x_{14}-x_{12}^{-1}x_{23}-x_{13}^{-1}x_{34}=
-1+(1+q+q^2)(1+q)-q(1+q+q^2)=-q^3,
$$
(where we have used the notation $x_{km}:=\sigma_1(q)_{km}$),
hence
$$
\sigma_1^{-1}(q)=\left(\begin{smallmatrix}
1&-(1+q+q^2)&q(1+q+q^2)&-q^3\\
0&1&-(1+q)&q\\
0&0&1&-1\\
0&0&0&1
\end{smallmatrix}\right),\,\,\sigma_1^{-1}(q^{-1}) =\left(\begin{smallmatrix}
1&-(1+q^{-1}+q^{-2})&q^{-1}(1+q^{-1}+q^{-2})&-q^{-3}\\
0&1&-(1+q^{-1})&q^{-1}\\
0&0&1&-1\\
0&0&0&1
\end{smallmatrix}\right),
$$
$$
\sigma_2(q)=(\sigma_1^{-1}(q^{-1}))^\sharp =
\left(\begin{smallmatrix}
 1    &0          &0              &0\\
-1    &1         &0                &0\\
q^{-1}&-(1+q^{-1})&1                 &0\\
-q^{-3}&q^{-1}(1+q^{-1}+q^{-2})&-(1+q^{-1}+q^{-2})&1
\end{smallmatrix}\right)
$$
and
$\sigma_2^{\Lambda}=\Lambda^\sharp(\sigma_1^{-1}(q^{-1}))^\sharp=\sigma_2^{\lambda}$
(see (\ref{TW42})).

{\bf For} $n=5$ if we put (see (\ref{La(q)}) and (\ref{cond_q}))
$q^{-3}=\frac{\lambda_2\lambda_4}{\lambda_1\lambda_5},\quad
q^{-4}=\frac{\lambda_3^2}{\lambda_1\lambda_5}\,\,$ we get
$$
\sigma_1^{\Lambda}=\sigma_1(q)\Lambda\quad\text{\,and\,}\quad
\sigma_2^{\Lambda}=\Lambda^\sharp(\sigma_1^{-1}(q^{-1}))^\sharp,
$$
where
$$
\sigma_1(q)=\left(\begin{smallmatrix}
1&(1+q)(1+q^2)&(1+q^2)(1+q+q^2)&(1+q)(1+q^2)&1\\
0&1&1+q+q^2&1+q+q^2&1\\
0&0&1&1+q&1\\
0&0&0&1&1\\
0&0&0&0&1
\end{smallmatrix}\right),\quad
\Lambda=\left(\begin{smallmatrix}
\lambda_1&0&0&0&0\\
0&\lambda_2&0&0&0\\
0&0&\lambda_3&0&0\\
0&0&0&\lambda_4&0\\
0&0&0&0&\lambda_5
\end{smallmatrix}\right),
$$
$$
\sigma_1(q)^{-1}= \left(\begin{smallmatrix}
1&-(1+q)(1+q^2)&q(1+q)(1+q+q^2)&-q^3(1+q)(1+q^2)&q^6\\
0&1&-(1+q+q^2)&q(1+q+q^2)&-q^3\\
0&0&1&-(1+q)&q\\
0&0&0&1&-1\\
0&0&0&0&1
\end{smallmatrix}\right).
$$
Setting
$\gamma=(\lambda_1\lambda_2\lambda_3\lambda_4\lambda_5)^{1/5}$ in
(\ref{TW5})) we have
$$
\sigma_1\mapsto\sigma_1^{ \lambda}=\left(\begin{smallmatrix}
\lambda_1 &
(1+\frac{\gamma^2}{\lambda_2\lambda_4})(\lambda_2+\frac{\gamma^3}{\lambda_3\lambda_4})
&
(\frac{\gamma^2}{\lambda_3}+\lambda_3+\gamma)(1+\frac{\lambda_1\lambda_5}{\gamma^2})&
(1+\frac{\lambda_2\lambda_4}{\gamma^2})
(\lambda_3+\frac{\gamma^3}{\lambda_2\lambda_4})&\frac{\gamma^3}{\lambda_1\lambda_5}
\\
0&\lambda_2&\frac{\gamma^2}{\lambda_3}+\lambda_3+\gamma&
\frac{\gamma^3}{\lambda_1\lambda_5}+\lambda_3+\gamma&\frac{\gamma^3}
{\lambda_1\lambda_5}
\\
0&0&\lambda_3&\frac{\gamma^3}{\lambda_1\lambda_5}+\lambda_3
&\frac{\gamma^3}{\lambda_1\lambda_5}\\
0&0&0&\lambda_4&\lambda_4\\
0&0&0&0&\lambda_5
\end{smallmatrix}\right)=
$$
$$
=\left(\begin{smallmatrix} 1 &
(1+\frac{\gamma^2}{\lambda_2\lambda_4})(1+\frac{\gamma^3}{\lambda_2\lambda_3\lambda_4})
& \left(1 +
\frac{\gamma}{\lambda_3}+\left(\frac{\gamma}{\lambda_3}\right)^2\right)
(1+\frac{\lambda_1\lambda_5}{\gamma^2})&
\frac{\lambda_3}{\lambda_4}(1+\frac{\lambda_2\lambda_4}{\gamma^2})
(1+\frac{\gamma^3}{\lambda_2\lambda_3\lambda_4})&\frac{\gamma^3}{\lambda_1\lambda_5^2}
\\
0&1&1+\frac{\gamma}{\lambda_3}+\left(\frac{\gamma}{\lambda_3}\right)^2
& \frac{\lambda_3}{\lambda_4}\left(1+\frac{\gamma}{\lambda_3}
+\frac{\gamma^3}{\lambda_1\lambda_3\lambda_5}\right)
&\frac{\gamma^3}{\lambda_1\lambda_5^2}
\\
0&0&1&\frac{\lambda_3}{\lambda_4}\left(1+\frac{\gamma^3}{\lambda_1\lambda_3\lambda_5}\right)
&\frac{\gamma^3}{\lambda_1\lambda_5^2}\\
0&0&0&1&\frac{\lambda_4}{\lambda_5}\\
0&0&0&0&1
\end{smallmatrix}\right)\Lambda.
$$
We show the equivalence of our representation with the Tuba-Wenzl
representation by finding  some inversible  matrix $C\in {\rm
Mat}(5,{\mathbb C})$ such that
$$
\sigma_1^\lambda=C^{-1}\sigma_1^\Lambda C.
$$
Indeed, using (\ref{La(q)}) and (\ref{cond_q}) we get
$$
\Lambda(q)= {\rm
diag}(1,\frac{\lambda_2\lambda_4}{\lambda_1\lambda_5},
\frac{\lambda_3^2}{\lambda_1\lambda_5},\frac{\lambda_2\lambda_4}
{\lambda_1\lambda_5},1)={\rm
diag}(1,q^{-3},q^{-4},q^{-3},1),
$$
we get $q^{-3}=\frac{\lambda_2\lambda_4}{\lambda_1\lambda_5},\quad
q^{-4}=\frac{\lambda_3^2}{\lambda_1\lambda_5}\,\,$ (hence
$q^{-1}=\frac{\lambda_3^2}{\lambda_2\lambda_4}$). Recalling that
$\gamma=(\lambda_1\lambda_2\lambda_3\lambda_4\lambda_5)^{1/5}$ we
conclude that
$$
\frac{\gamma^2}{\lambda_2\lambda_4}=\left(\frac{\lambda_1\lambda_5}
{\lambda_2\lambda_4}\right)^{2/5}\left(\frac{\lambda_3^2}
{\lambda_2\lambda_4}\right)^{1/5}=q^{6/5}q^{-1/5}=q,\quad
\frac{\gamma^3}{\lambda_1\lambda_3\lambda_5}=
\frac{\lambda_2\lambda_4}{\gamma^2}=q^{-1},
$$
$$
\frac{\gamma^3}{\lambda_2\lambda_3\lambda_4}=\left(\frac{\lambda_1\lambda_5}
{\lambda_2\lambda_4}\right)^{3/5}\left(\frac{\lambda_2\lambda_4}{\lambda_3^2}
\right)^{1/5}=q^{9/5}q^{1/5}=q^2,\quad\frac{\lambda_1\lambda_5}{\gamma^2}=
\frac{\gamma^3}{\lambda_2\lambda_3\lambda_4}=q^2,
$$
$$
\quad \frac{\gamma}{\lambda_3}=\left(\frac{\lambda_1\lambda_5}
{\lambda_2\lambda_4}\right)^{1/5}\left(\frac{\lambda_2\lambda_4}{\lambda_3^2}
\right)^{2/5}=q^{3/5}q^{2/5}=q,\quad
\frac{\gamma^3}{\lambda_1\lambda_5^2}=\frac{\gamma^3}{\lambda_1\lambda_3\lambda_5}
\frac{\lambda_3}{\lambda_5}=q^{-1}\frac{\lambda_3}{\lambda_5},
$$
hence
$$
\sigma_1^{ \lambda}=\left(\begin{smallmatrix}
1&(1+q)(1+q^2)&(1+q^2)(1+q+q^2)&\frac{\lambda_3}
{\lambda_4}(1+q^{-1})(1+q^2)&q^{-1}\frac{\lambda_3}{\lambda_5}\\
0&1&1+q+q^2&\frac{\lambda_3}{\lambda_4}(1+q+q^{-1})&q^{-1}\frac{\lambda_3}
{\lambda_5}\\
0&0&1&\frac{\lambda_3}{\lambda_4}(1+q^{-1})&q^{-1}\frac{\lambda_3}{\lambda_5}\\
0&0&0&1&\frac{\lambda_4}{\lambda_5}\\
0&0&0&0&1
\end{smallmatrix}\right)\Lambda
$$
$$
=C_4^{-1}(q)\left(\begin{smallmatrix}
1&(1+q)(1+q^2)&(1+q^2)(1+q+q^2)&(1+q)(1+q^2)&q^{-1}\frac{\lambda_3}{\lambda_5}\\
0&1&1+q+q^2&1+q+q^2&q^{-1}\frac{\lambda_3}{\lambda_5}\\
0&0&1&1+q&q^{-1}\frac{\lambda_3}{\lambda_5}\\
0&0&0&1&q^{-1}\frac{\lambda_3}{\lambda_5}\\
0&0&0&0&1
\end{smallmatrix}\right)\Lambda C_4(q)=
$$
$$
C_5^{-1}(q)C_4^{-1}(q)\left(\begin{smallmatrix}
1&(1+q)(1+q^2)&(1+q^2)(1+q+q^2)&(1+q)(1+q^2)&1\\
0&1&1+q+q^2&1+q+q^2&1\\
0&0&1&1+q&1\\
0&0&0&1&1\\
0&0&0&0&1
\end{smallmatrix}\right)\Lambda
C_4(q)C_5(q)=C^{-1}\sigma_1^{\Lambda}C,
$$
where $C=C_4(q)C_5(q)$ and
$$
 C_{4}(q)={\rm
diag}(1,1,1,q^{-1}\frac{\lambda_3}{\lambda_4},1),\,\,\,C_5(q)={\rm
diag}(1,1,1,1,q^{-1}\frac{\lambda_3}{\lambda_5}).
$$
 Finally we have
$$
\sigma_1^{
\lambda}=C^{-1}\sigma_1^{\Lambda}C,\quad\text{where}\quad C={\rm
diag}(1,1,1,q^{-1}\frac{\lambda_3}{\lambda_4},q^{-1}\frac{\lambda_3}{\lambda_5})
$$
and hence $\sigma_1^{\lambda}$ should be as follows
$\sigma_1^{\lambda}=C^{-1}\sigma_2^{\Lambda}C$. \qed\end{pf*}
\section{Representations of $B_3$ and $q-$Pascal triangle }
\begin{pf*}{Proof of Theorem 1.} Let us first consider the case $n=1$.
 We have
$$
\sigma_1^{\Lambda}=\sigma_1(1)\Lambda=\left(\begin{smallmatrix}
\lambda_0&\lambda_1\\
0&\lambda_1&
\end{smallmatrix}\right),
\text{\,where\,} \sigma_1(1)=\left(\begin{smallmatrix}
1&1\\
0&1
\end{smallmatrix}\right),\,\,\Lambda=\left(\begin{smallmatrix}
\lambda_0&0\\
0&\lambda_1
\end{smallmatrix}\right),
$$
$$
\sigma_2^{\Lambda}=\Lambda^\sharp
\sigma_2(1)=\left(\begin{smallmatrix}
\lambda_1&0\\
-\lambda_0&\lambda_0
\end{smallmatrix}\right),\text{\,where\,} \sigma_2(1)=\left(\begin{smallmatrix}
1&0\\
-1&1
\end{smallmatrix}\right),\,\,
\Lambda^\sharp=\left(\begin{smallmatrix}
\lambda_1&0\\
0&\lambda_0
\end{smallmatrix}\right),
$$
hence
$$
\sigma_1^{\Lambda}\sigma_2^{\Lambda}=
\lambda_0\lambda_1\left(\begin{smallmatrix}
0&1\\
-1&1
\end{smallmatrix}\right)=\lambda_0\lambda_1\sigma_1(1)\sigma_2(1),
$$
$$
\sigma_2^{\Lambda}\sigma_1^{\Lambda}= \left(\begin{smallmatrix}
\lambda_0\lambda_1&\lambda_1^2\\
-\lambda_0^2&0
\end{smallmatrix}\right)=\left(\begin{smallmatrix}
\lambda_1&0\\
0&\lambda_0
\end{smallmatrix}\right)
\left(\begin{smallmatrix}
1&1\\
-1&0
\end{smallmatrix}\right)
\left(\begin{smallmatrix}
\lambda_0&0\\
0&\lambda_1
\end{smallmatrix}\right)=\Lambda^\sharp\sigma_1(1)\sigma_2(1)\Lambda,
$$
$$
\sigma_1^{\Lambda}\sigma_2^{\Lambda}\sigma_1^{\Lambda}=
\sigma_2^{\Lambda}\sigma_1^{\Lambda}\sigma_2^{\Lambda}
=\lambda_0\lambda_1\left(\begin{smallmatrix}
0&\lambda_1\\
-\lambda_0&0
\end{smallmatrix}\right)=\lambda_0\lambda_1\left(\begin{smallmatrix}
0&1\\
-1&0
\end{smallmatrix}\right)
\left(\begin{smallmatrix}
\lambda_0&0\\
0&\lambda_1
\end{smallmatrix}\right)=\lambda_0\lambda_1S\Lambda.
$$
This is (\ref{Br_n(q)}) for $n=1$. Let us  show that
(\ref{Br_n(q)}) holds for general $n\in{\mathbb N}$.

We first show that (\ref{Br_n(q)}): $
\sigma_1^{\Lambda}\sigma_2^{\Lambda}\sigma_1^{\Lambda}=
\sigma_2^{\Lambda}\sigma_1^{\Lambda}\sigma_2^{\Lambda}=\lambda_0\lambda_nS(q)\Lambda,
$ is equivalent with
\begin{equation}
\label{Br_n2(q)}
\sigma_1(q)\Lambda(q)\sigma_2(q)=S(q)\sigma_1^{-1}(q),\quad
\sigma_1(q)\Lambda(q)\sigma_2(q)=\sigma_2^{-1}(q)S(q).
\end{equation}
In fact (\ref{Br_n(q)}) is equivalent with
$$
\sigma_1^{\Lambda}\sigma_2^{\Lambda}=
\lambda_0\lambda_nS(q)\Lambda(\sigma_1^{\Lambda})^{-1}\text{\,\and\,\,}
\sigma_2^{\Lambda}\sigma_1^{\Lambda}=\lambda_0\lambda_nS(q)\Lambda
(\sigma_2^{\Lambda})^{-1}.
$$
We  have
\begin{equation}
\label{Br_n3(q)} \sigma_1^{\Lambda}\sigma_2^{\Lambda}=
\lambda_0\lambda_nS(q)\Lambda(\sigma_1^{\Lambda})^{-1}=
\lambda_0\lambda_nS(q)\sigma_1^{-1}(q),
\end{equation}
$$
\sigma_2^{\Lambda}\sigma_1^{\Lambda}=\lambda_0\lambda_nS(q)\Lambda
(\sigma_2^{\Lambda})^{-1}= \lambda_0\lambda_nS(q)\Lambda
\sigma_2(q)^{-1} (\Lambda^{\sharp})^{-1}=
$$
\begin{equation}
\label{Br_n4(q)} \Lambda^{\sharp}S(q)
\sigma_2(q)^{-1}\lambda_0\lambda_n(\Lambda\Lambda^\sharp)^{-1}\Lambda=
 \Lambda^{\sharp}S(q)
\sigma_2(q)^{-1}\Lambda(q)^{-1}\Lambda,
\end{equation}
(where we used the relation
$\Lambda\Lambda^\sharp=\lambda_0\lambda_n\Lambda(q)$, see
(\ref{cond_q}) and $S(q)\Lambda=\Lambda^\sharp S(q)$). On the
other hand we get
\begin{equation}
\label{Br_n5(q)}
\sigma_1^{\Lambda}\sigma_2^{\Lambda}=\sigma_1(q)\Lambda\Lambda^\sharp\sigma_2(q)=
\lambda_0\lambda_n\sigma_1(q)\Lambda(q)\sigma_2(q).
\end{equation}
Comparing (\ref{Br_n3(q)}) with (\ref{Br_n5(q)}) we conclude that
the first equality in (\ref{Br_n(q)}) and the first part of
(\ref{Br_n2(q)}) are equivalent. Further we have
$$
\sigma_2^{\Lambda}\sigma_1^{\Lambda}=\Lambda^\sharp\sigma_2(q)\sigma_1(q)\Lambda.
$$
Comparing (\ref{Br_n4(q)}) with the latter equation  we conclude
that the second equality in (\ref{Br_n(q)}) and the second part of
(\ref{Br_n2(q)}) are equivalent.

To prove (\ref{Br_n2(q)}) {\bf for general} $n\in{\mathbb N}$, we
give in Lemma \ref{l.Hum2} the explicit formulas for
$\sigma_1^{-1}(q),\,\,\sigma_2(q)$ and $\sigma_2^{-1}(q)$ (compare
with Lemma 4.1, Section 3). Let us recall also the notation
(see(\ref{La(q)}))
 $$
q_n=q^{\frac{(n-1)n}{2}},\,\,n\in{\mathbb N}.
 $$
\begin{lem}
\label{l.Hum2} Let the operator
$\sigma_1(q)=(\sigma_1(q)_{km})_{0\leq k,m\leq n}$ be defined by $
\sigma_1(q)_{km}=C_{n-k}^{n-m}(q)$. Then for the operators
$\sigma_1^{-1}(q),\,\, \sigma_2(q)$ and $\sigma_2^{-1}(q)$ we have
respectively
\begin{equation}
\label{si1(q)} \sigma_1(q)_{km}=C_{n-k}^{n-m}(q),\quad
\sigma_1^{-1}(q)_{km}=(-1)^{k+m}q_{m-k}C_{n-k}^{n-m}(q)
\end{equation}
and
\begin{equation}
\label{si2(q)}
\sigma_2(q)_{km}=(-1)^{k+m}q^{-1}_{k-m}C_{k}^{m}(q^{-1}),\quad
\sigma_2^{-1}(q)_{km}=C_{k}^{m}(q^{-1}).
\end{equation}
\end{lem}
\begin{pf}
The equality
$\sigma_1^{-1}(q)_{km}=(-1)^{k+m}q_{m-k}C_{n-k}^{n-m}(q)$ is
equivalent with
\begin{equation}
\label{Bin1(q)}
\sum_{r=k}^{n}\sigma_1(q)_{kr}\sigma_1^{-1}(q)_{rm}
=\sum_{r=k}^{n}C_{n-k}^{n-r}(q) (-1)^{r+m}q_{m-r}
C_{n-r}^{n-m}(q)=\delta_{km},
\end{equation}
and the equality $\sigma_2^{-1}(q)_{km}=C_{k}^{m}(q^{-1})$ is
equivalent with the following
\begin{equation}
\label{Bin2(q)}
\sum_{r=k}^{n}\sigma_2(q)_{kr}\sigma_2^{-1}(q)_{rm}
=\sum_{r=k}^{n}(-1)^{k+r}q_{k-r}^{-1}C_{k}^{r}(q^{-1})
C_{r}^{m}(q^{-1})=\delta_{km}.
\end{equation}
The identities (\ref{Bin1(q)}) and (\ref{Bin2(q)}) hold however by
(\ref{Bin1[q]}) (see the proof in Section 8).
 \qed\end{pf}
We shall prove now (\ref{Br_n2(q)}) for general $n\in{\mathbb
 N}$:
$$
\sigma_1(q)\Lambda(q)\sigma_2(q)=S(q)\sigma_1^{-1}(q),
$$
\begin{equation}
\label{Bin(q)1}
\text{i.e\,\,}\sum_{r=0}^{n}\sigma_1(q)_{kr}\Lambda(q)_{rr}\sigma_2(q)_{rm}=
S(q)_{k,n-k}
\sigma_1^{-1}(q)_{n-k,m}
\end{equation}
and
$$
\sigma_1(q)\Lambda(q)\sigma_2(q)=\sigma_2^{-1}(q)S(q),
$$
\begin{equation}
\label{Bin(q)2}
\text{i.e\,\,}\sum_{r=k}^{n}\sigma_1(q)_{kr}\Lambda(q)_{rr}\sigma_2(q)_{rm}
= \sigma_2^{-1}(q)_{k,n-m}S(q)_{n-m,m}.
\end{equation}
 Using (\ref{La(q)}) and (\ref{S(q)}) we have
$$
S(q)=(S(q)_{km})_{0\leq k,m\leq n},\quad
S(q)_{km}=q^{-1}_k(-1)^k\delta_{k+m,n},
$$
$$
\Lambda(q)={\rm diag}\left(q_{rn}
\right)_{r=0}^n,\text{\,\,where\,\,}
q_{rn}^{-1}:=\frac{q_n}{q_rq_{n-r}}.
$$
Then by (\ref{si1(q)}), (\ref{si2(q)}) we get
\begin{gather}
\label{Bin(q)3}
\sum_{r=0}^{n}\sigma_1(q)_{kr}\Lambda(q)_{rr}\sigma_2(q)_{rm}=
\sum_{r=0}^{n}C_{n-k}^{n-r}(q)\frac{q_{r}q_{n-r}}{q_{n}}
(-1)^{r+m} q^{-1}_{r-m} C_{r}^{m}(q^{-1}),\\
(S(q)\sigma_1^{-1}(q))_{km}=S(q)_{k,n-k}
\sigma_1^{-1}(q)_{n-k,m}=q^{-1}_{k}(-1)^{k}(-1)^{n-k+m}q_{m-n+k}
C_{k}^{n-m}(q)\nonumber\\
\label{Bin(q)4} =
(-1)^{n+m}\frac{q_{m-n+k}}{q_{k}}C_{k}^{n-m}(q),\\
\label{Bin(q)5}
 (\sigma_2^{-1}(q)S(q))_{km}=
\sigma_2^{-1}(q)_{k,n-m}S(q)_{n-m,m}=C_{k}^{n-m}(q^{-1})q_{n-m}^{-1}(-1)^{n-m}.
\end{gather}
Using (\ref{Bin(q)3}), (\ref{Bin(q)4}) and (\ref{Bin2[q]}) we get
(\ref{Bin(q)1}). To prove that
$(S(q)\sigma_1^{-1}(q))_{km}=(\sigma_2^{-1}(q)S(q))_{km}$ it is
sufficient to show that
\begin{equation}\label{C(q{-1})}
C_{n}^{k}(q)=q^{-1}_{kn}C_{n}^{k}(q^{-1})=\frac{q_n}{q_kq_{n-k}}C_{n}^{k}(q^{-1})
\quad
0\leq k\leq n.
\end{equation}
Indeed if (\ref{C(q{-1})}) holds we have
$$
C_{k}^{n-m}(q)=q^{-1}_{k,n-m}C_{k}^{n-m}(q^{-1})=\frac{q_k}{q_{m-n+k}q_{n-m}}C_{k}^{n-m}(q^{-1}).
$$
Comparing (\ref{Bin(q)4}) and (\ref{Bin(q)5}) we conclude that
$(S(q)\sigma_1^{-1}(q))_{km}=(\sigma_2^{-1}(q)S(q))_{km}$.

To prove (\ref{C(q{-1})}) it is sufficient to use the definition
(\ref{C_n^k(q)}) of $C_n^k(q)$
$$
C_n^k(q)=\frac{(1-q)(1-q^2)...(1-q^n)}{(1-q)(1-q^2)...
(1-q^k)(1-q)(1-q^2)...(1-q^{n-k})},
$$
$$
C_n^k(q^{-1})=\frac{(1-q^{-1})(1-q^{-2})...(1-q^{-n})}{(1-q^{-1})
(1-q^{-2})...(1-q^{-k})
(1-q^{-1})(1-q^{-2})...(1-q^{-(n-k)})},
 $$
 recall that $q_k=q^{\frac{(k-1)k}{2}}=q^{1+2+...+k-1}$ and observe
 that
 $\frac{q_{n+1}}{q_{k+1}q_{n+1-k}}=\frac{q_n}{q_kq_{n-k}}.$ Hence
 (\ref{Bin(q)1}) and (\ref{Bin(q)2}) are both proven, which
 implies (\ref{Br_n2(q)}).
\qed\end{pf*}
\begin{rem}
We illustrate the identity (\ref{Br_n2(q)}) and Lemma \ref{l.Hum2}
for $n=2,3,4$.
\end{rem}
 {\bf For} $n=2$ we have by (\ref{si_1(q)}),
(\ref{si_2(q)}) and (\ref{S(q)})
$$
\sigma_1(q) = \left(\begin{smallmatrix}
1&1+q&1\\
0&1&1\\
0&0&1
\end{smallmatrix}\right),\quad
\sigma_1^{-1}(q) =\left(\begin{smallmatrix}
1&-(1+q)&q\\
0&1&-1\\
0&0&1
\end{smallmatrix}\right),\quad S(q)=\left(\begin{smallmatrix}
0&0&1\\
0&-1&0\\
q^{-1}&0&0
\end{smallmatrix}\right),
$$
$$
\sigma_2(q) =\left(\begin{smallmatrix}
      1&0&0\\
-1&1&0\\
 q^{-1}&-(1+q^{-1})&1
\end{smallmatrix}\right),\quad\sigma_2^{-1}(q)=
\left(\begin{smallmatrix}
      1&0&0\\
1&1&0\\
1&(1+q^{-1})&1
\end{smallmatrix}\right),\quad
\Lambda(q)=\left(\begin{smallmatrix}
1&0&0\\
0&q^{-1}&0\\
0&0&1
\end{smallmatrix}\right).
$$
Finally we get
$$
\sigma_1(q)\Lambda(q)\sigma_2(q) = \left(\begin{smallmatrix}
1&1+q&1\\
0&1&1\\
0&0&1
\end{smallmatrix}\right)
\left(\begin{smallmatrix}
1&0&1\\
0&q^{-1}&0\\
0&0&1
\end{smallmatrix}\right)
\left(\begin{smallmatrix}
      1&0&0\\
-1&1&0\\
 q^{-1}&-(1+q^{-1})&1
\end{smallmatrix}\right)=
\left(\begin{smallmatrix} 0&0&1\\
0&-1&1\\
q^{-1}&-(1+q^{-1})&1
\end{smallmatrix}\right),
$$

$$
S(q)\sigma_1^{-1}(q)=\left(\begin{smallmatrix}
0&0&1\\
0&-1&0\\
q^{-1}&0&0
\end{smallmatrix}\right)\left(\begin{smallmatrix}
1&-(1+q)&q\\
0&1&-1\\
0&0&1
\end{smallmatrix}\right)=\left(\begin{smallmatrix} 0&0&1\\
0&-1&1\\
q^{-1}&-(1+q^{-1})&1
\end{smallmatrix}\right),
$$
and
$$
\sigma_2(q)^{-1}S(q)= \left(\begin{smallmatrix}
      1&0&0\\
1&1&0\\
1&(1+q^{-1})&1
\end{smallmatrix}\right)
 \left(\begin{smallmatrix}
0&0&1\\
0&-1&0\\
q^{-1}&0&0
\end{smallmatrix}\right)=\left(\begin{smallmatrix} 0&0&1\\
0&-1&1\\
q^{-1}&-(1+q^{-1})&1
\end{smallmatrix}\right).
$$
This is (\ref{Br_n2(q)}) for $n=2$.

{\bf For} $n=3$ we have by (\ref{si_1(q)}), (\ref{si_2(q)}) and
(\ref{S(q)})
$$
\sigma_1(q)=\left(\begin{smallmatrix}
1&(1+q+q^2)&(1+q+q^2)&1\\
0&1&(1+q)&1\\
0&0&1&1\\
0&0&0&1
\end{smallmatrix}\right),\quad
\sigma_1^{-1}(q)=\left(\begin{smallmatrix}
1&-(1+q+q^2)&q(1+q+q^2)&-q^3\\
0&1&-(1+q)&q\\
0&0&1&-1\\
0&0&0&1
\end{smallmatrix}\right),
$$
$$
\quad\sigma_1^{-1}(q^{-1}) =\left(\begin{smallmatrix}
1&-(1+q^{-1}+q^{-2})&q^{-1}(1+q^{-1}+q^{-2})&-q^{-3}\\
0&1&-(1+q^{-1})&q^{-1}\\
0&0&1&-1\\
0&0&0&1
\end{smallmatrix}\right),\quad
S(q)=\left(\begin{smallmatrix}
0     &0     &0 &1\\
0     &0     &-1&0\\
0     &q^{-1}&0 &0\\
q^{-3}&0     &0 &0
\end{smallmatrix}\right),
$$
$$
\sigma_2(q)=(\sigma_1^{-1}(q^{-1}))^\sharp =
\left(\begin{smallmatrix}
 1    &0          &0              &0\\
-1    &1         &0                &0\\
q^{-1}&-(1+q^{-1})&1                 &0\\
-q^{-3}&q^{-1}(1+q^{-1}+q^{-2})&-(1+q^{-1}+q^{-2})&1
\end{smallmatrix}\right),
$$
$$
\sigma_2^{-1}(q)= \left(\begin{smallmatrix}
 1    &0          &0              &0\\
1    &1         &0                &0\\
1&(1+q^{-1})&1                 &0\\
1&(1+q^{-1}+q^{-2})&(1+q^{-1}+q^{-2})&1
\end{smallmatrix}\right),\quad \Lambda(q)=\left(\begin{smallmatrix}
1&0&0&0\\
0&q^{-2}&0&0\\
0&0&q^{-2}&0\\
0&0&0&1
\end{smallmatrix}\right).
$$

We verify that (\ref{Br_n2(q)}) holds, moreover that
$$
\sigma_1(q)\Lambda(q)\sigma_2(q)=
S(q)\sigma_1^{-1}(q)=\left(\begin{smallmatrix}
0&0&0&1\\
0&0&-1&1\\
0&q^{-1}&-(1+q^{-1})&1\\
-q^{-3}&q^{-1} (1+q^{-1}+q^{-2})&-(1+q^{-1}+q^{-2})&1
\end{smallmatrix}\right).
$$
 Indeed we have
$$
\sigma_1(q)\Lambda(q)\sigma_2(q) = \left(\begin{smallmatrix}
1&(1+q+q^2)&(1+q+q^2)&1\\
0&1&(1+q)&1\\
0&0&1&1\\
0&0&0&1
\end{smallmatrix}\right)\left(\begin{smallmatrix}
1&0&0&0\\
0&q^{-2}&0&0\\
0&0&q^{-2}&0\\
0&0&0&1
\end{smallmatrix}\right)\times
$$
$$
 \left(\begin{smallmatrix}
 1    &0          &0              &0\\
-1    &1         &0                &0\\
q^{-1}&-(1+q^{-1})&1                 &0\\
-q^{-3}&q^{-1}(1+q^{-1}+q^{-2})&-(1+q^{-1}+q^{-2})&1
\end{smallmatrix}\right)=\left(\begin{smallmatrix}
0&0&0&1\\
0&0&-1&1\\
0&q^{-1}&-(1+q^{-1})&1\\
-q^{-3}&q^{-1} (1+q^{-1}+q^{-2})&-(1+q^{-1}+q^{-2})&1
\end{smallmatrix}\right).
$$

$$
S(q)\sigma_1^{-1}(q)= \left(\begin{smallmatrix}
0&0&0&1\\
0&0&-1&0\\
0&q^{-1}&0&0\\
-q^{-3}&0&0&0
\end{smallmatrix}\right)\left(\begin{smallmatrix}
1&-(1+q+q^2)&q(1+q+q^2)&-q^3\\
0&1&-(1+q)&q\\
0&0&1&-1\\
0&0&0&1
\end{smallmatrix}\right)
$$
$$
=\left(\begin{smallmatrix}
0&0&0&1\\
0&0&-1&1\\
0&q^{-1}&-(1+q^{-1})&1\\
-q^{-3}&q^{-1} (1+q^{-1}+q^{-2})&-(1+q^{-1}+q^{-2})&1
\end{smallmatrix}\right)=\sigma_1(q)\Lambda(q)\sigma_2(q),
$$
and
$$
\sigma_2(q)^{-1}S(q)= \left(\begin{smallmatrix}
1&0&0&0\\
1&1&0&0\\
1&(1+q^{-1})&1&0\\
1&(1+q^{-1}+q^{-2})&(1+q^{-1}+q^{-2})&1
\end{smallmatrix}\right)\left(\begin{smallmatrix}
0&0&0&1\\
0&0&-1&0\\
0&q^{-1}&0&0\\
-q^{-3}&0&0&0
\end{smallmatrix}\right)
$$
$$
=\left(\begin{smallmatrix}
0&0&0&1\\
0&0&-1&1\\
0&q^{-1}&-(1+q^{-1})&1\\
-q^{-3}&q^{-1} (1+q^{-1}+q^{-2})&-(1+q^{-1}+q^{-2})&1
\end{smallmatrix}\right)=\sigma_1(q)\Lambda(q)\sigma_2(q).
$$
This is (\ref{Br_n2(q)})  for n=3. {\bf For} $n=4$ we have
$$
\sigma_1(q)=\left(\begin{smallmatrix}
1&(1+q)(1+q^2)&(1+q+q^2)(1+q^2)&(1+q)(1+q^2)&1\\
0&1&(1+q+q^2)&(1+q+q^2)&1\\
0&0&1&(1+q)&1\\
0&0&0&1&1\\
0&0&0&0&1
\end{smallmatrix}\right),\quad \Lambda(q)=\left(\begin{smallmatrix}
1&0     &0     &0    &0\\
0&q^{-3}&0     &0    &0\\
0&0     &q^{-4}&0    &0\\
0&0     &0     &q^{-3}&0\\
0&0     &0     &0    &1
\end{smallmatrix}\right)
$$
$$
\sigma_1^{-1}(q^{-1}) =\left(\begin{smallmatrix}
1&-(1+q^{-1})(1+q^{-2})&q^{-1}(1+q^{-1}+q^{-2})(1+q^{-2})&-q^{-3}(1+q^{-1})(1+q^{-2})&q^{-6}\\
0&1&-(1+q^{-1}+q^{-2})&q^{-1}(1+q^{-1}+q^{-2})&-q^{-3}\\
0&0&1&-(1+q^{-1})&q^{-1}\\
0&0&0&1&-1\\
0&0&0&0&1
\end{smallmatrix}\right),
$$
$\sigma_2(q)=(\sigma_1^{-1}(q^{-1}))^\sharp =$
$$
\left(\begin{smallmatrix}
 1    &0          &0              &0&0\\
-1    &1         &0                &0&0\\
q^{-1}&-(1+q^{-1})&1                 &0&0\\
-q^{-3}&q^{-1}(1+q^{-1}+q^{-2})&-(1+q^{-1}+q^{-2})&1&0\\
q^{-6}&-q^{-3}(1+q^{-1})(1+q^{-2})&q^{-1}(1+q^{-1}+q^{-2})
(1+q^{-2})&-(1+q^{-1})(1+q^{-2})&1
\end{smallmatrix}\right),
$$
$$
\sigma_2^{-1}(q)= \left(\begin{smallmatrix}
 1    &0          &0              &0&0\\
1    &1         &0                &0&0\\
1&(1+q^{-1})&1                 &0&0\\
1&(1+q^{-1}+q^{-2})&(1+q^{-1}+q^{-2})&1&0\\
1&(1+q^{-1})(1+q^{-2})&(1+q^{-1}+q^{-2})(1+q^{-2})&(1+q^{-1})(1+q^{-2})&1
\end{smallmatrix}\right),
$$
hence
$\sigma_1(q)\Lambda(q)\sigma_2(q)=S(q)\sigma_1^{-1}(q)=\sigma_2^{-1}(q)S(q)$.

\section{Combinatorial identities for $q-$binomial coefficients}
\begin{lem}
The following identities hold
\begin{equation}
 \label{Bin1[q]}
\sum_{i=0}^{n}C_{m}^{i}(q) (-1)^{i+j}q_{i-j}
C_{i}^{j}(q)=\sum_{i=0}^{n}(-1)^{i+m}q_{m-i}C_{m}^{i}(q)
C_{i}^{j}(q)=\delta_{mj},
\end{equation}
\begin{equation}
 \label{Bin2[q]}
\sum_{r=0}^{n}C_{n-k}^{n-r}(q)\frac{q_{r}q_{n-r}}{q_{n}}
(-1)^{n-r} q^{-1}_{r-m}
C_{r}^{m}(q^{-1})=\frac{q_{k-(n-m)}}{q_{k}}C_{k}^{n-m}(q).
\end{equation}
\end{lem}
\begin{rem}
For $q=1$ (\ref{Bin1[q]}) and (\ref{Bin2[q]}) reduce to the well
known identities (\ref{Bin1}) and (\ref{Bin2}) (see
\cite[p.4]{Rio68} and \cite[p.8 eq. (5)]{Rio68}):
 \begin{equation}\label{Bin1}
\sum_{i=0}^n(-1)^{i+m}\left(\begin{smallmatrix}
 m\\
 i
\end{smallmatrix}\right)\left(\begin{smallmatrix}
 i\\
 j
\end{smallmatrix}\right)=
\sum_{i=0}^n(-1)^{i+j}\left(\begin{smallmatrix}
 m\\
 i
\end{smallmatrix}\right)\left(\begin{smallmatrix}
 i\\
 j
\end{smallmatrix}\right)=
\delta_{mj},
\end{equation}
\begin{equation}
\label{Bin2}
 \sum_{i=0}^n(-1)^i\left(\begin{smallmatrix}
 m\\
 i
\end{smallmatrix}\right)
\left(\begin{smallmatrix}
 n-i\\
 n-j
\end{smallmatrix}\right)=
\sum_{i=0}^n(-1)^i\left(\begin{smallmatrix}
 m\\
 i
\end{smallmatrix}\right)
\left(\begin{smallmatrix}
 n-i\\
 j-i
\end{smallmatrix}\right)=
\left(\begin{smallmatrix}
 n-m\\
 j
\end{smallmatrix}\right).
\end{equation}
\end{rem}
\begin{pf}
We prove the identities by induction. For $n=0$ we have in both
cases $1=1$. Let (\ref{Bin1[q]}) holds for $n\in{\mathbb N}.$ We
prove that this holds then for $n+1$ i.e.
$$
\sum_{i=0}^{n+1}C_{m}^{i}(q) (-1)^{i+j}q_{i-j} C_{i}^{j}(q)
=\delta_{mj}\quad 0\leq m,j\leq n+1.
$$
For $0\leq m,j\leq n$ this hold by the assumption. It is
sufficient to consider $m=n+1$. We have by (\ref{GP1})
$$
\sum_{i=0}^{n+1}C_{n+1}^{i}(q) (-1)^{i+j}q_{i-j}C_{i}^{j}(q)=
\sum_{i=0}^{n+1}(-1)^{i+j}q_{i-j}\left(C_{n}^{i-1}(q)+q^iC_{n}^{i}(q)\right)
C_{i}^{j}(q)=
$$
$$
\sum_{i=0}^{n+1}(-1)^{i+j}q_{i-j}C_{n}^{i-1}(q)C_{i}^{j}(q)+\sum_{i=0}^{n+1}C_{n}^{i}(q)
(-1)^{i+j}q^iq_{i-j}C_{i}^{j}(q)=
$$
$$
\sum_{i=0}^{n+1}(-1)^{i+j}q_{i-j}C_{n}^{i-1}(q)\left(C_{i-1}^{j-1}(q)+q^jC_{i-1}^{j}
(q)\right)+
\sum_{i=0}^{n+1}C_{n}^{i}(q)
(-1)^{i+j}q^iq_{i-j}C_{i}^{j}(q)=
$$
$$
\sum_{i=0}^{n+1}(-1)^{i+j}q_{i-j}C_{n}^{i-1}(q) C_{i-1}^{j-1}(q)+
$$
$$
\sum_{i=0}^{n+1}(-1)^{i+j}q_{i-j}C_{n}^{i-1}(q) q^jC_{i-1}^{j}(q)+
\sum_{i=0}^{n+1}C_{n}^{i}(q) (-1)^{i+j}q^iq_{i-j}C_{i}^{j}(q)=
\delta_{n,j-1}=\delta_{n+1,j}
$$
since the sum of the last two terms gives $0$ by
$q^jq_{i+1-j}=q^iq_{i-j}$. The latter relation follows from
$q_{n+1}=q^nq_n$.

Let us suppose that (\ref{Bin2[q]}) is true for $n\in{\mathbb N}$.
We prove that then this holds for $n+1$ i.e.
$$
\sum_{r=0}^{n+1}(-1)^{n+1-r}C_{n+1-k}^{n+1-r}(q)\frac{q_{r}q_{n+1-r}}{q_{n+1}}
 q^{-1}_{r-m}C_{r}^{m}(q^{-1})=\frac{q_{k-(n+1-m)}}{q_{k}}C_{k}^{n+1-m}(q).
$$
Indeed, by (\ref{GP1}) the left hand side of the latter equation
is equal to
$$\sum_{r=0}^{n+1}(a_r+b_r)(c_r+d_r)=\sum_{r=0}^{n+1}[a_r(c_r+d_r)+b_rc_r+b_rd_r]:=$$
$$
\sum_{r=0}^{n+1}(-1)^{n+1-r}\left(C_{n-k}^{n-r}(q)+q^{n+1-r}C_{n-k}^{n+1-r}(q)\right)
\frac{q_{r}q_{n+1-r}}{q_{n+1}}  q^{-1}_{r-m}\times
$$
$$
\left(C_{r-1}^{m-1}(q^{-1})+q^{-m}C_{r-1}^{m}(q^{-1})\right)
=\sum_{r=0}^{n+1}(-1)^{n+1-r}C_{n-k}^{n-r}(q)\frac{q_{r}q_{n+1-r}}{q_{n+1}}
 q^{-1}_{r-m} C_{r}^{m}(q^{-1})
$$
$$
+\sum_{r=0}^{n+1}(-1)^{n+1-r}q^{n+1-r}C_{n-k}^{n+1-r}(q)
\frac{q_{r}q_{n+1-r}}{q_{n+1}} q^{-1}_{r-m}C_{r-1}^{m-1}(q^{-1})
$$
$$
+ \sum_{r=0}^{n+1}(-1)^{n+1-r}q^{n+1-r}C_{n-k}^{n+1-r}(q)
\frac{q_{r}q_{n+1-r}}{q_{n+1}} q^{-1}_{r-m}C_{r-1}^{m}(q^{-1}).
$$
Since $q_r=q_{r-1}q^{r-1}$ and $q_{n+1}=q_nq^n$  we have
\begin{equation}
\label{[1]}
q^{n+1-r}\frac{q_{r}q_{n+1-r}}{q_{n+1}}=\frac{q_{r-1}q_{n-(r-1)}}{q_{n}}=
\frac{q_{s}q_{n-s}}{q_{n}}.
\end{equation}
Setting $s=r-1$  we get by the  assumption of the induction
(\ref{Bin2[q]})
$$
\sum_{r=0}^{n+1}b_rc_r=
\sum_{s=0}^{n}(-1)^{n-s}C_{n-k}^{n-s}(q)\frac{q_{s}q_{n-s}}{q_{n}}
 q^{-1}_{s-(m-1)}
 C_{s}^{m-1}(q^{-1})
$$
$$
 =\frac{q_{k-[n-(m-1)]}}{q_k}C_{k}^{n-(m-1)}(q)=
\frac{q_{k-(n+1-m)}}{q_k}C_{k}^{n+1-m}(q).
$$
We prove that
$\sum_{r=0}^{n+1}\left[a_r(c_r+d_r)+b_rd_r\right]=0$. Indeed
$\sum_{r=0}^{n+1}a_r(c_r+d_r)=$
$$
\sum_{r=0}^{n+1}(-1)^{n+1-r}C_{n-k}^{n-r}(q)\frac{q_{r}q_{n+1-r}}{q_{n+1}}
 q^{-1}_{r-m} C_{r}^{m}(q^{-1})=
$$
$$
-\sum_{r=0}^{n}(-1)^{n-r}C_{n-k}^{n-r}(q)\frac{q_{r}q_{n-r}}{q_{n}}
 q^{-r}q^{-1}_{r-m} C_{r}^{m}(q^{-1}).
$$
If we set $s=r-1$,  use (\ref{[1]}) and $
q_{s+1-m}q^{m}=q_{s-m}q^{s-m}q^m=q_{s-m}q^s$ we get
$$
\sum_{r=0}^{n+1}b_rd_r=\sum_{r=0}^{n+1}(-1)^{n+1-r}q^{n+1-r}C_{n-k}^{n+1-r}(q)
\frac{q_{r}q_{n+1-r}}{q_{n+1}}
q^{-1}_{r-m}q^{-m}C_{r-1}^{m}(q^{-1})=
$$
$$
\sum_{r=0}^{n}(-1)^{n-s}C_{n-k}^{n-s}(q)
\frac{q_{s}q_{n-s}}{q_{n}} q^{-1}_{s+1-m}q^{-m}C_{r-1}^{m}(q^{-1})
$$
hence $\sum_{r=0}^{n+1}\left[a_r(c_r+d_r)+b_rd_r\right]=0$ and
(\ref{Bin2[q]}) is proven for general $n\in{\mathbb N}$.
\qed\end{pf}

\section{A $q-$analogue of the results of E.~Ferrand}
 Denote by $\Phi(q)=\Phi_n(q)$ the endomorphism of the
space ${\mathbb C}^n[X]$ of polynomials of degree $n$ with complex
coefficients, which maps a polynomial $p(x)$ to the polynomial
$p_q(X+1),$ where for $p(X)=\sum_{k=0}^na_kX^k$ we define
$$
p(X)\stackrel{\Phi_n(q)}{\mapsto}p_q(1+X):=\sum_{k=0}^na_k(1+X)^k_q.
$$
 Denote by $\Psi(q)=\Psi_n(q)$ the endomorphism of ${\mathbb
C}^n[X]$ which maps  a polynomial $p(X)$ to the following
polynomial
$$
p(X)\stackrel{\Psi_n(q)}{\mapsto}
\sum_{k=0}^na_kq_k(1-X)^{n-k}_{q^{-1}}X^k,
$$
compare with the expression for $\Psi$ (see Section 3)
$$
p(X)\stackrel{\Psi}{\mapsto}(1-X)^np\left(\frac{X}{1-X}\right)=
\sum_{k=0}^na_k(1-X)^{n-k}X^k.
$$
\begin{thm} The endomorphisms
$\Phi(q)$ and $\Psi(q)$ satisfy  braid-like relation
$\Phi(q)\Psi(q)\Phi(q)=\Psi(q)\Phi(q)\Psi(q)$.
\end{thm}
\begin{pf} We have by (\ref{GP2})
$$
X^k\stackrel{\Phi(q)}{\mapsto}
(1+X)^k_q=(1+X)(1+qX)...(1+q^{k-1}X)
=\sum_{r=0}^kq^{r(r-1)/2}C_k^r(q)x^r,
$$
hence $ \Phi_{rk}(q)=q^{r(r-1)/2}C_k^r(q)=q_rC_k^r(q)$ and by
(\ref{si_1(q)}) and (\ref{D_n(q)}) we conclude that
\begin{equation}
\label{Phi_n(q)}
\Phi(q)=\Phi_n(q)=D_n(q)\sigma_1^s(q)=(\sigma_1(q)D_n^\sharp(q))^s.
\end{equation}
Indeed $\sigma_1(q)_{km}=C_{n-k}^{n-m}(q)$ hence
$\sigma_1^s(q)_{km}=C_{m}^{k}(q)$ (we recall that
$a^s_{ij}=a_{n-j,n-i}$).
For the operator $\Psi_n(q)$ we get
\begin{equation}
\label{Fer3(q)}
X^k\stackrel{\Psi(q)}{\mapsto}q_{n-k}(1-X)^{n-k}_{q^{-1}}X^k
\end{equation}
$$
=q_{n-k}\sum_{s=0}^{n-k}(-1)^sq_s^{-1}C_{n-k}^s(q^{-1})X^sX^k
=q_{n-k}\sum_{r=k}^{n}(-1)^{r+k}q_{r-k}^{-1}C_{n-k}^{r-k}(q^{-1})X^r,
$$
hence
$\Psi_{rk}(q)=(-1)^{r+k}q_{n-k}q_{r-k}^{-1}C_{n-k}^{r-k}(q^{-1})$
and by (\ref{si_2(q)}) and (\ref{D_n(q)}) we conclude that
\begin{equation}
\label{Psi_n(q)} \Psi(q)=\Psi_n(q)=\sigma_2^s(q)D_n^s(q)=
(D_n(q)\sigma_2(q))^s.
\end{equation}
Indeed $\sigma_2(q)_{rk}=
(-1)^{r+k}q_{r-k}^{-1}C_{r}^{k}(q^{-1}),$ hence
$\sigma_2^s(q)_{rk}= (-1)^{r+k}q_{r-k}^{-1}C_{n-k}^{n-r}(q^{-1})$,
 and $D_n^s(q)={\rm
diag}(q_{n-r})_{r=0}^n$. We note that
\begin{equation}\label{Phi^{-1}(q)}
(\Phi^{-1}(q)p)(X)=p_{q^{-1}}(X-1).
\end{equation}
To finish the proof we use Remark 4.4, Section 3, representation
(\ref{RepD(q)}). \qed\end{pf}

In the particular cases for $n=2$ and $n=3$ we get
\begin{equation}
\label{Phi_2(q)}
\Phi_2(q)=D_2(q)\sigma_1^s(q)=\left(\begin{smallmatrix}
1&0&0\\
0&1&0\\
0&0&q
\end{smallmatrix}\right)\left(\begin{smallmatrix}
1&1&1\\
0&1&1+q\\
0&0&1
\end{smallmatrix}\right)=
\left(\begin{smallmatrix}
1&1&1\\
0&1&1+q\\
0&0&q
\end{smallmatrix}\right),
\end{equation}
\begin{equation}
\label{Phi_3(q)} \Phi_3(q)=D_3(q)\sigma_1^s(q) =
\left(\begin{smallmatrix}
1&0&0&0\\
0&1&0&0\\
0&0&q&0\\
0&0&0&q^3
\end{smallmatrix}\right)
\left(\begin{smallmatrix}
1&1&1&1\\
0&1&1+q&1+q+q^2\\
0&0&1&1+q+q^2\\
0&0&0&1
\end{smallmatrix}\right)=
\left(\begin{smallmatrix}
1&1&1&1\\
0&1&1+q&1+q+q^2\\
0&0&q&q(1+q+q^2)\\
0&0&0&q^3
\end{smallmatrix}\right),
\end{equation}
\begin{equation}
\label{Psi_2(q)}\Psi_2(q)=\sigma_2^s(q)D_2^s(q)=\left(\begin{smallmatrix}
1          &0 &0\\
-(1+q^{-1})&1 &0\\
q^{-1}     &-1&1
\end{smallmatrix}\right)
\left(\begin{smallmatrix}
q&0&0\\
0&1&0\\
0&0&1
\end{smallmatrix}\right)=
\left(\begin{smallmatrix}
q&0&0\\
-(1+q)&1&0\\
1&-1&1
\end{smallmatrix}\right),
\end{equation}
\begin{equation}
\label{Psi_3(q)} \Psi_3(q)= \left(\begin{smallmatrix}
 1    &0          &0              &0\\
-(1+q^{-1}+q^{-2})    &1         &0                &0\\
q^{-1}(1+q^{-1}+q^{-2})&-(1+q^{-1})&1                 &0\\
-q^{-3}&q^{-1}&-1&1
\end{smallmatrix}\right)\left(\begin{smallmatrix}
q^3&0&0&0\\
0&q&0&0\\
0&0&1&0\\
0&0&0&1
\end{smallmatrix}\right)= \left(\begin{smallmatrix}
q^3     &0&0&0\\
-q(1+q+q^2)&q&0&0\\
(1+q+q^2)&-(1+q)&1&0\\
- 1&1&-1&1
\end{smallmatrix}\right).
\end{equation}
{\it Acknowledgements.} {The second author would like to thank the
Max-Planck-Institute of Mathematics and
the Institute of Applied Mathematics, University of Bonn
for the hospitality. The partial financial support by the DFG
project 436 UKR 113/87 is gratefully acknowledged.}

\end{document}